\documentclass[a4paper,reqno,12pt]{scrartcl}

\usepackage[english]{babel}

\usepackage{amsmath}
\usepackage{cases}
\usepackage{amssymb}
\usepackage{amsthm}
\usepackage{amsfonts}
\usepackage{graphicx}
\usepackage{subcaption}
\usepackage{adjustbox}
\usepackage{enumitem}
\usepackage{color}
\usepackage[text={16.5cm,22.5cm},centering]{geometry}
\usepackage{manfnt}
\usepackage{bm}
\usepackage{pifont}
\usepackage{tikz}
\usepackage{booktabs}
\usepackage{array}
\usepackage{mathrsfs}
\usepackage{calligra}
\usepackage[T1]{fontenc}
\usepackage{thmtools}

\newcommand{\capdot}{\mathbin{\ooalign{\hfil$\cap$\hfil\cr\hfil$\cdot$\hfil\cr}}}

\usetikzlibrary{knots, calc}

\numberwithin{equation}{section}

\setcapindent{0pt} 

\newcommand{\hash}{\mathbin{\#}}

\newcommand{\ssum}[1]{%
  (#1)
}

\newtheorem{Lemma}{Lemma}[section]
\newtheorem{Prop}{Proposition}[section]
\newtheorem{ClaimL}{Claim}[Lemma]

\newtheorem*{claim*}{Claim}
\newtheorem{Cor}{Corollary}
\newtheorem{Question}{Question}

\theoremstyle{definition}
\newtheorem{Def}[Lemma]{Definition} 

\newtheorem{Remark}[Lemma]{Remark}
\newtheorem{Fact}[Lemma]{Fact}

\newcommand{\bigstackarrowww}[1]{%
  \mathrel{
    \begin{tikzpicture}[baseline={(current bounding box.center)}]
      \node (a) at (0,0) {\LARGE$\Rightarrow$};
      \node at (a.north) [yshift=1.2ex] {\normalsize$#1$};
    \end{tikzpicture}
  }
}

\newcommand{\bigstackarrowws}{%
  \mathrel{
    \begin{tikzpicture}
      \node (a) at (0,0) {\LARGE$\Rightarrow$};
    \end{tikzpicture}
  }
}

\newcommand{\Fan}{\operatorname{Fan}}

\newcommand{\cordop}[1]{\mathop{\bm{\circlearrowleft}}\!\left(#1\right)}

\renewcommand{\vert}{\operatorname{vert}}
\newcommand{\ray}{\operatorname{ray}}

\newcommand{\MIN}{\operatorname{MIN}}

\newcommand{\Cross}{\operatorname{Cross}}

\newcommand{\Card}{\operatorname{Card}}

\newcommand{\id}{\operatorname{id}}

\DeclareMathOperator{\dom}{dom}

\DeclareMathOperator{\pr}{pr}

\renewcommand{\mid}{\,:\,}
\newenvironment{enumerate-(a)}{\begin{enumerate}[label={\upshape (\alph*)}, leftmargin=2pc]}{\end{enumerate}}
\newenvironment{enumerate-(a)-r}{\begin{enumerate}[label={\upshape (\alph*)}, leftmargin=2pc,resume]}{\end{enumerate}}
\newenvironment{enumerate-(a)-5}{\begin{enumerate}[label={\upshape (\alph*)}, leftmargin=2pc,start=5]}{\end{enumerate}}
\newenvironment{enumerate-(A)}{\begin{enumerate}[label={\upshape (\Alph*)}, leftmargin=2pc]}{\end{enumerate}}
\newenvironment{enumerate-(A)-r}{\begin{enumerate}[label={\upshape (\Alph*)}, leftmargin=2pc,resume]}{\end{enumerate}}
\newenvironment{enumerate-(i)}{\begin{enumerate}[label={\upshape (\roman*)}, leftmargin=2pc]}{\end{enumerate}}
\newenvironment{enumerate-(i)-r}{\begin{enumerate}[label={\upshape (\roman*)}, leftmargin=2pc,resume]}{\end{enumerate}}
\newenvironment{enumerate-(I)}{\begin{enumerate}[label={\upshape (\Roman*)}, leftmargin=2pc]}{\end{enumerate}}
\newenvironment{enumerate-(I)-r}{\begin{enumerate}[label={\upshape (\Roman*)}, leftmargin=2pc,resume]}{\end{enumerate}}
\newenvironment{enumerate-(1)}{\begin{enumerate}[label={\upshape (\arabic*)}, leftmargin=2pc]}{\end{enumerate}}
\newenvironment{enumerate-(1)-r}{\begin{enumerate}[label={\upshape (\arabic*)}, leftmargin=2pc,resume]}{\end{enumerate}}

\newenvironment{enumerate-(star)}{\begin{enumerate}[label={\upshape{(\( \star_{ \arabic*} \))}}, leftmargin=2pc]}{\end{enumerate}}

\renewcommand{\a}{\alpha}

\newcommand{\f}{\varphi}

\newcommand{\R}{\mathbb{R}} 
 
\newcommand{\N}{\mathbb{N}} 

\newcommand{\DD}{\mathscr{D}}

\newcommand{\MM}{\mathcal{M}}

\newcommand{\RR}{\mathcal{R}}
\newcommand{\QQ}{\mathcal{Q}}

\newcommand{\OO}{\mathcal{O}} 
\newcommand{\mSS}{\mathscr{S}} 
 
\newcommand{\XX}{\mathcal{X}} 
\newcommand{\KK}{\mathscr{K}} 
\newcommand{\JJ}{\mathscr{J}}

\newcommand{\cS}{\mathcal{S}} 
\newcommand{\ocirc}[1]{\overset{\circ}{#1}{}}

\renewcommand{\int}{\operatorname{int}}

\allowdisplaybreaks

\usepackage[colorlinks=true, linkcolor=blue, citecolor=blue, urlcolor=blue]{hyperref}

\mathchardef\mhyphen="2D
\newcommand{\rest}{\!\restriction\!}

\newcommand{\St}{\operatorname{St}}

\newcommand{\es}{\varnothing}
\renewcommand{\ge}{\geqslant}
\renewcommand{\le}{\leqslant}

\newcommand{\e}{\varepsilon}
\renewcommand{\d}{\delta}

\renewcommand{\Cup}{\bigcup}
\newcommand{\la}{\langle}
\newcommand{\ra}{\rangle}
\newcommand{\cat}{{}^\frown}

\newcommand{\intersect}{\operatorname{\mathbf{i}}}

\usepackage[labelfont=bf]{caption}

\newcommand{\oD}{\overset{\circ}{\square}{}}

\usepackage{stmaryrd} 
\newcommand{\dint}[1]{\llbracket #1\rrbracket}
\newcommand{\openint}[1]{\left] #1 \right[}
\setkomafont{title}{\normalfont\Large\bfseries}
\title{Additivity of Crossing Number via Restricted Reidemeister Moves} 

\author{Vadim Weinstein}

\begin{document}

\maketitle

\begin{abstract}
  We define a set of restricted Reidemeister moves and show that if
  $K$ is obtained from $K_0\hash K_1$ using those moves, then the
  crossing number of $K$ is at least $c(K_0)+c(K_1)$.  We also explore
  topological interpretations of this result.
\end{abstract}

\tableofcontents




\section{Introduction}
\label{sec:Intro}

It is a long-standing open question whether the minimal crossing
number of a composite knot is equal to the sum of the minimal crossing
numbers of its factors.  The best result for general knots so far has
been due to \cite{Lackenby2009} showing that for any finite number of
knots $K_1,\dots,K_n$ the inequality
$c(K_1)+\cdots+c(K_n)\le 152c(K_1\hash \cdots\hash K_n)$ holds.  It
has also been established that for a certain class of links,
$c(L_1\hash L_2)\ge c(L_1)+3$ \cite{Diao2004}. Equality
$c(K_1\hash K_2)=c(K_1)+c(K_2)$ has been previously established when
$K_1$ and $K_2$ are both alternating knots \cite{Kauffman1987} or both
torus knots~\cite{Diao2004}.

Given a knot diagram $K$, denote by $\chi(K)$ the number of crossings
in that diagram. Denote by $\RR$ the set of classical Reidemeister
moves and Reidemeister equivalence between diagrams as
$K\sim_{\RR}K'$.  Then $c(K)$ is defined as the minimim of
$\{\chi(K')\mid K'\sim_\RR K\}$.  Now the question of additivity can
be reformulated as follows.  Is it possible that
$K'\sim_\RR K_0\hash K_1$ and $\chi(K')<c(K_0)+c(K_1)$?  We will
define a set of restricted Reidemeister moves $\RR'$ (in fact several
versions) and show that if $K'\sim_{\RR'}K$, then
$\chi(K')\ge c(K_0)+c(K_1)$.

A key part of the proof is bounding the number of certain type of crossings
(which we call \emph{hybrid} crossings)
using a graph-theoretic lemma in Section \ref{ssec:Matrices}.

\section*{Acknowledgements}
The author was generously supported by a European Research Council
Advanced Grant (ERC AdG, ILLUSIVE: Foundations of Perception
Engineering, 101020977).  The author is grateful to Steven LaValle for
encouragement and support, Kalle Timperi for early discussions, and
Tapani Hyttinen for comments on the proof. Special thanks are due to
Marc Lackenby for insightful suggestions and comments on the
manuscript.

\section{Preliminaries}
\label{sec:Preli}

In this section we will define basic notation used in the paper
and the fundamental concepts of genericity and general position,
knot and link diagrams, centered and partial diagrams, restricted Reidemeister
moves and equivalence, and crossing number.

\subsection{Basics}
\label{ssec:Basics}

If $a,b\in\R^n$, denote by $[a,b]$ the line segment connecting them.
For $a,b\in\R$, denote the open interval by
$\openint{a,b}=\{t\mid a<t<b\}$. In $\R$ these sets are empty, if
$b<a$.  Denote by $I=[0,1]\subset\R$ the unit interval and by
$\square=I\times I$ the unit square, $\oD$ is its interior and
$\partial \square$ its boundary.  Let $o^-,o^+,o\in \square$ be the
points $o^-=(\tfrac{1}{2},0)$, $o^+=(\tfrac{1}{2},1)$, and
$o=(\tfrac{1}{2},\tfrac{1}{2})$. The point $o$ is called the
\emph{critical point}.  Denote by $I^+=[\tfrac{1}{2},1]$ and
$I^-=[0,\tfrac{1}{2}]$ and $\square^{\pm}=I\times I^{\pm}$.  We work
in the piecewise linear category and all maps are assumed to be
piecewise linear with respect to the linear structure of $\R^n$ unless
stated otherwise. A \emph{generic immersion} of a $1$-manifold $D$
into $\R^2$ is a map $f\colon D\to \R^2$ which is locally an embedding
(every point of the domain has a neighborhood restricted to which the
map is an embedding), $f\rest\partial D$ is an embedding, and $f$ has
only finitely many multiplicity points each of which is a transversal
double point. The set of these double points is denoted $\Cross(f)$.
The cardinality of this set is denoted~$\chi(f)$.

By $\dom(f)$ denote the domain of a function $f$ and by $|f|$ its
range.  Often, however, we identify a function $f$ with its range
$|f|$ and might denote, for example $f\cap g\subset A$ instead of
$|f|\cap |g|\subset A$, or $f\rest C\subset g$ instead of
$|f\rest C|\subset |g|$. It should always be clear from the context
when this is the case. If $f$ and $g$ are maps which agree on
$\dom(f)\cap \dom(g)$, denote by $f\cup g$ the map with domain
$\dom(f)\cup\dom(g)$ such that $(f\cup g)(x)=f(x)$ for $x\in\dom(f)$
and $(f\cup g)(x)=g(x)$ for $x\in\dom(g)$. Sometimes, by abuse of
notation, we write the concatenation of curves using $\cup$ instead of
$\cat$, when the intended meaning is clear from the context.
The notation $C=A\sqcup B$ means that $C=A\cup B$ and $A\cap B=\es$.

The composition of functions is denoted $f\circ g$ but sometimes
``$\circ$'' is ignored and it is denoted~$fg$.  If $A$ is a set,
denote by $fA=f[A]=f(A)=|f\rest A|=\{f(x)\mid x\in A\}$ the image
of~$A$ (usually $f[A]$ is preferred).  By $\id_A$ we denote the
identity map $A\to A$. The identity map \emph{fixes $A$
  pointwise}. Sometimes we say that a function $f$ fixes a set $A$
pointwise to mean that $f\rest A=\id_A$.

By $\ocirc{U}$, $\partial U$, and $\bar U$ denote the interior,
closure, and the boundary of a set~$U\subset\R^n$ respectively.
By $\Card(A)$ denote the cardinality of the set $A$.  
If $\gamma_1$ and $\gamma_2$ are mappings from $I$ into some space $X$
such that $\gamma_1(1)=\gamma_2(0)$, we denote by
$\gamma_1\cat\gamma_2$ the map $\gamma\colon I\to \R^2$ given by
$\gamma(t)=\gamma_1(2t)$ if $t<1/2$, and $\gamma(t)=\gamma_2(2t-1)$
otherwise.  By $B(x,r)$ and $\ocirc{B}(x,r)$ denote the closed and
open balls respectively centered at $x$ with radius~$r$.  

\subsection{Genericity}

\begin{Def}\label{def:BasicImmersion}
  A \emph{polyhedron} in $\R^n$ is a set which is the realization of a
  finite simplicial complex in $\R^n$ compatible with the linear
  structure of~$\R^n$. A map $f\colon X\to Y$ between polyhedra $X$
  and $Y$ is an \emph{immersion} if for all $x\in X$ there is a
  polyhedral neighborhood $U$ of $x$ (i.e. $U$ is a polyhedron and
  $x\in \ocirc{U}$) such that $f\rest U$ is an embedding. A
  point $x\in |f|$ is a \emph{multiplicity point} if
  $\Card(f^{-1}(x))>1$. It is a \emph{double point} if
  $\Card(f^{-1}(x))=2$.
  Suppose $T\colon D\to \R^2$ is an immersion
  where $D$ is a $1$-manifold with boundary. 
  Suppose that $x\in\R^2$ is a double point,
  Then $x$ is
  \emph{transversal} if there is a compact polyhedral
  neighborhood $U$ of $x$, a homeomorphism $h\colon U\to [-1,1]^2$,
  and disjoint closed connected sets $C_1,C_2\subset D$ each of which
  contains one of the points of $T^{-1}(x)$ in their interior such
  that $|h\circ (T\rest C_1)|=[-1,1]\times \{0\}$ and
  $|h\circ (T\rest C_2)|=\{0\}\times [-1,1]$.
\end{Def}

\begin{Def}\label{def:generic}
  Let $D$ a $1$-manifold with boundary.  We say that an immersion $g$
  of $D$ into $\R^2$ is \emph{generic} if there are at most finitely
  many multiplicity points, all of them are transversal double points
  and none of them is at $g[\partial D]$. These double points are
  called \emph{crossings} and the set of crossings is denoted
  $$\Cross(g)=\{x\in\R^2\mid x\text{ is a double point of }g\}.$$
\end{Def}

\subsection{Centered Knot Diagrams}
\label{ssec:KD}

\begin{figure}
  \centering
  
  \begin{subfigure}[t]{0.31\textwidth}
    \centering
    \fbox{\includegraphics[width=0.8\linewidth]{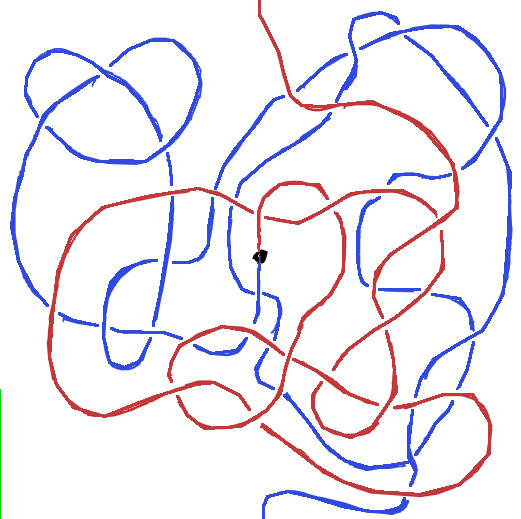}}
    \caption{Another knot t-diagram $K$.}
  \end{subfigure}
  \hfill
  \centering
    \begin{subfigure}[t]{0.31\textwidth}
    \centering
    \fbox{\includegraphics[width=0.8\linewidth]{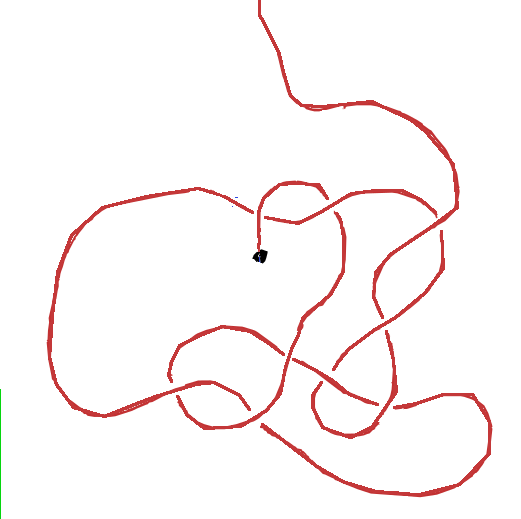}}
    \caption{The upper partial knot diagram $K^{+}$ of the knot in (a).}
  \end{subfigure}
  \hfill
    \centering
    \begin{subfigure}[t]{0.31\textwidth}
    \centering
    \fbox{\includegraphics[width=0.8\linewidth]{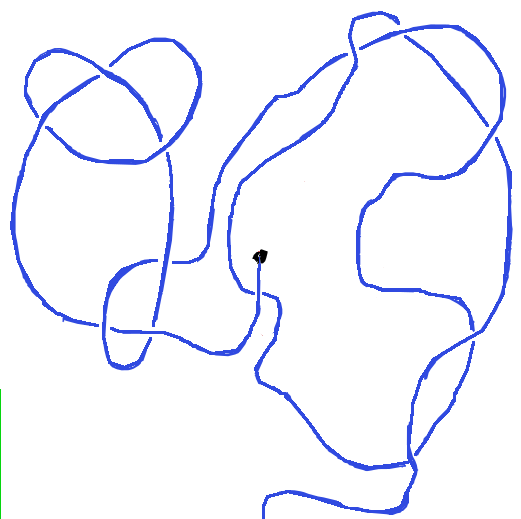}}
    \caption{The lower partial knot diagram $K^{-}$ of the knot in (a)}
  \end{subfigure}
  \caption{(a): A centered knot diagram $K$ and the corresponding
    upper (b) and lower (c) partial knot diagrams $K^+$ and $K^-$.
    Conversely, the diagram }
  \label{fig:CenteredIntro}
\end{figure}
A \emph{knot diagram (k.d.)} is a pair $(K,\ell)$ such that
$K\colon I\to \square$ is a generic immersion with
$K(0)=o^-,K(1)=o^+$, $K(t)\in \oD$ for all $0<t<1$, and such that the
one-sided derivatives of $K$ at $t=0$ and $t=1$ are vertical, and
$\ell\colon\Cross(K)\to I^2$ is such that $\ell(x)=(a,b)$ where
$\{a,b\}= K^{-1}(x)$ for all $x$ (it is an ordering of the inverse
image of $x$ under~$K$). The diagram is \emph{centered}, if also
$K(\tfrac{1}{2})=o$, $o\notin\Cross(K)$, and the derivative of $K$ at
$t=\tfrac{1}{2}$ exists and is vertical.  Note that for a piecewise
linear map the derivatives exist for all but finitely many points and
the one-sided derivatives at the end-points always exist. Denote
\begin{equation}
  K^{\pm}=K\rest I^{\pm}.\label{eq:Kpm}
\end{equation}
Conventionally in illustrations, $K^+$ will
be drawn with \emph{\textcolor{red}{red}} color and $K^-$ with
\emph{\textcolor{blue}{blue}}. See Figure~\ref{fig:CenteredIntro}.
More generally, let us define partial knot diagrams:

\subsection{Partial Knot Diagrams}
\label{ssec:Partial}

An \emph{upper (lower) partial knot diagram} is a pair $(P,\ell)$
(resp. $(Q,\ell)$) where $P$ (resp. $Q$) is a generic immersion
$P\colon I^{+} \to \square$ (resp. $Q\colon I^-\to\square$) such that
$P(1)=o^+$ (resp. $Q(0)=o^-$) and the one-sided derivative at $o^+$
(resp. $o^-$) is vertical. The function
$\ell\colon \Cross(P)\to (I^+)^2$ (resp.
$\ell\colon \Cross(Q)\to (I^-)^2$) is a crossing labeling function
just like for knot diagrams. The partial diagram is \emph{centered},
if also $P(\tfrac{1}{2})=o$ (resp. $Q(\tfrac{1}{2})=o$) and the
one-sided derivative at $\tfrac{1}{2}$ is vertical.  Clearly, if $K$
is a centered knot diagram, then $K^{\pm}$ are upper and lower
centered partial knot diagrams respectively.

\begin{Remark}
  The concept of partial knot diagrams is reminiscient of knotoids
  introduced in~\cite{Turaev2012Knotoids}. The author wishes to thank
  Prof. Marc Lackenby for pointing this out.
\end{Remark}

\subsection{Restricted Reidemeister Moves}
\label{ssec:RReid}

\begin{figure}
  \centering
  \includegraphics[width=\textwidth]{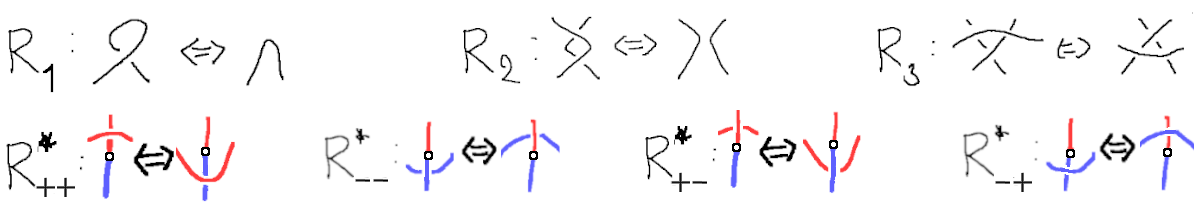}
  \caption{Reidemeister moves and sliding moves}
  \label{fig:RReid}
\end{figure}

The set of (classical) Reidemeister moves is denoted $\RR$.  It
contains the moves $R_0,R_1,R_2,R_3$. The moves $R_i$ for
$i\in \{1,2,3\}$ are depicted in the top half of Figure \ref{fig:RReid}.  The
move $R_0$ is an application of a self-homeomorphism of $\square$
which is the identity on a neighbourhood of
$\partial\square\cup \{o\}$. The equivalence relation ``the knot diagram
$K$ is obtained from the knot diagram $K'$ using Reidemeister moves''
is denoted $K\sim_{\RR}K'$. The corresponding \emph{restricted
  Reidemeister moves} $R^*_i$ are such that $R^*_0=R_0$, and for
$i\in \{1,2,3\}$ the move $R^*_i$ is like $R_i$ but has to be
performed in a neighborhood which does not contain the critical
point~$o$.  Let $\RR^*=\{R^*_0,R^*_1,R^*_2,R^*_3\}$ be the set of
\emph{non-sliding restricted Reidemeister moves}.  The equivalence
relation corresponding to $\RR^*$ is denoted~$\sim^*_{\RR}$.

The \emph{sliding moves} $R^*_{++}$, $R^*_{--}$, $R^*_{+-}$, and
$R^*_{-+}$ are moves in which an arc slides under or over the critical
point (Figure~\ref{fig:RReid} below). For $\pm\in \{-+\}$, the move
$R^*_{\pm -}$ allows a subarc of $K^{\pm}$ to slide \emph{under} the critical
point, and the move $R^*_{\pm +}$ allows a subarc of $K^{\pm}$ to
slide over the critical point.

We define four (extended) sets of restricted Reidemeister moves:
\begin{align*}
  \RR^*_{opp1}&=\RR^*\cup \{R^*_{++},R^*_{--}\}&&
  \RR^*_{same1}=\RR^*\cup \{R^*_{+-},R^*_{--}\}\\
  \RR^*_{opp2}&=\RR^*\cup \{R^*_{+-},R^*_{-+}\}&&
  \RR^*_{same2}=\RR^*\cup \{R^*_{++},R^*_{-+}\}\\
\end{align*}
In the sets $\RR^*_{opp1}$ and $\RR^*_{opp2}$, the red and blue arcs
(corresponding to $K^+$ and $K^-$ respectively) are allowed to slide on
the \emph{opposite} sides of the critical point. In $\RR^*_{opp1}$ the red
arc slides \emph{over} and the blue arc \emph{under}, and in $\RR^*_{opp2}$
vice versa. In $\RR^*_{same1}$ and $\RR^*_{same2}$ they slide on the \emph{same}
side of the critical point. In $\RR^*_{same1}$ both slide \emph{under} and in $\RR^*_{same2}$
both slide \emph{over}.
The corresponding equivalence relations are respectively denoted
$$\sim^{o1}_{\RR}, \quad \sim^{o2}_{\RR}, \quad \sim^{s1}_{\RR},\quad \text{and}\quad\sim^{s2}_{\RR}.$$

\begin{Remark}\label{rem:EquivImpl}
  For $X\in \{o1,o2,s1,s2\}$ and any knot diagrams $K$ and $K'$ we have
  $$K\sim^*_\RR K'\Rightarrow K\sim^X_{\RR}K'\Rightarrow K\sim_{\RR}K'.$$
  A \emph{mirror image} of a knot diagram $K=(K,\ell)$ is $\mu(K)=(K,\ell_\mu)$
  where $\ell_\mu$ flips every label given by~$\ell$, i.e. if $\ell(x)=(a,b)$, then
  $\ell_{\mu}=(b,a)$. It is now clear that for all knot diagrams $K$ and $K'$,
  $$K\sim^{o1}_\RR K'\iff \mu(K)\sim^{o2}_{\RR} \mu(K')\qquad \text{ and }\qquad K\sim^{s1}_\RR K'\iff \mu(K)\sim^{s2}_{\RR} \mu(K').$$
\end{Remark}

\subsection{Composite Knot Diagrams}
\label{ssec:CompositeKD}

Let $H=I\times \left\{\tfrac{1}{2}\right\}$ be the horizontal line
passing through the middle of~$\square$. 
Suppose that $K$ is a centered knot diagram. If $K\cap H=\{o\}$,
we say that $K$ is a \emph{composite knot diagram}. For $\pm\in\{-,+\}$,
let $h^{\pm}\colon I^{\pm}\to I$ be the (unique) linear orientation preserving
homeomorphism. Let $\widehat h^{\pm}\colon \square^{\pm}\to \square$ be given
by $\widehat h^{\pm}(x,y)=(x,h^{\pm}(y))$. Further, let $h^\pm_*\colon I^2\to I^2$
be given by $h^{\pm}_*(a,b)=(h^{\pm}(a),h^{\pm}(b))$.
If $K=(K,\ell)$ is a composite k.d., define $K_0=(K_0,\ell_0)$ and $K_1=(K_1,\ell_1)$
be the \emph{factors} of $K$ defined by
\begin{align*}
  K_0&=\widehat h^{-}\circ K\circ (h^{-})^{-1}&&\ell_0=h^-_*\circ\ell\circ (\widehat h^-)^{-1}\\
  K_1&=\widehat h^{+}\circ K\circ (h^{+})^{-1}&&\ell_1=h^+_*\circ\ell\circ (\widehat h^+)^{-1}.
\end{align*}
If $K_0$ and $K_1$ are factors of $K$, denote $K=K_0\hash K_1$.
\begin{Remark}\label{rem:MirrorComp}
  $K=K_0\hash K_1$ if and only if $\mu(K)=\mu(K_0)\hash \mu(K_1)$.
\end{Remark}

\subsection{Crossing Number}
\label{ssec:CrN}

The number of double points in a (partial) knot diagram $K$ is denoted $\chi(K)$. As usual, define
\begin{equation}
  c(K)=\min\{\chi(K')\mid K'\sim_\RR K\}.\label{eq:MinCrossingNumber}
\end{equation}

For (partial) knot diagrams $P$ and $Q$ denote by $\intersect(P,Q)$ the
number of crossings between $P$ and $Q$ not counting $o$ or the boundary:
$$\intersect(P,Q)=\Card(P\cap Q\setminus (\partial\square\cup \{o\})).$$
Here we assume that $P$ and $Q$ are \emph{mutually in general position
  (modulo $o$ and $\partial\square$)}, i.e.
$\Cross(P)\cap\Cross(Q)=\es$ and all intersections between $P$ and $Q$
(except for $o$ and those on the boundary) are transversal double
points.

\begin{Remark}\label{rem:mirror}
  For all knot diagrams $K$, $c(K)=c(\mu(K))$ and
  $\chi(K)=\chi(K^+)+\chi(K^-)+\intersect(K^+,K^-)$.
\end{Remark}

\section{Results}
\label{sec:Results}

In this section we present the main results of this paper.

\subsection{The Main Result}

Recall the definitions of the equivalence relaions $\sim^X_\RR$ for
$X\in \{*,o1,o2,s1,s2\}$ from Section~\ref{ssec:RReid}.

\begin{restatable}{Thm}{MainThm}\label{thm:Main}
  Suppose that $K=K_0\hash K_1$ and $X\in \{*,o1,o2,s1,s2\}$.
  If $\widehat K\sim^X_\RR K$, then
  $\chi(\widehat K)\ge c(K_0)+c(K_1)$.
\end{restatable}

By Remarks \ref{rem:EquivImpl}, \ref{rem:MirrorComp}, and
\ref{rem:mirror}, the case $X=*$ follows from any other case, and the
cases $X=o2$ and $X=s2$ follow by symmetry from $X=o1$ and $X=s1$
respectively. So it will be enough to consider the cases $X=o1$
and~$X=s1$.

\subsection{Topological interpretations}
\label{ssec:TopInt}

The case $X=s1$ has an intuitive topological interpretation. Let
$C=\square\times I=I^3$ be the unit cube. Let
$\hat o^{\pm}=(o^{\pm},\tfrac{1}{2})$ be the middpoints of the upper
and lower faces respectively, and let $\hat o=(o,\tfrac{1}{2})$ be the
midpoint of the cube. Let $f^+=(\tfrac{1}{2},\tfrac{1}{2},1)$ be the
midpoint of the frontal face of the cube. Let
$W=I\times\left\{\tfrac{1}{2}\right\}\times I$ be the horizontal
middle cross-section of the cube.  A \emph{knot} is an embedding
$\KK\colon I\to C$ such that $\KK(0)=\hat o^-$, $\KK(1)=\hat o^+$,
$\KK(t)\in\ocirc{C}$ for $0<t<1$, and the one-sided derivatives at
$t=0$ and $t=1$ are vertical.  It is \emph{centered}, if
$\KK\left(\tfrac{1}{2}\right)=\hat o$ and the derivative at
$\tfrac{1}{2}$ is vertical. It is \emph{composite}, if
$\KK\cap W=\{\hat o\}$. It is in \emph{general position}, if
$\pr\circ\KK$ is a generic immersion where $\pr\colon C\to\square$ is
the projection along the third coordinate $\pr(x,y,z)=(x,y)$.  For a
centered knot $\KK$, denote $\KK^{\pm}=\KK\rest I^{\pm}$.  Let
$C^\pm=\square^\pm\times I$ be the upper and lower halves of the cube.
Let $g^{\pm}\colon C^{\pm}\to C$ be the linear homeomorphism given by
$g(x,y,z)=(x,h(y),z)$ where $h$ is as defined in
Section~\ref{ssec:CompositeKD}. If $\KK$ is composite, let
$\KK_0=g^-\circ\KK^-$ and $\KK_1=g^+\circ \KK^+$.  If the knot $\KK$
is in general position, define $\DD(\KK)=(K,\ell)$ be the knot diagram
such that $K=\pr\circ\KK$ and $\ell(x)$ orders the set $K^{-1}(x)$
according to the distance from $\square\times \{0\}$.  It is now clear
that for a composite $\KK$, $\DD(\KK)$ is a composite knot diagram,
and $\DD(\KK_0)=\DD(\KK)_0$ and $\DD(\KK_1)=\DD(\KK)_1$. We can now
formulate a corollary to Theorem~\ref{thm:Main}:

\begin{Cor}\label{cor:Topo1}
  Suppose $\KK$ is a centered composite knot in general position.  Let
  $K_0=\DD(\KK_0)$ and $K_1=\DD(\KK_1)$.  Suppose that
  $h\colon C\to C$ is any self-homeomorphism of $C$ which fixes
  pointwise the boundary $\partial C$ and the line segment
  $[\hat o,f^+]$. Let $K'=\DD(h\KK)$.  Then $\chi(K')\ge c(K_0)+c(K_1)$.
\end{Cor}
\begin{proof}
  It is easy to see that $\DD(h\KK)$ can be obtained from
  $\DD(\KK)=K_0\hash K_1$ using moves in $\RR^*_{same1}$.
  So $K'\sim^{s1}_{\RR} K_0\hash K_1$. By Theorem~\ref{thm:Main},
  $\chi(K')\ge c(K_0)+c(K_1)$.
\end{proof}

The cases $X=o1,o2$ also have a topological interpretation, but it is
less intuitive. Let $f^-=(\tfrac{1}{2},\tfrac{1}{2},0)$ be the midpoint of the
back face of the cube (opposite of $f^+$) and let $W_d$ be the diagonal cross section
of the cube which equals to the convex hull of the set
$\{(0,0,1),(1,0,1),(0,1,0),(1,1,0)\}$. We then have:

\begin{Cor}\label{cor:Topo2}
  Suppose $\KK$ is a centered composite knot in general position such
  that $\KK\cap W_d=\{\hat o\}$. Let $K_0=\DD(\KK_0)$ and
  $K_1=\DD(\KK_1)$.  Suppose that $h\colon C\to C$ is any
  self-homeomorphism of $C$ such that
  $h[W_d]\cap [f^-,f^+]=\{\hat o\}$. Let $K'=\DD(h\KK)$.  Then
  $\chi(K')\ge c(K_0)+c(K_1)$.
\end{Cor}
\begin{proof}
  It is not hard to see that $\DD(h\KK)$ can be obtained from
  $\DD(\KK)$ using moves in $\RR^*_{opp1}$.
  So $K'\sim^{o1}_{\RR} K_0\hash K_1$. By Theorem~\ref{thm:Main},
  $\chi(K')\ge c(K_0)+c(K_1)$.
\end{proof}

\begin{Remark}
  The result of Corollary \ref{cor:Topo2} is \emph{ad hoc}. Note, for example,
  that the set of homeomorphisms which satisfy the premise of the corollary
  is not closed under composition. See Question~\ref{q:1}
\end{Remark}

\subsection{Composition of Links}

A version of Theorem \ref{thm:Main} also holds for links in the
following way.  A \emph{link diagram} is a triple $L=(L,M,\ell)$ where
$L\colon D\to \square$ is a generic immersion, $M\in\N$,
$D=I\cup (M\times S^1)$ is the disjoint union of $I$ with $M$ many
circles, $L\rest I=(L\rest I,\ell\rest\Cross(L\rest I))$ is a knot
diagram, and for all $t\in M\times S^1$, $L(t)\in \oD$.  It is
\emph{centered}, if $L\rest I$ is centered. It is \emph{composite}, if
$L\cap H=\{o\}$. Note that it follows that $(L\rest I)\cap H=\{o\}$ is
composite (recall the definition of $H$ in the beginning of
Section~\ref{ssec:CompositeKD}).  A \emph{partitioned link} is a tuple
$(L,M^+,M^-,\ell)$ such that $(L,M^+\cup M^-,\ell)$ is a link and
$M=M^+\sqcup M^-$.  For a composite $L$, there is a unique partition
$M=M^+\sqcup M^-$ such that
\begin{equation}
  L\rest (I^{\pm}\cup (M^{\pm}\times S^1))\subset \square^{\pm}\label{eq:partitionLink}
\end{equation}
Given a partitioned link $L$, let $L^\pm=(L\rest (I^{\pm}\cup (M^{\pm}\times S^1)))$.

Let $L=(L,M,\ell)$ be a composite link. Let $\bar L=(L,M^+,M^-,\ell)$ be the
associated partitioned link with the partition witnessing
\eqref{eq:partitionLink}, and let the homeomorphisms
$$h^{\pm}_L\colon I^{\pm}\cup (M^{\pm}\times S 1)\to I\cup (M\times S^1)$$
be defined by $h^{\pm}_L(t)=h^{\pm}(t)$, if $t\in I^{\pm}$ (where
$h^{\pm}$ is as defined in Section~\ref{ssec:CompositeKD}) and
otherwise $h^{\pm}_L(m,t)=(m,t)$, i.e.
$h^{\pm}_L\rest(M^{\pm}\times S^1)$ is the standard inclusion into
$M\times S^1$. 
Define the \emph{factors} $L_0=(L_0,M_0\ell_0)$ and
$L_1=(L_1,M_1,\ell_1)$ of $L$ as
\begin{align*}
  L_0&=\widehat h^{-}\circ L\circ (h^{-}_L)^{-1}&&\ell_0=h^-_*\circ\ell\circ (\widehat h^-)^{-1}&&M_0=M^-\\
  L_1&=\widehat h^{+}\circ K\circ (h^{+}_L)^{-1}&&\ell_1=h^+_*\circ\ell\circ (\widehat h^+)^{-1}&&M_1=M^+.
\end{align*}
where $\widehat h^{\pm}$ and $h^{\pm}_*$ are as in
Section~\ref{ssec:CompositeKD}.
If $L_0$ and $L_1$ are factors of $L$, denote $L=L_0\hash L_1$.
The Reidemeister moves, crossing
number, and minimal crossing number for links are defined just as for
knots. Given a partiioned link $L$, define the restricted Reidemeister
moves $R^{^*}_{\pm\mp}$ just as for knots replacing $K^{\pm}$ with $L^{\pm}$.
For example in $R^*_{++}$ for links, the subarc of $L^+$ can only pass \emph{over}
the critical point~$o$. This enables to generalize the relations
$\sim^X_\RR$, $X\in \{*,o1,o2,s1,s2\}$ for links. We can now state a
version of Theorem~\ref{thm:Main} for links:

\begin{restatable}{Thm}{MainThmLinks}\label{thm:MainLinks}
  Suppose that $L=L_0\hash L_1$ and $X\in \{*,o1,o2,s1,s2\}$.
  If $L'\sim^X_\RR \bar L$, then $\chi(L')\ge c(L_0)+c(L_1)$.
\end{restatable}

We will dedicate this paper to proving Theorem~\ref{thm:Main} and explain
in Section~\ref{sec:Links} how to adapt the proof to links.

\begin{Cor}[Topological interpretation for links]
  Corollaries \ref{cor:Topo1} and \ref{cor:Topo2} hold for
  link compositions instead of knot compositions.
\end{Cor}

\subsection{A Version for Knots Diagrams on $S^2$}

\begin{figure}
  \centering
  \begin{subfigure}{0.3\textwidth}    
    \centering
    \fbox{\includegraphics[width=\linewidth]{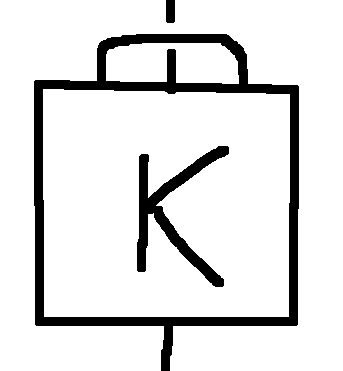}}
  \end{subfigure}
  $\bigstackarrowws$
  \begin{subfigure}{0.3\textwidth}
    \centering
    \fbox{\includegraphics[width=\linewidth]{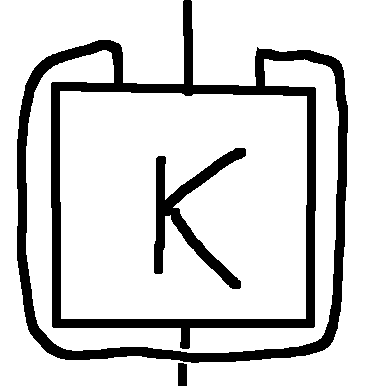}}
  \end{subfigure}
  \caption{Crossing the infinity point on $S^2$ translated to a tangle diagram.}
  \label{fig:NonSliding}
\end{figure}

The results and proofs in this paper are stated for knots presented as
tangle diagrams, but they apply for standard knots through a standard
transformation. 
Parametrize the unit circle as
$S^1=\{e^{i\theta}\mid 0\le \theta\le 2\pi\}$. Denoet the upper and lower halves of the
circle by $S^1_+=\{e^{i\theta}\mid 0\le \theta\le \pi\}$ and
$S^1_{-}=\{e^{i\theta}\mid \pi\le\theta\le 2\pi\}$.
Think of $S^2=\R^2\cup \{\infty\}$.
Let
$\R_{\ge 0}=\{x\in\R\mid x\ge 0\}$ and
$\R_{\le 0}=\{x\in\R\mid x\le 0\}$.  A \emph{(classical) knot diagram} is a generic embedding
$K\colon S^1\to S^2$ equipped with the under and over crossing information.
W.l.o.g. assume that $\infty\notin |K|$.
The notions $c(K)$ and $\chi(K)$ are defined in the same way as previously.
Suppose that $K$ satisfies the following:
$$K(t)\in \R\times \{0\}\iff t\in \{e^{i0},e^{i\pi}\},\quad K(e^{i0})=(0,0),\quad\text{ and }\quad\ K(e^{i\pi})=(1,0).$$
Then the diagram is called \emph{composite}. In this case, let
\begin{align*}
  K_0&=K\cap (R\times \R_{\le 0})\cup ([0,1]\times \{0\}),\ \text{and}\\
  K_1&=K\cap (R\times \R_{\ge 0})\cup ([0,1]\times \{0\}).
\end{align*}
After reparametrization (cf. Section~\ref{ssec:Reparam}) these are
\emph{factors} of~$K$, and we denote $K=K_0\hash K_1$. Define
restricted Reidemeister moves for such diagrams in the same way as in
Section~\ref{ssec:RReid} with the difference that instead of one
critical point $o$ we have two critical points, $(0,0)$ and $(1,0)$.
To be more formal, let $\RR^c$ be the set of those (classical)
Reidemeister moves which are performed in a neighborhood which does
not contain either of the critical points.  A \emph{sliding move} is
performed in a neighborhood which contains exactly one of the critical
points and in such a move an arc slides over or under that point.  Let
$R^c_{\pm \mp}$ be the sliding move where a subarc of
$K\rest S^1_{\pm}$ is only allowed to cross \emph{over} (if $\mp=+$)
or \emph{under} (if $\mp=-$) a critical point. Then define
$\RR^c_{opp1}$, $\RR^c_{opp2}$, $\RR^c_{same1}$, and $\RR^c_{same2}$
by analogy with $\RR^*_{opp1}$, $\RR^*_{opp2}$ etc. This induces the
relations $\sim^{o1}_{\RR}$, $\sim^{o2}_{\RR}$, $\sim^{s1}_{\RR}$, and
$\sim^{s2}_{\RR}$ for classical knot diagrams.  We can then formulate:
\begin{Cor}
  Suppose $K=K_0\hash K_1$ is a composite (classical) knot diagram on
  $S^2$ and $X\in \{*,o1,o2,s1,s2\}$.  If $\widehat K\sim^{X}_{\RR}K$,
  then $\chi(\widehat K)\ge c(K_0)+c(K_1)$.
\end{Cor}
\begin{proof}
  By a standard argument one can assume that the point $(1,0)$ is
  never crossed, i.e. all sliding moves over/under that point can be
  replaced by a sequence of other moves. Instead of crossing $(1,0)$
  one can move the arc over the other side of $S^2$, passing
  through~$\infty$. W.l.o.g. one can assume that a neighborhood of
  $\infty$ is never occupied by $K$ (before or after a restricted
  Reidemeister move). Cutting open at that neighborhood and projecting
  onto $\square$, one obtains knot diagram as in the previous sections
  (as well as the rest of the paper). Crossing the ``other'' critical
  point becomes a move as shown on Figure~\ref{fig:NonSliding}.
  This translates the problem back to the one presented above
  and so we can apply Theorem~\ref{thm:Main}.
\end{proof}

\section{Basics of Tangle Diagrams}
\label{sec:BasicsOfTangles}

We will need to deal with tangle diagrams which are generalizations of
knot and link diagrams defined in Section~\ref{sec:Results}.

\subsection{Definition of Tangle Diagrams}
\label{sec:DefT}

A \emph{tangle diagram} (a \emph{t-diagram} for short) is a tuple
$T=(T,N,M,\ell)$ such that
\begin{enumerate}[label={\upshape ($T_{\arabic*}$)}, leftmargin=3pc]
\item $N$ and $M$ are disjoint finite sets. 
\item $T$ is a generic immersion $T\colon D_T\to \square$ where
  $D_T=(I\times N)\cup (S^1\times M)$ is the disjoint
  union of $\Card(N)$ many unit intervals and $\Card(M)$
  many circles. 
\item \label{def:Tangle3} $|T|$ is disjoint from the corners of the square,
  $|T|\cap \{(0,0),(0,1),(1,0),(1,1)\}=\es$.
\item \label{def:Tangle4} $T^{-1}[\partial \square]\subset \{0,1\}\times N$, i.e. a point is
  mapped into the boundary $\partial\square$ only if it is an endpoint of one of
  the strands, or equivalently a boundary point of the manifold $D_T$.
\item \label{def:Tangle5} $\ell\colon \Cross(T)\to D_T^2$ is a labeling function which tells
  which arc is above the other in a crossing.  For each
  $x\in\Cross(T)$, we have $\ell(x)=(a,b)$ where $a$ and $b$ are such that
  $\{a,b\}= T^{-1}(x)$
\end{enumerate}

Note that it follows from the definition of transversality that
$\Cross(T)\cap \partial \square=\es$. As for knot and link diagrams,
the number of crossings in $T$ is denoted $\chi(T)=\Card(\Cross(T))$.
The restrictions $T\rest (I\times \{i\})$ for $i\in N$ are the
\emph{strands} of~$T$ and the restrictions $T\rest (S^1\times \{j\})$
for $j\in M$ are the \emph{loops}. If $T=(T,N,M,\ell)$ is a t-diagram
and $S=T\rest (I\times \{i\})$ is a strand of $T$, then
$(S,\{i\},\es,\ell_S)$ is also a t-diagram where
$\ell_S=\ell_T\rest \Cross(S)$. 
Similarly, if $j\in M$,
$(S,\es,\{j\},\ell_S)$ is a t-diagram where
$S=T\rest (S^1\times \{j\})$.  Denote by
$E(T)=|T|\cap \partial \square$ the set of endpoints.  A t-diagram is
\emph{loopless} if $M=\es$. In this case we may drop $M$ from the
notation and consider triples $(T,N,\ell)$. A knot diagram is a
special case of a loopless t-diagram $T$ with $N=1$, and a link
diagram is a special case of a (not-necessarily-loopless) t-diagram
with $N=1$. As mentioned in Section~\ref{ssec:Basics}, we often identify
$T$ with its range $|T|$ and drop the vertical bars from the notation.

Note that in \ref{def:Tangle4} we did not require equivalence (we
required ``only if'' as opposed to ``if and only if''), so there can
be loose ends in a t-diagram, i.e. boundary points of $D_T$ which are
mapped into the interior of $\square$.  Thus, the partial diagrams
defined in Section \ref{ssec:Partial} are also special cases of
t-diagrams. A t-diagram is called \emph{complete} if there are no
loose ends and $T^{-1}[\partial\square]=\partial D_T$.  For a complete
t-diagram we have $\Card(E(T))=2\Card(N)$.

\subsection{Operations on Tangle Diagrams}

\begin{Def}\label{def:AxisAligned}
  A rectangle $A\subset \R^2$ is \emph{axis-aligned} if it is of the
  form $[a,b]\times [c,d]$ for some $a<b$ and $c<d$. A homeomorphism
  $h\colon A\to B$ between axis-aligned rectangles
  $A=[a,b]\times [c,d]$ and $A'=[a',b']\times [c',d']$ is called
  \emph{axis-aligned} if there are increasing homeomorphisms
  $h_1\colon [a,b]\to [a',b']$ and $h_2\colon [c,d]\to [c',d']$ such
  that $h(x,y)=(h_1(x),h_2(y))$ for all $(x,y)\in A$.
\end{Def}

\begin{Fact}\label{fact:AxisAligned}
  Given two axis-aligned rectangles $A,A'$ there is a unique linear
  axis-aligned homeomorphism between them. \qed
\end{Fact}

\begin{Def}\label{def:G-tangle}
  Let $G$ be an axis-aligned rectangle. A tuple $(T,N,M,\ell)$ is a
  \emph{$G$-tangle diagram} (\emph{$G$-t-diagram} for short) if
  $T\subset G$ and $(h\circ T,N,M,\ell\circ h^{-1})$ is a t-diagram
  where $h$ is the unique linear axis-aligned homeomorphism
  $h\colon G\to\square$.
\end{Def}

\begin{Def}\label{def:TcapG-tangle}
  Suppose $T=(T,N,M,\ell)$ is a t-diagram and $G\subset\square$ an
  axis-aligned rectangle. We say that $T$ is in \emph{general
    position} with respect to $G$ if $\Cross(T)\cap \partial G=\es$,
  $T$ does not intersect $G$ in the corners, and each intersection of
  $T$ with $\partial G$ is transversal. In this case, $T^{-1}[G]$ is a
  finite union of nondegenerate connected closed submanifolds of
  $D_T$.  Thus, there exist $N'$ and $M'$, and a homeomorphism
  $\eta\colon (I\times N')\cup (S^1\times M')\to D_T$. Define
  $$T\capdot G$$
  to be the $G$-t-diagram
  $(T\circ\eta,N',M',\eta^{-1}\circ (\ell\rest G))$.  Note that the
  choice of $\eta$ is not unique, but it does not matter for our
  purposes. The choice of $N'$ and $M'$ is unique.
\end{Def}

\subsection{Type of a Tangle}

Denote the vertical and horizontal open edges of the square as follows:
\begin{align*}
  V^-&=\{(0,t)\mid 0<t<1\}&&=\text{vertical left edge}\\
  V^+&=\{(1,t)\mid 0<t<1\}&&=\text{vertical right edge}\\
  H^-&=\{(t,0)\mid 0<t<1\}&&=\text{horizontal lower edge}\\
  H^+&=\{(t,1)\mid 0<t<1\}&&=\text{horizontal upper edge}
\end{align*}
Note that the square's corner points are not included in these sets.

\begin{Def}\label{def:Type}
  Let $n_1,n_2,n_3,n_4\in\N$. We say that a t-diagram $T=(T,N,M,\ell)$ is
  \emph{of type} $(n_1,n_2,n_3,n_4)$ if
  $$
    \Card(T\cap H^+)=n_1,\
    \Card(T\cap V^-)=n_2,\ 
    \Card(T\cap H^-)=n_3,\text{ and }
    \Card(T\cap V^+)=n_4.
  $$
\end{Def}

Recall that the t-diagram avoids the corners and there are no
multiplicity points on $\partial\square$, so we necessarily have
$n_1+n_2+n_3+n_4=\Card(E(T))$, and if $T$ is complete, then also
$=2\Card(N)$.

\begin{figure}
  \centering
  \includegraphics[width=0.4\textwidth]{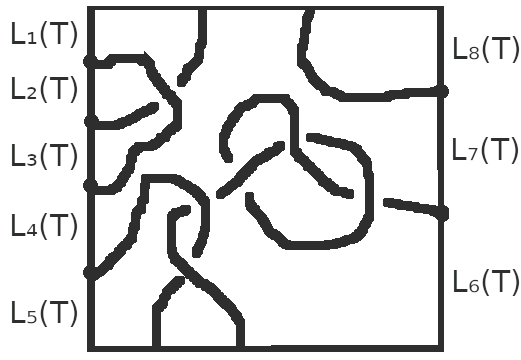}
  \caption{The sets $L_1(T),\dots,L_m(T)$ of a given t-diagram
    illustrated.  They list the connected components of the vertical
    edges of the square from which the endpoints of the tangle have
    been removed.  The enumeration is counterclockwise starting from
    the top left corner.}
  \label{fig:L_T}
\end{figure}

\begin{Def}\label{def:Region}
  Suppose $T=(T,N,M,\ell)$ is a t-diagram
  and let
  $E_V(T)=T\cap (V^-\cup V^+)$ be the set of all its endpoints which
  are on the vertical edges of~$\square$.
  The set $(V^-\cup V^+)\setminus E_V(T)$ is the union of
  disjoint open subarcs of $V^-\cup V^+$.
  Let $L_1(T),\dots,L_{m}(T)$
  be the counterclockwise enumeration of these arcs so that $L_1(T)$
  is the one whose closure contains the left top corner $(0,1)$.
  See Figure~\ref{fig:L_T} for illustration.
\end{Def}

\begin{Def}\label{def:MutualGP}
  A special case of this was already considered in Section~\ref{ssec:CrN}.
  Two t-diagrams $T$ and $S$ are in \emph{mutual general position} if
  $E(S)\cap E(T)=\es$, $T\cap\Cross(S)=S\cap \Cross(T)=\es$, there are
  no more than finitely many $x\in T\cap S$ and each such $x$ is a
  transversal double point between an arc of $T$ and an arc of~$S$
  unless $x$ is a loose point in which case we require that it is a
  loose point of both. For example, if $K$ is a knot diagram
  (Section~\ref{ssec:KD}), then $K^+$ and $K^-$ are in mutual general
  position, even though $o\in K^+\cap K^-$ is not a transversal double
  point.
  When two t-diagrams $T$ and $S$ are in mutual general
  position, denote
  \begin{equation}
      \label{def:Intersect}
      \intersect(T,S)=\Card(\ocirc{T}\cap \ocirc{S}) 
  \end{equation}
  where by $\ocirc T$ in this context we mean the restriction of $T$
  to $D_T\setminus \partial D_T=$ the interior points of the manifold
  precisely because we do not want $\intersect(\cdot,\cdot)$ to count
  the possible loose points.  We always assume by default that all
  t-diagrams are in mutual general position in a given context.
\end{Def}

\subsection{Applying a Homeomorphism to a Tangle}
\label{sec:ApplyHomeo}

Suppose $h\colon \square\to\square$ is a homeomorphism and
$(T,N,M,\ell)$ a t-diagram. Then we denote by $hT$ the t-diagram
$(h\circ T,N,M,\ell\circ h^{-1})$. Important special cases are the
applications of the moves $R_0$ and $R^*_0$,
(Section~\ref{ssec:RReid}).

\subsection{Reparametrization}
\label{ssec:Reparam}

Sometimes t-diagrams need to be reparametrized. An example of this was
the restriction of a t-diagram to a rectangle $T\capdot G$,
Definition~\ref{def:TcapG-tangle}. Another example is smoothing. A
\emph{local smoothing} of a t-diagram is a local change at a crossing
as depicted in Figure~\ref{fig:LocalSmoothing}. If t-diagrams
$(T,N,M,\ell)$ and $(T',N',M',\ell')$ are the same t-diagram before
and after smoothing, the number of strands and loops may or may not
change, i.e. we might have $N\ne N'$ and/or $M\ne M'$, or they might
not change etc. We will not give the technical details of the
reparametrization and assume that the reader understands how it can be
done in each specific case.  A smoothing operation may connect two
strands which were originally oriented in opposite directions. Then
one of them changes orientation after merging with the other
strand. This is not a problem, as in this context we do not care about
the orientation of strands and loops in t-diagrams.

Similarly for curves. We introduced the notations
$\gamma_1\cat \gamma_2$ and $f\cup g$ in the beginning of this
section, but sometimes the domain of the curves may not be exactly
equal to $I$, or they may be oriented in the ``wrong'' direction,
i.e. $\gamma_1(0)=\gamma_2(0)$ and we want to connect them at their
beginning points. It should always be clear from the
context how to (re-)parametrize the curves in a proper manner. 

In above examples the orientation is not of crucial importance.
Sometimes, however, we use orientation and related properties as a
technical tool.  For example in the definition of a \emph{centered
  knot diagram} (see Sections~\ref{ssec:KD} and \ref{ssec:Partial}) we
required that the points $0,\frac{1}{2},1\in I$ are always mapped
respectively to the bottom, the middle, and the top of the square
which indirectly implies an orientation.

\section{Bounding the Number of Non-Hybrid Crossings}
\label{sec:NonHybrid}

To prove Theorem \ref{thm:Main}, we will first bound the number of the
``non-hybrid crossings'', i.e. the crossings in $\widehat K^+$ and
$\widehat K^-$. It is done in this section and is the easier part
of the proof of Theorem \ref{thm:Main}.

\subsection{The Minimal Ray Intersection Number}

\begin{Def} \label{def:MINRay}
  A \emph{ray} is curve without self-intersections in
  $\square$ which has one endpoint at $o$ and another at $\partial\square$.
  For a t-diagram $T$, let
  \begin{align*}
    \MIN_{\ray}(T)&=\min\{\intersect(P, J)\mid J \text{ is a ray}\}.
  \end{align*}
  Call it the \emph{minimal ray intersection number.}
  As usual, we only consider rays which are in general position with
  respect to~$T$. Note that by the definition of $\intersect(\cdot,\cdot)$, \eqref{def:Intersect},
  the point $o$ is not counted in $\intersect(T,J)$. 
\end{Def}

Recall that by $\sim^*_0$ we denote the equivalance relation on
(partial) knot diagrams which is induced by the $R_0$-moves alone,
i.e.  $K\sim^*_0 K'$ iff there is a self-homeomorphism
$h\colon \square\to\square$ such that $h$ is the identity on
some neighborhood of $\partial \square\cup \{o\}$ and $K'=hK$.

\begin{Lemma}[Invariance under $R_0$]\label{lemma:invariance_under_R0}
  Let $K$ be a centered knot diagram and $\widehat K\sim^*_0 K$. Then
  $\MIN_{\ray}(\widehat K^{\pm})=\MIN_{\ray}(K^{\pm})$.
\end{Lemma}
\begin{proof}
  A sequence of $R_0$-moves is equivalent to one $R_0$-move,
  because the composition of witnessing homeomorphisms is a witnessing
  homeomorphism. So suppose that $\widehat K$ has been obtained from $K$
  using such a homeomorphism~$h$.  Recall that by the
  definition of the $R_0$-move, $h$ keeps $o$ fixed, $J$ is a
  ray if and only if $hJ$ is a ray and
  $\intersect(K^{\pm},J)=\intersect(hK^{\pm},hJ)$. So the statement follows from the
  obvious equality $hK^{\pm}=\widehat K^{\pm}$.
\end{proof}

\subsection{The Bound on $\chi(K^\pm)$}
\label{ssec:Bound}

Recall the definition of a composite knot diagram from
Section~\ref{ssec:CompositeKD} and the definition of the equivalence
relations $\sim^X_\RR$ for $X\in \{*,o1,o2,s1,s2\}$ from
Section~\ref{ssec:RReid}. Recall also the notations $\hat o$ and
$\hat o^{\pm}$ from Section~\ref{ssec:TopInt}.

\begin{Prop}\label{prop:KodownJ}
  Suppose $K$ is a composite centered knot diagram and suppose $\widehat K\sim^X_\RR K$
  for some $X\in \{*,o1,o2,s1,s2\}$. Let $K_i$, $i\in \{0,1\}$, be the factors of $K$,
  i.e. $K=K_0\hash K_1$.
  Then there are knot diagrams $\widehat K_0$ and $\widehat K_1$ such that
  $\widehat K_0\sim_{\RR} K_0$, $\widehat K_1\sim_{\RR}K_1$,
  and the following equalities hold:
  \begin{align}
    \chi(\widehat K_1)&= \chi(\widehat K^+)+\MIN_{\ray}(\widehat K^+)\label{def:IMPORTANT1}\\
    \chi(\widehat K_0)&= \chi(\widehat K^-)+\MIN_{\ray}(\widehat K^-)\label{def:IMPORTANT2}
  \end{align}
\end{Prop}
\begin{proof}
  The case $X=*$ will follow from any of the other cases.  Suppose
  that $X=o1$. Then the subarcs of $K^+$ are not allowed to slide
  \emph{under} the critical point $o$. Let $\KK$ and $\widehat\KK$ be
  knots in $C=I^3$ such that $\DD(\KK)=K$ and
  $\DD(\widehat\KK)=\widehat K$.  Then there is a boundary-fixing
  homeomorphism $h\colon C\to C$ such that
  $\widehat\KK=h\widehat\KK$. By the assumption, we can assume that
  $h$ fixes the tube $\tau=B(o,\e)\times I^{-}$ for some sufficiently
  small~$\e$. We can also assume that there is a neighborhood $U$ of
  $\partial C$ which is fixed pointwise by~$h$.  Let $J$ be a ray
  such that $\intersect(K^+,J)=\MIN_{\ray}(K^+)$. By modifying $J$ in
  a neighborhood of $\partial\square$, we can assume without loss of
  generality that the other endpoint of $J$ is at $o^-$ and that
  $(J,\ell_{\es})$ satisfies the conditions of being a lower partial
  knot diagram without crossings.  It is now easy to find a curve
  $\JJ$ in $C$ such that one endpoint of $\JJ$ is at $\hat o$, the
  other endpoint is on $\partial C$, $\JJ$ is contained in the set
  $\tau\cup U$, and $\pr\JJ=J$. Let $\KK_1=\KK^+\cup \JJ$ and
  $\widehat\KK_1=\widehat\KK^+\cup \JJ$.  We now have
  $\widehat KK_1=h\KK_1$. On the other hand, it is also clear that
  $\DD(\KK_1)$ is equivalent to the composition of $K_1$ with a
  trivial knot, so $\DD(\KK_1)\sim_\RR K_1$ and hence also
  $\DD(\widehat \KK_1)\sim_{\RR} K_1$.  Let
  $\widehat K_1=\DD(\widehat \KK_1)$. So we have
  $\widehat K_1\sim_{\RR}K_1$.  To compute the number of crosings in
  $\widehat K_1$, note that $\widehat K_1^-=J$, so by
  Remark~\ref{rem:mirror},
  $$\chi(\widehat K_1)=\chi(\widehat K_1^+)+\chi(J)+\intersect(\widehat K^+_1,J)$$
  But $\widehat K_1^+=\widehat K^+$, $\chi(J)=0$, and
  $\intersect(\widehat K^+,J)=\MIN_{\ray}(\widehat K^+)$, so we have
  $$\chi(\widehat K_1)=\chi(\widehat K^+)+\MIN_{\ray}(\widehat K^+)$$
  as required. By a symmetric argument, obtain $\widehat K_0$
  with $\widehat K_0\sim_{\RR}K_0$ which satisfies \eqref{def:IMPORTANT2}.
  The cases $X\in \{o2,s1,s2\}$ are obtained by symmetric arguments.
\end{proof}

\begin{Cor}[A Bound for Non-Hybrid Crossings]\label{cor:Bound}
  Suppose that $K$ is a composite centered knot diagram and
  $K=K_0\hash K_1$. Suppose that $\widehat K\sim^X_{\RR} K$
  for some $X\in \{*,o1,o2,s1,s2\}$.
  Then
  \begin{align}
    \chi(\widehat K^+)+\MIN_{\ray}(\widehat K^+)&\ge c(K_1)\label{def:IMPORTANT11}\\
    \chi(\widehat K^-)+\MIN_{\ray}(\widehat K^-)&\ge c(K_0)\label{def:IMPORTANT22}
  \end{align}
\end{Cor}
\begin{proof}
  Let $\widehat K_0$ and $\widehat K_1$ be as given by Proposition~\ref{prop:KodownJ}.
  Since $\widehat K_0\sim_{\RR}K_0$ and $\widehat K_0\sim_{\RR}K_0$,
  it follows (from \eqref{eq:MinCrossingNumber}) that
  $\chi(\widehat K_0)\ge c(K_0)$ and $\chi(\widehat K_1)\ge c(K_1)$.
  The statement follows now from \eqref{def:IMPORTANT1} and~\eqref{def:IMPORTANT2}.
\end{proof}

\section{Crossings Between Overlapping Tangle Diagrams}
\label{sec:BoundingTangles}

We have now established a lower bound on the number of non-hybrid
crossings in~$\hat K$.  Next, we will have to establish a lower bound
for the number of hybrid crossings
$\intersect(\widehat K^-,\widehat K^+)$. Together they will lead to
the proof of Theorem~\ref{thm:Main} using the equation in
Remark~\ref{rem:mirror}.  This will be done in
Section~\ref{sec:HybridCrossings}. In this section we prove
preliminary results needed. The main objective is to
prove Proposition \ref{prop:Hybrid} 
which gives a lower bound for the number of crossings between two
``overlapping'' tangle diagrams under certain conditions. To prove
it, we will need a lemma about
integer-valued matrices (Section~\ref{ssec:Matrices},
Lemma~\ref{lemma:KeyLemma}) and a lemma about tangle diagrams which
enables the application of Lemma~\ref{lemma:KeyLemma} to the context
of t-diagrams (Section \ref{ssec:TangleDiagrams},
Lemma~\ref{lem:ConnectingLemma1}).

\subsection{Multigraphs, Matrices and Zigzag Paths}
\label{ssec:Matrices}
 
Denote by $\N$ natural numbers including zero and by
$\N_+=\N\setminus \{0\}$ natural numbers excluding zero.  For
$n,m\in\N$ denote by $\dint{n,m}$ the discrete set
$\{n,n+1,\dots,m-1,m\}$.  By convention, if $m<n$, then $\dint{n,m}$
is the empty set and if $m=n$, then it is the singleton
$\{n\}=\{m\}$. Recall that $\Card(A)$ the cardinality of the set $A$,
e.g. for $n\le m$ we have
\begin{equation}
  \label{eq:dintcard}
  \Card(\dint{n,m})=m-n+1.
\end{equation}
The $ij$-th element of an $(n\times n)$-matrix $e$ is expressed as
$e(i,j)$ where $i,j\in \dint{1,n}$. Let $\MM^s_n(\N)$ be the set of
symmetric $\N$-valued $(n\times n)$-matrices,
$$\MM^s_n(\N)=\{x\in \N^{(n\times n)}\mid \forall i,j\in\dint{1,n}(x(i,j)=x(j,i))\}.$$
A useful intuitive interpretation of a square
$\N$-valued matrix is that it represents a multigraph. In this
interpretation, the value $e(i,j)$ is the number of connections from
node $i$ to node $j$.  This is, intuitively, how the results of this section
are going to be applied in the context of t-diagrams (see
Definition~\ref{def:eTS}).

Throughout this section we fix $k\in\N_+$ and let $n=2k$. Let
$f\colon \dint{1,n}\to \dint{1,n}$ be the function defined by
$f(i)=n-i+1$. Note that
\begin{equation}
  f(f(i))=i\label{eq:involution}
\end{equation}
for all $i\in \dint{1,n}$, in particular $f$ is a bijection. We will also need
the fact that
\begin{equation}
  i\in \dint{2,k-1}\iff f(i)\in \dint{k+2,n-1}.\label{eq:involution2}
\end{equation}
For an
$(n\times n)$-matrix $e$, denote by $s_i(e)$ the sum of the $i$-th row of~$e$,
\begin{equation}
  s_i(e)=\sum_{j=1}^n e(i,j).\label{eq:sumofstrands}
\end{equation}
In the multigraph interpretation, $s_i(e)$ is the total number of
out-going connections from $i$.

\begin{Def}\label{def:PropertyStar}
  A matrix $e\in \MM^s_n(\N)$ is said to satisfy property $(*)$ if
  the following conditions hold:
  \begin{enumerate}[label={\upshape ($*_{\arabic*}$)}, leftmargin=3pc]
  \item For all $i\in \dint{2,k-1}$, $s_i(e)=s_{f(i)}(e)$ \label{item:star2}
  \item If $k>1$, then
    $s_1(e)=s_{f(1)}(e)+1$ and $s_{k}(e)=s_{f(k)}(e)+1$
    \label{item:star3}
  \item For all $i\in \dint{1,n}$, $e(i,i)$ is even. \label{item:star4}
  \end{enumerate}
\end{Def}

\begin{Def}\label{def:ZigZagpath}
  Given $e\in \MM^s_n(\N)$, a \emph{zigzag $e$-path} is a sequence
  $(i_1,\dots,i_p)\subset \dint{1,n}$ such that for all
  $q\in \dint{1,p-1}$, $e(f(i_q),i_{q+1})>0$.
\end{Def}

For $a,b\in\MM^s_n(\N)$, denote $a\le b$ if $a(i,j)\le b(i,j)$ holds
for all $i,j\in \dint{1,n}$, and $a<b$ if $a\le b$ and
$a\ne b$. 


\begin{Lemma}\label{lemma:ZigZag}
  Suppose that $e\in\MM^s_n(\N)$ satisfies
  property~$(*)$.  Then there is a zigzag $e$-path $(i_1,\dots,i_p)$
  such that $i_1=f(1)$ and $i_p \in \{k,f(k)\}$.
\end{Lemma}
\begin{proof}
  In case $k=1$, let $P=(i_1)$ be the path of length $p=1$, in which
  $i_1=f(1)=n=2$.  Then $P$ is trivially a zigzag $e$-path, and
  $i_1=f(1)=f(k)=i_p$, so we can assume that $k>1$.
  Suppose for a contradiction that there is no zigzag path as
  required.  We will inductively define an infinite sequence
  $(x_i)_{i\in\N}\subset\MM^s_n(\N)$ such that $x_1>x_2>\cdots$.  This
  is a contradiction, because there cannot be infinite strictly descending
  sequences in $\MM^s_n(\N)$ due to the well-foundedness of $(\N,<)$.
  If $m\in \dint{1,n}\setminus \{k,k+1\}$, we say that
  $x\in \MM^s_n(\N)$ satisfies the \emph{$m$-modified property} $(*)$
  called $(*^m)$ if:
  \begin{enumerate}[label={\upshape ($*^m_{\arabic*}$)}, leftmargin=3pc]
  \item For all $i\in \dint{1,n}\setminus \{m,f(m),k,k+1\}$, $s_i(x)=s_{f(i)}(x_q)$ \label{item:star2q}
  \item $s_{f(m)}(x)=s_{m}(x)+1$ and $s_{k}(x)=s_{f(k)}(x)+1$.
    \label{item:star3q}
  \item (same as \ref{item:star4}): For all $i\in \dint{1,n}$, $x(i,i)$ is even. \label{item:star4q}
  \end{enumerate}
  Let $P_1=(i_1)$ and $X_1=(x_1)$, in which $i_1=f(1)=n$ and
  $x_1=e$. Now using \eqref{eq:involution} and \eqref{eq:involution2},
  it is easy to verify that $(*^{i_1})=(*^n)$ is the same as $(*)$, so
  it applies to $e$.  It is also clear that $P_1$ is a zigzag
  $e$-path.  Suppose that we have defined $P_q=(i_1,\dots,i_q)$ and
  $X_q=(x_1,\dots,x_q)$ such that $P$ is a zigzag $e$-path,
  $x_1>x_2>\dots$, and $x_q$ satisfies $(*^m)$ for $m=i_q$ (note that
  by the counterassumptions, $i_q\notin \{k,k+1\}$). Then by
  \ref{item:star3q} for $x_q$ we have $s_{f(m)}(x_q)=s_{m}(x_q)+1$.
  Since $s_{m}(x_q)$ is non-negative,
  $s_{f(m)}(x_q)>0$. By \eqref{eq:sumofstrands} this means that
  there must exist $m'\in \dint{1,n}$ with
  \begin{equation}
    x_q(f(m),m')>0.\label{eq:positive}  
  \end{equation}
  Let $i_{q+1}=m'$. Since $e=x_1>x_q$, also $e(f(m),m')>0$ holds, so
  $P_{q+1}=(i_1,\dots,i_q,i_{q+1})$ is a zigzag $e$-path. Note that by our
  counterassumption,
  \begin{equation}
    \label{eq:mnotk}
    m'\notin \{k,f(k)\}. 
  \end{equation}
  Let $x_{q+1}$ be defined as follows
  \begin{numcases}{x_{q+1}(i,j)=}
    x_q(i,j)-1&if $(i,j)\in \{(f(m),m'),(m',f(m))\}$ and $f(m)\ne m'$,\nonumber\\
    x_q(i,j)-2&if $(i,j)\in \{(f(m),m'),(m',f(m))\}$ and $f(m)=m'$,\nonumber\\
    x_q(i,j)&otherwise.\nonumber
  \end{numcases}
  Clearly $x_{q}>x_{q+1}$. By \eqref{eq:positive}, \ref{item:star4q}, and the symmetry
  of $x_q$ we have $x_{q+1}\in\MM^s_n(\N)$. It remains to verify that
  $x_{q+1}$ satisfies $(*^{m'})$. The condition \hyperref[item:star4q]{$(*^{m'}_{3})$} for $x_{q+1}$
  follows trivially from the assumption that
  \ref{item:star4q} holds for $x_{q}$. So let us prove
  \hyperref[item:star2q]{$(*^{m'}_1)$} and \hyperref[item:star3q]{$(*^{m'}_2)$} for $x_{q+1}$.
  Suppose first that $f(m)\ne m'$. Then the definition of $x_{q+1}$ means that the number $1$
  has been subtracted from exactly two rows in $x_q$ to obtain $x_{q+1}$.
  These rows are $f(m)$ and $m'$. This means that
  \begin{align}
    s_{f(m)}(x_{q+1})&=s_{f(m)}(x_{q})-1 \label{obs1}\\
    s_{m'}(x_{q+1})&=s_{m'}(x_q)-1 \label{obs2}\\
    s_i(x_{q+1})&=s_i(x_q)&&\text{for all }i\notin \{f(m),m'\} \label{obs3}
  \end{align}
  By the induction hypothesis we have
  \begin{align}
    s_{f(m)}(x_q)&=s_{m}(x_q)+1, \label{ind_hyp}\\
    s_{k}(x_q)&=s_{f(k)}(x_q)+1\label{ind_hyp3}\\
    s_{i}(x_q)&=s_{f(i)}(x_q)&& \text{for all }i \notin  \{m,f(m),k,f(k)\}.\label{ind_hyp2}
  \end{align}
  To verify \hyperref[item:star2q]{$(*^{m'}_1)$}, suppose that
  \begin{equation}
    i\in \dint{1,n}\setminus \{m',f(m'),k,f(k)\}\label{eq:choiceofi}
  \end{equation}
  If $i\notin \{m,f(m)\}$,
  then
  $$s_i(x_{q+1})\stackrel{\eqref{obs3}}{=}s_i(x_{q})\stackrel{\eqref{ind_hyp2}}{=}s_{f(i)}(x_q)\stackrel{\eqref{obs3}}{=}s_{f(i)}(x_{q+1}).$$
  Condition \eqref{obs3} can be applied in the last equality, because
  $i\notin \{f(j),j\}$ implies $f(i)\notin \{f(j),j\}$ for any $j$ by
  \eqref{eq:involution}. Suppose $i=m$. Note that then $m\ne m'$ by
  the choice of $i$ \eqref{eq:choiceofi}. Also then
  $i\ne f(m)$, so we have
  $$s_i(x_{q+1})\stackrel{\eqref{obs3}}{=}s_i(x_{q})\stackrel{\eqref{ind_hyp}}{=}s_{f(i)}(x_q)-1\stackrel{\eqref{obs1}}{=}s_{f(i)}(x_{q+1}).$$
  Suppose that $i=f(m)$. By the choice of $i$ this implies again
  $m\ne m'$, because $i$ cannot be equal to $f(m')$. Note that by \eqref{eq:involution}
  we have now $f(i)=m$, so:
  $$s_i(x_{q+1})\stackrel{\eqref{obs1}}{=}s_i(x_{q})-1\stackrel{\eqref{ind_hyp}}{=}s_{f(i)}(x_q)\stackrel{\eqref{obs3}}{=}s_{f(i)}(x_{q+1}).$$
  This proves \hyperref[item:star2q]{$(*^{m'}_1)$} for $x_{q+1}$.  Now consider
  \hyperref[item:star3q]{$(*^{m'}_2)$}. We are still working under the assumption that $m'\ne f(m)$.
  Since $m,m'\notin \{k,f(k)\}$, it also follows
  that $f(m)\notin \{k,f(k)\}$ and symmetrically that $k\notin \{m,f(m),m',f(m')\}$, so
  $$s_k(x_{q+1})\stackrel{\eqref{obs3}}{=}s_k(x_q)\stackrel{\eqref{ind_hyp3}}{=}s_{f(k)}(x_q)+1\stackrel{\eqref{obs3}}{=}s_{f(k)}(x_{q+1}).$$
  This verifies the second part of \ref{item:star3q}. For the first part, first consider
  the case $m\ne m'$. Then $f(m)\ne f(m')$ and by \eqref{eq:mnotk} and the assumption
  $m'\ne f(m)$ we have
  $$s_{f(m')}(x_{q+1})\stackrel{\eqref{obs3}}{=}s_{f(m')}(x_{q})\stackrel{\eqref{ind_hyp2}}{=}s_{m'}(x_q)
  \stackrel{\eqref{obs2}}{=} s_{m'}(x_{q+1})+1.
  $$
  Suppose that $m=m'$. Then $f(m)=f(m')$ and we have
  $$s_{f(m')}(x_{q+1})\stackrel{\eqref{obs1}}{=}s_{f(m')}(x_{q})-1\stackrel{\eqref{ind_hyp}}{=}s_{m'}(x_q)
  \stackrel{\eqref{obs2}}{=} s_{m'}(x_{q+1})+1.
  $$
  We have proved $(*^{m'})$ for $x_{q+1}$ in the case $m'\ne f(m)$.
  Let us now assume that $f(m)=m'$. Then the definition of $x_{q+1}$ means
  that the number $2$ has been subtracted frmo exactly one row of $x_q$ to obtain $x_{q+1}$
  that row being $m'=f(m)$. This means:
  \begin{align}
    s_{m'}(x_{q+1})&=s_{m'}(x_q)-2=s_{f(m)}(x_q)-2\label{obs1_}\\
    s_i(x_{q+1})&=s_i(x_q)&&\text{for all }i\ne m'. \label{obs3_}
  \end{align}
  Substituting $m'$ for $f(m)$ in \eqref{ind_hyp}--\eqref{ind_hyp2} 
  the induction hypothesis gives:
  \begin{align}
    s_{m'}(x_q)&=s_{m}(x_q)+1, \label{ind_hyp_}\\
    s_{k}(x_q)&=s_{f(k)}(x_q)+1\label{ind_hyp3_}\\
    s_{i}(x_q)&=s_{f(i)}(x_q)&& \text{for all }i \notin  \{m,m',k,f(k)\}.\label{ind_hyp2_}
  \end{align}
  To verify \hyperref[item:star2q]{$(*^{m'}_1)$}, suppose that
  \begin{equation}
    i\in \dint{1,n}\setminus \{m',f(m'),k,f(k)\}.\label{eq:choiceofi}
  \end{equation}
  The assumption $m'=f(m)$ imples $f(m')=m$, so we also have $i\ne m$.
  By \eqref{ind_hyp2_} we have $s_{i}(x_q)=s_{f(i)}(x_q)$. Since $i\notin \{m',f(m')\}$,
  also $f(i)\notin \{m',f(m')\}$, so \eqref{obs3_} implies the desired equality.
  For  \hyperref[item:star3q]{$(*^{m'}_2)$}, consider:
  \begin{align*}
    s_{f(m')}(x_{q+1})&=s_{m}(x_{q+1})\stackrel{\eqref{obs3_}}{=}s_{m}(x_q)
    \stackrel{\eqref{ind_hyp_}}{=}s_{m'}(x_q)-1\stackrel{\eqref{obs1_}}{=}s_{m'}(x_{q+1})+1.
  \end{align*}
  The part for $k$ follows as above because again $m'=f(m)\notin \{k,f(k)\}$.
\end{proof}

\begin{Def}\label{def:XiMu}
  Let $w\in \MM^s_n(\N)$ be defined for all $i\le j\in \dint{1,n}$ by
  \begin{numcases}{w(i,j)=}
    |j-i|& if $i\le j\le k$ or $k<i\le j$\label{def:w1}\\
    |f(j)-i|+1& if  $i\le k<j$\label{def:w2}
  \end{numcases}
  and if $j<i$, then let $w(i,j)=w(j,i)$.
  For any $e\in\MM^s_n(\N)$, let
  \begin{equation}
    \mu(e)=\sum_{i=1}^k\sum_{j=k+1}^ne(i,j)\label{eq:defmu}
  \end{equation}
  and
  \begin{equation}
    \xi(e)=\sum_{i=1}^n\sum_{j=i+1}^ne(i,j)w(i,j).\label{eq:defxi}  
  \end{equation}
\end{Def}

We are ready to prove the following key lemma.

\begin{Lemma}\label{lemma:KeyLemma}
  Suppose that $e\in\MM^s_n(\N)$ satisfies property $(*)$.
  Then $\xi(e)\ge \mu(e)+k-1$.
\end{Lemma}
\begin{proof}
  Let us first consider the case $k=1$.
  Then $e$ and $w$ are $(2\times 2)$-matrices, in fact
  $$w=\left(
    \begin{bmatrix}
      0&1\\
      1&0\\
    \end{bmatrix}
  \right)$$
  and
  $$\xi(e)=\sum_{i=1}^2\sum_{j=i+1}^2e(i,j)w(i,j)=e(1,2)=\mu(e)=\mu(e)+1-1=\mu(e)+k-1$$
  as requried. So we can assume that $k>1$.

  Denote $\cS=\{(i,j)\in \dint{1,n}^2\mid i<j\}$ and
  define $\OO,\XX\subset \cS$ to be the sets
  $\OO=\{(i,j)\in \cS\mid i, j\le k\text{ or }i,j>k\}$,
  $\XX=\cS\setminus \OO$. We now have
  $$\xi(e)=\sum_{(i,j)\in\cS}e(i,j)w(i,j)\quad\text{ and }\quad \mu(i,j)=\sum_{(i,j)\in \XX}e(i,j).$$  
  So it is enough to show that
  \begin{equation}
    \sum_{(i,j)\in\cS}e(i,j)w(i,j)\ge \sum_{(i,j)\in \XX}e(i,j) + k-1.\label{eq:SSXXForm}  
  \end{equation}
  Let
  $\pi\colon \dint{1,n}\to\dint{1,n}$ be defined by
  $$\pi(i)=
  \begin{cases}
    i&\text{if }i\le k\\
    f(i)&\text{otherwise.}\\
  \end{cases}
  $$
  Note that since
  \begin{equation}
    \label{eq:Flip}
    \forall i\in \dint{1,n}(i\le k\iff f(i)>k)
  \end{equation}
  $\pi$ is a kind of projection in the sense that it satisfies for all $i\in \dint{1,n}$
  \begin{equation}
    \label{eq:projection}
    \pi(i)=\pi(\pi(i))\quad\text{ and by \eqref{eq:involution} also }\quad \pi(f(i))=\pi(i).
  \end{equation}
  Let $\rho,\sigma\in\MM^s_n(\N)$ be defined by
  $\rho(i,j)=|\pi(i)-\pi(j)|$ and
  $$\sigma(i,j)=
  \begin{cases}
    0&\text{if }(i,j)\in \OO\text{ or }(j,i)\in\OO\text{ or }i=j\\
    1&\text{otherwise.}
  \end{cases}
  $$
  \begin{ClaimL}\label{cl:w_is_rho_p_s}
    $w=\rho+\sigma$.
  \end{ClaimL}
  \begin{proof}
    We want to show that for all $i,j\in \dint{1,n}$,
    $w(i,j)=\rho(i,j)+\sigma(i,j)$.  By symmetry it is enough to
    consider $i\le j$. If $i=j$, then by definitions, $w(i,j)=\rho(i,j)=\sigma(i,j)=0$.
    So we can assume that $(i,j)\in \cS$.
    Note that
    $$|f(j)-i|=|n-j+1-i|=|n-i+1-j|=|f(i)-j|.$$
    The definition of $w$ can be re-written for all $(i,j)\in\cS$ as
    $$
    w(i,j)=
    \begin{cases}
      |j-i|&\text{if }(i,j)\in\OO\\
      |f(j)-i|&\text{otherwise.} 
    \end{cases}
    $$
    Also note that for all $(i,j)\in\OO$ we have
    \begin{equation}
      |j-i|=|n-n+1-1+j-i|=|n-i+1-(n-j+1)|=|f(i)-f(j)|.\label{eq:Symm}    
    \end{equation}
    We can now easily check the claim case by case:
    \begin{align*}
      &(i,j)\in \OO\land i,j\le k:\\
      &\phantom{=}w(i,j)=|j-i|=|\pi(j)-\pi(i)|=\rho(i,j)+0=\rho(i,j)+\sigma(i,j)\\
      &(i,j)\in \OO\land i,j> k:\\
      &\phantom{=}w(i,j)=|j-i|\stackrel{\eqref{eq:Symm}}{=}|f(j)-f(i)|\stackrel{\eqref{eq:Flip}}{=}|\pi(j)-\pi(i)|=\rho(i,j)+0=\rho(i,j)+\sigma(i,j)\\
      &(i,j)\in \XX\land i\le k<j:\\
      &\phantom{=}w(i,j)=|f(j)-i|+1=|\pi(j)-\pi(i)|+1=\rho(i,j)+1=\rho(i,j)+\sigma(i,j)\tag*\qedhere
    \end{align*}
  \end{proof}
  Given two matrices $a$ and $b$, denote by $ab$ their Hadamard product, i.e. the matrix
  such that $(ab)(i,j)=a(i,j)b(i,j)$. 
  Let us also introduce temporarily the following shorthand notation for a matrix $a$ and a set $A$:
  $\displaystyle{\ssum{A}a}=\sum_{(i,j)\in A}a(i,j)$. Then the statement \eqref{eq:SSXXForm} becomes
  \begin{equation}
    \label{eq:ToProve}
    \ssum{\cS}ew\ge \ssum{\XX}e+k-1
  \end{equation}
  Let $(i_1,\dots,i_p)$ be the zigzag $e$-path given by Lemma~\ref{lemma:ZigZag}.
  Let $\QQ\subset\cS$ be defined by $\QQ=\QQ_+\cup\QQ_-$ where
  $$\QQ_+=\{(f(i_q),i_{q+1})\mid q\in \dint{1,p}\}\ \text{ and }\ \QQ_-= \{(i_{q+1},f(i_q))\mid q\in \dint{1,p}\}.$$
  By the linearity properties of the sum, Claim~\ref{cl:w_is_rho_p_s}, and definitions of
  $\rho$ and $\sigma$ we have
  \begin{align*}
    \ssum{\cS}ew&=\ssum\cS e(\rho+\sigma)=\ssum{\OO}e\rho + \ssum{\XX}e(\rho+1) \\
                &=\ssum{\OO\setminus\QQ}e\rho + \ssum{\OO\cap\QQ}e\rho + \ssum{\XX\setminus \QQ}e(\rho+1)+ \ssum{\XX\cap \QQ}e(\rho+1)\\
                &=\ssum{\OO\setminus\QQ}e\rho + \ssum{\OO\cap\QQ}e\rho + \ssum{\XX\setminus \QQ}e\rho+\underline{\ssum{\XX\setminus \QQ}e}+ \ssum{\XX\cap \QQ}e\rho+\underline{\ssum{\XX\cap \QQ}e}\\
                &=\ssum{\OO\setminus\QQ}e\rho + \ssum{\OO\cap\QQ}e\rho + \ssum{\XX\setminus \QQ}e\rho+ \ssum{\XX\cap \QQ}e\rho+\underline{\ssum{\XX}e}\\
                &\ge  \ssum{\OO\cap\QQ}e\rho + \ssum{\XX\cap \QQ}e\rho+\ssum{\XX}e\\
                &=  \ssum{\QQ\cap \cS}e\rho +\ssum{\XX}e.
  \end{align*}
  At the inequality we just threw away the terms
  $\ssum{\OO\setminus\QQ}e\rho$ and $\ssum{\XX\setminus \QQ}e\rho$.
  The underlined terms are highlighted to make it easier to see which
  ones were combined to obtain what.  In the last step we use the
  equality $\OO\cup \XX=\cS$ which follows from the definitions of
  $\OO$ and $\XX$. So we are left to show that
  \begin{equation}
    \label{eq:lefttoshow}
    \ssum{\QQ\cap \cS}e\rho\ge k-1.
  \end{equation}
  We have
  $$\ssum{\QQ\cap \cS}e\rho=\ssum{\QQ_+\cap\cS}e\rho + \ssum{\QQ_-\cap\cS}e\rho.$$
  If $(i,j)\in\QQ_+\cap \cS$, then for some $q\in \dint{1,p}$,
  $i=f(i_q)$, $j=i_{q+1}$, and $f(i_q)\le i_{q+1}$. In this case
  $$e(i,j)\rho(i,j)=e(f(i_q),i_{q+1})|\pi(f(i_q))-\pi(i_{q+1})|\ge |\pi(f(i_q))-\pi(i_{q+1})|\stackrel{\eqref{eq:projection}}{=}|\pi(i_q)-\pi(i_{q+1})|.$$
  The inequality follows from Definition \ref{def:ZigZagpath} of a zigzag path. 
  If $(i,j)\in \QQ_-\cap \cS$ then
  symmetrically for some $q\in \dint{1,p}$,
  $j=f(i_q)$, $i=i_{q+1}$, and $f(i_q)\ge i_{q+1}$. Then by the symmetricity of $e$ and $\rho$,
  $$e(i,j)\rho(i,j)=e(j,i)\rho(j,i)=e(f(i_q),i_{q+1})|\pi(f(i_q))-\pi(i_{q+1})|\ge |\pi(i_q)-\pi(i_{q+1})|.$$
  If $q\in \dint{1,p}$, then either $(f(i_q),i_{q+1})\in\QQ_+\cap\cS$
  or $(i_{q+1},f(i_q))\in\QQ_-\cap\cS$, but not both, because
  $\cS$ contains exactly one of $(i,j)$ or $(j,i)$ for $i\ne j$. All in all,
  \begin{equation}\ssum{\QQ_+\cap\cS}e\rho + \ssum{\QQ_-\cap\cS}e\rho\ge \sum_{q=1}^p|\pi(i_q)-\pi(i_{q+1})|
    \ge \left|\sum_{q=1}^p\pi(i_q)-\pi(i_{q+1})\right|=|\pi(i_1)-\pi(i_p)|.
    \nonumber
    \end{equation}
  But $(i_1,\dots,i_p)$ was the zigzag $e$-path given by Lemma \ref{lemma:ZigZag},
  so $i_1=f(1)=n$ and $i_p\in \{k,f(k)\}$. By the definition of $\pi$, we have
  $\pi(i_1)=1$ and $\pi(i_p)=k$, so 
  $|\pi(i_1)-\pi(i_p)|=|1-k|=k-1$ and we are done.
\end{proof}

\subsection{Minimal Splitting Curves}
\label{ssec:TangleDiagrams}


We have already defined \emph{rays} which were used to count certain
minimal intersection numbers. In this section we will define a related
notion. A \emph{splitting curve} is an embedding
$J\colon I\to \square$ such that
$J(t)\in \ocirc{\square}\iff t\notin \{0,1\}$.  In other words a
splitting curve $J$ can be seen as a (complete) tangle diagram
$(J,\{1\},\es,\ell_{\es})$ or a knot diagram $(J,\ell_{\es})$ where
$\ell_{\es}=\es$ is the empty labeling function. Thus,
Definition~\ref{def:MutualGP} applies also to splitting curves, so we
can talk about a splitting curve being in general position with
respect to a t-diagram (i.e. they are in mutual general position), or
another splitting curve. Also the notation $\intersect(\cdot,\cdot)$
of Definition \ref{def:MutualGP}\eqref{def:Intersect} generalizes to
splitting curves. Like with t-diagrams, we assume by default that
splitting curves are in general position w.r.t. to t-diagrams in each
given context.

\begin{Def}\label{def:ijSplitting}
  Recall Definition~\ref{def:Region} of $L_1(T),\dots,L_m(T)$. Given
  a t-diagram $T$ and a splitting curve $J$, if the starting
  point of $J$ is in some $L_i(T)$ and the ending point is in some
  $L_j(T)$, then we call $J$ an \emph{$ij$-splitting curve} (of $T$).
\end{Def}

\begin{Def}\label{def:MatrixV_T}
  Given a t-diagram $T=(T,N,M,\ell)$ denote by $\MIN_{ij}(T)$ the
  minimum number of intersections that an $ij$-splitting curve can
  have with $T$, formally
  \begin{equation}
    \label{eq:defMINij}
    \MIN_{ij}(T)=\min \{\intersect(T,J)\mid J\text{ is an $ij$-splitting curve of $T$}\}.
  \end{equation}
  As per our assumption, the curves $J$ in \eqref{eq:defMINij} range
  only over those that are in general position w.r.t.~$T$. Now define
  the \emph{minimal intersection matrix} for a t-diagram $T$,
  denoted $v_T$ to be given by $v_T(i,j)=\MIN_{ij}(T)$.
\end{Def}

A splitting curve $J$ is called \emph{vertical} if $J(0)\in H^+$ and
$J(1)\in H^-$, i.e. it passes from the top horizontal edge and goes to
the bottom horizontal edge splitting the t-diagram vertically. We
introduce  a special notation for the \emph{vertical minimum intersection
  number} of a t-diagram $T$ defined as
\begin{equation}
  \label{eq:VertMIN}
  \MIN_{\vert}(T)=\min \{\intersect(T,J)\mid J\text{ is a vertical splitting curve of }T\}.
\end{equation}
Again, we note that in this definition we only consider splitting
curves in general position with respect to~$T$.  Recall Definition
\ref{def:Type} of \emph{type} of a t-diagram.

\begin{Lemma}\label{lemma:UpperBoundOnMIN}
  If $T$ is a t-diagram of type $(n_1,n_2,n_3,n_4)$, then
  $\MIN_{\vert}(T)\le \min\{n_2,n_4\}$.
\end{Lemma}
\begin{proof}
  Suppose w.l.o.g. that $n_2\le n_4$. Let $\e$ be so small that the
  $\e$-neighborhood $U$ of $V^-$ is such that $U\cap T$ is a disjoint
  union of $n_2$ horizontal lines. This is possible without loss of
  generality by the piecewise linearity assumption. Let $J$ be a
  vertical straight line inside of $U$.  Then
  $\intersect(T,J)=n_2$ proving
  $\MIN_{\vert}(T)\le n_2=\min\{n_2,n_4\}$.
\end{proof}

Recall that for $(n\times n)$-matrices $a$ and $b$ we denote $a\ge b$
if for all $i,j\in \dint{1,n}$ we have $a(i,j)\ge b(i,j)$.  We are now
ready to state the main lemma of this section which is the key
connection between t-diagrams and matrices of
Section~\ref{ssec:Matrices}.

\begin{Lemma}\label{lem:ConnectingLemma1}
  Fix $k\in\N_+$ and $n=2k$. Suppose that $T=(T,N,M,\ell)$ is a
  t-diagram of type $(1,k-1,1,k-1)$ with $\MIN_{\vert}(T)=k-1$.  Let
  $w$ be the matrix of Definition~\ref{def:XiMu}. Then $v_T\ge w$.
\end{Lemma}

The proof is postponed to the end of this section.

\begin{Def}\label{def:SubsetStar}
  Let $T$ and $S$ be t-diagrams, 
  and $\e>0$. We write $S\subset^\e T$
  to mean that there is $C\subset\Cross(T)$ such that the following conditions hold:
  \begin{enumerate-(a)}
  \item $|S|\subset |T|\cup \Cup_{c\in C}\bar B(c,\e)$, \label{item:subsetstar1}
  \item For all $c\in C$, $\Cross(S)\cap \bar B(c,\e)=\es$.
    \label{item:subsetstar2}
  \end{enumerate-(a)}
  We write $S\subset^\e_{\max} T$ if in addition to \ref{item:subsetstar1} and \ref{item:subsetstar2}, all endpoints of $T$ are endpoints of~$S$.
\end{Def}

\begin{Def}\label{def:mBridge}
  A t-diagram $T$ is \emph{trivial} if $\Cross(T)=\es$.  An
  \emph{$\e$-smoothing} of a t-diagram $T$ is a trivial t-diagram $S$ with
  $S\subset^\e T$. It is \emph{maximal} if $E(S)=E(T)$, i.e.  if
  $S\subset^\e_{\max}T$. An $m$-\emph{bridge} for some $m\in\N$ is a
  trivial t-diagram of type $(0,m,0,m)$ in which every strand connects
  $V^-$ to~$V^+$.  The $0$-bridge is the empty t-diagram and is of type
  $(0,0,0,0)$.
\end{Def}

\begin{Lemma}\label{lemma:CrossingsWithSmoothings}
  Suppose $T$ and $S$ are two t-diagrams and
  $S'\subset^\e S$.  If $\e$ is sufficiently small, then
  $\intersect(T,S)\ge \intersect(T,S')$.
\end{Lemma}
\begin{proof}
  Suppose $x\in \Cross(S)$. Then by general position there is $\e_x$
  such that $B(x,\e_x)\cap T=\es$. Let
  $\e=\min\{\e_x\mid x\in\Cross(S)\}$.  The set $\Cross(S)$ is finite
  so the minimum exists. Suppose $y\in S'\cap T$.  Since $y\in T$,
  $y\notin B(x,\e)$ for all $x\in\Cross(S)$, so by the definition of
  $\subset^\e$, $y\in S$. Thus $y\in S\cap T$. This shows that
  $S'\cap T\subset S\cap T$ which implies the statement.
\end{proof}

\begin{figure}
  \centering
  \fbox{\includegraphics[width=0.20\linewidth]{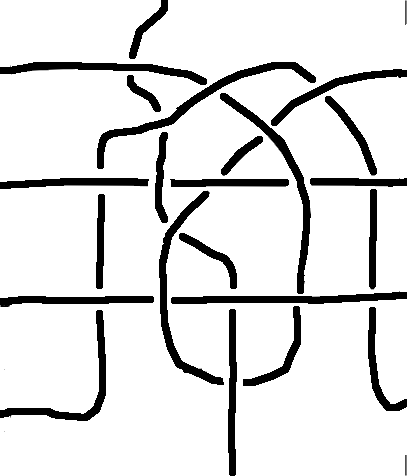}}
  \raisebox{14\height}{$\Rightarrow$}
  \fbox{\includegraphics[width=0.20\linewidth]{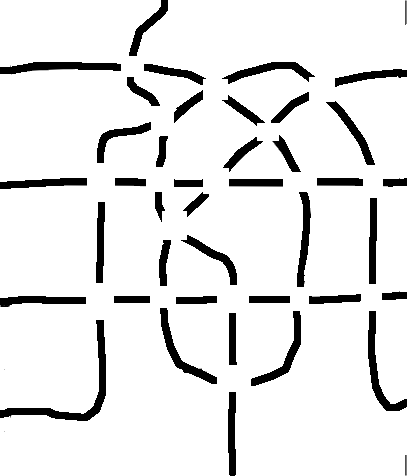}}
  \raisebox{14\height}{$\Rightarrow$}
  \fbox{\includegraphics[width=0.20\linewidth]{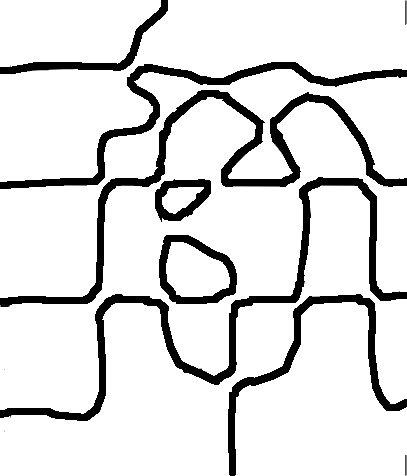}}
  \raisebox{14\height}{$\Rightarrow$}
  \fbox{\includegraphics[width=0.20\linewidth]{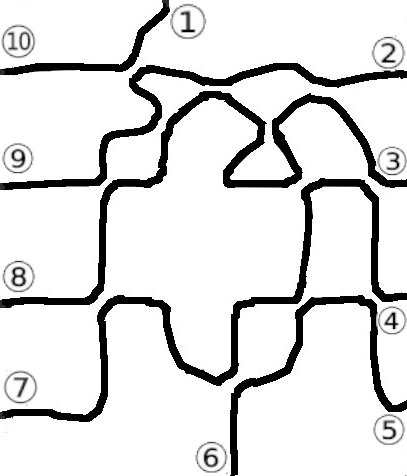}}
  \caption{Illustrating the existence of a maximal $\e$-smoothing. We
    start with a t-diagram of type $(1,4,1,4)$ (left). Then all the
    crossings are replaced by local smoothings. The third diagram from the left is a maximal $\e$-smoothing (all endpoints are preserved). So is the right-most one which is also loopless (connected
    components disconnected from the boundary were removed). This figure happens also to illustrate
    Lemma~\ref{lemma:SoneStwo}. In this case $k=5$, and the circled
    numbers in the right-most diagram correspond to the clockwise
    enumeration of $s_1,\dots,s_{10}$, and the diagram represents
    $\Sigma_2$.  For example we can see that $s_3$ is connected to
    $s_{10-3+1}=s_{8}$, as required by Lemma~\ref{lemma:SoneStwo}.}
  \label{fig:Smoothings}
\end{figure}

\begin{Remark}\label{rem:MaximalSmoothingAlwaysExists}
  It is standard to construct a maximal $\e$-smoothing of a tangle
  diagram for sufficiently small $\e$. This process is illustrated in
  Figure~\ref{fig:Smoothings}.
\end{Remark}

Clearly, the relations $\subset^\e$ and $\subset^\e_{\max}$
are transitive for any fixed $\e$.
The relation $S'\subset^\e_{\max} S$ does not
imply that $S'$ is an $\e$-smoothing of $S$,
because an $\e$-smoothing has no crossings
while $S'$ can still have them. Nonetheless,
there is a transitivity-like relationship between the two notions:

\begin{Lemma}\label{lemma:SmoothingTrans}
  Suppose $S'\subset^\e_{\max} S$ and $S^*$ is a maximal $\e$-smoothing of $S'$.
  Then $S^*$ is a maximal $\e$-smoothing of~$S$.
\end{Lemma}
\begin{proof}
  It is easy to check from the definitions that $\Cross(S')\subset\Cross(S)$.
  Let $C'\subset \Cross(S')$ be such that it witnesses $S^*\subset^\e S'$.
  Let $C$ witness $S'\subset^\e S$. Then $C\cup C'$ witnesses $S^*\subset^\e S$. By transitivity of $\subset^\e_{\max}$, also $S'\subset^\e_{\max} S$.
\end{proof}


\begin{Lemma}\label{lemma:MINofBridge}
  If 
  $\beta$ is an $m$-bridge, then $\MIN_{\vert}(\beta)=m$.
\end{Lemma}
\begin{proof}
  If $J$ is a vertical splitting curve, then by the Jordan Curve
  Theorem every strand of $\beta$ must intersect $J$ at least once,
  so $\MIN_{\vert}(\beta)\ge m$. The other direction follows
  from Lemma~\ref{lemma:UpperBoundOnMIN}.
\end{proof}


\begin{Lemma}\label{lemma:Hex}
  Suppose $T$ is a trivial t-diagram. Then exactly one of the
  following holds: (1) there is a vertical splitting curve $J$ with
  $C\cap J=\es$ or (2) there is a strand $s$ of $T$ which connects
  $V^-$ to $V^+$.
\end{Lemma}
\begin{proof}
  This proof is based on \cite{gale1979hex}. By tilting the square a
  bit to a rhombus it can be filled with a hexagonal grid so that it
  constitutes a board for the game of Hex. Assume that the hexagonal
  grid is so fine that for every hexagon $\sigma$ the set $|T|\cap \St(\sigma)$
  is path-connected, where $\St(\sigma)$ is the union of all the hexagons
  that touch $\sigma$, including $\sigma$ itself.  This is possible because $T$
  is piecewise linear and therefore locally path connected. Use
  compactness to find a uniform bound on the size of the
  hexagons. Note that since $T$ is trivial, $|T|\cap \St(\sigma)$ is in
  fact a subarc of some strand of~$T$.  Color the hexagons black if
  they intersect $T$ and white if they do not.  By the Hex
  Theorem~\cite{gale1979hex} there is a white path of hexagons from
  top to bottom (from $H^+$ to $H^-$) if and only if there is no black
  path from left to right (from $V^-$ to $V^+$).  In the former case,
  we find a vertical splitting curve disjoint from $T$ by following
  the hexagonal white path from top to bottom. In the latter case let
  $\sigma_1,\dots,\sigma_n$ be the sequence of black hexagons from $V^-$ to
  $V^+$ so that $\sigma_i$ touches $\sigma_{i+1}$ for $i\in \dint{1,n-1}$.  For
  each $i\in \dint{1,n}$, let $x_i\in \sigma_i\cap |T|$.  By the
  assumption, for each $i$ there is a subarc of a strand in $T$
  connecting $x_{i-1}$, $x_i$, and $x_{i+1}$. By the triviality
  of $T$, all these arcs
  must be subarcs of the same strand.  The union of these subarcs
  connects $V^-$ to $V^+$ so it must, in fact, be equal to that
  strand.
\end{proof}

\begin{Lemma}\label{lemma:SequenceOfSmoothingsAndBridge}
  Suppose $S=(S,N,M,\ell)$ is a t-diagram of type $(0,l_1,0,l_2)$
  for some $l_1,l_2\in\N$, and $\MIN_{\vert}(S)=m$. Then there is a
  positive $\e_0$ such that for all positive $\e<\e_0$ there is an
  $m$-bridge $\beta$ with $\beta\subset^\e S$. In fact, there is a
  maximal 
  $\e$-smoothing $S^*$ of $S$ such that $\beta \subset S^*$,
  i.e. every strand of $\beta$ is a strand of $S^*$.
\end{Lemma}
\begin{proof}
  If $m=0$ then let $\beta=\es$ be the $0$-bridge. Clearly
  $\beta\subset^\e S$ for all $\e$ and $\beta\subset S^*$ for all
  smoothings $S^*$ of $S$, so by
  Remark~\ref{rem:MaximalSmoothingAlwaysExists} we are done. Since
  $T$ is of type $(0,l_1,0,l_2)$, the edge $H^+$ is outside of the
  tangle, i.e. $H^+\cap S=\es$.  Let $A$ be the connected component
  of $\square\setminus S$ which contains~$H^+$, and let
  $C=\bar A\cap S$. For later use, note that it follows that
  \begin{equation}
    \label{eq:CinBdA}
    C\subset \partial A\text{ and }C\cap A=\es.
  \end{equation}
  We will prove the lemma by induction on $m$ and each step is
  proved by induction on the cardinality of the set $\Cross(S)\cap C$.
  So we will prove the following by induction on
  $\Card(\Cross(S)\cap C)$:
  \begin{ClaimL}\label{claim:Smoothing}
    Assume that the statement of the lemma holds for $m=p$ and that
    $\MIN_{\vert}(S)=p+1$. Then there is $\e_0>0$ such that for all
    positive $\e<\e_0$ there is an $\e$-smoothing $S^*$ of $S$ and a
    $(p+1)$-bridge $\beta\subset S^*$.
  \end{ClaimL}
  \begin{proof}[Proof of Claim~\ref{claim:Smoothing}]
    Suppose first that $\Cross(S)\cap C=\es$.  Let us show that then
    $C$ is a union of strands of $S$. For a contradiction assume that
    there is some strand $s$ of $S$ which is only partially in $C$ and
    let $x\in C\cap s$ and $y\in s\setminus C$. Take a path $\gamma$
    from $H^+$ to $x$ which touches $S$ only at $x$.  Let
    $\eta\subset s$ be a path from $x$ to $y$ in $s$ connecting $x$ to
    $y$. If $\eta$ passes through a crossing, let $z$ be the first
    crossing $\eta$ encounters on the way from $x$ to $y$.  Obtain
    $\eta'$ by shifting $\eta$ slightly in the tubular neighborhood of
    $s$ in the direction toward $\gamma$, i.e. so that
    $\eta'\cap \gamma$ is non-empty.  Obtain $\eta''$ by modifying
    $\eta'$ slightly at its other endpoint so that $\eta'$ intersects
    $S$ for the first time at $z$. Then here is a path inside
    $|\gamma|\cup |\eta''|$ which goes from $H^+$ to $z$ without
    intersecting $S$ on the way. By definition of $C$ this means that
    $z\in C$, a contradiction.  Therefore there is no crossing of $S$
    on $\eta$. Now shift $\eta$ slightly, again in the direction
    towards $\gamma$, to obtain $\eta'$ so that there is a path in
    $|\gamma|\cup |\eta'|$ which connects $H^+$ to $y$, a
    contradiction again with the assumption that $y\notin C$.  Thus,
    we have that $C$ is in fact a trivial t-diagram consisting of
    strands of~$S$.

    We now claim that there is a strand $s$ of $C$ such that $s$
    connects $V^-$ to $V^+$. If not, then by Lemma~\ref{lemma:Hex}
    there is a vertical splitting curve $J_v$ in the complement of
    $C$. By the definition of $C$, $J_v$ cannot intersect $S$.  But
    then $\MIN_{\vert}(S)=0$, a contradiction (we assumed
    $\MIN_{\vert}(S)=p+1>p\ge 0$).
  
    Let $S_1$ be the t-diagram which consists of all other strands of $S$
    than $s$. Then $S_1\cap s=\es$, for otherwise we would find a
    crossing in~$C$.  We claim that $\MIN_{\vert}(S_1)=p$.  Let
    $J\colon I\to \square$ be a vertical splitting curve with
    $\intersect(S,J)=p+1$. Then at least one of the intersections must
    be on $s$ by the Jordan Curve Theorem, so $\intersect(S_1,J)\le p$
    witnessing $\MIN_{\vert}(S_1)\le p$. Suppose that $J'$ is an
    arbitrary vertical splitting curve. Let $t\in I$ be the maximum of
    the set $(J')^{-1}[s]$ (the latter set is finite by general
    position and non-empty by the same argument as above, so it has a
    maximum). Since $J'(t)$ is in $C$, there is a curve $\gamma$ from
    $H^+$ to $J'(t)$ which intersects $S$ for the first time at
    $J'(t)$. Let $J''=\gamma\cat (J'\rest [t,1])$.  Then $J''$
    intersects $s$ exactly once and $S_1$ at most as many
    times as $J'$, so
    $\intersect(S_1,J')\ge \intersect(S_1,J'')=\intersect(S,J'')-1\ge
    \MIN_{\vert}(S)-1=p$.  By the arbitrary choice of $J'$ this means
    that $\MIN_{\vert}(S_1)\ge p$. We have proved that
    $\MIN_{\vert}(S_1)=p$.

    By the induction hypothesis there is an $\e_0>0$ such that for all
    $\e<\e_0$ there is an $\e$-smoothing $S^*$ of $S_1$ and
    a $p$-bridge $\beta$ such that
    $\beta\subset S^*$. Let $\e_0'<\e_0$ be such that an
    $2\e'_0$-neighborhood of $s$ does not contain any crossing
    of~$S_1$ (which are the same as crossings of $S$). Let $\e<\e'_0$
    and let $\beta$ be a $p$-bridge with $\beta\subset^\e S_1$.  Let
    $\beta'=\beta\cup s$. Then $\beta'$ is a bridge of type
    $(0,p+1,0,p+1)$, and $\beta'\subset^\e S$.

\begin{figure}
  \centering
  \includegraphics[width=0.9\textwidth]{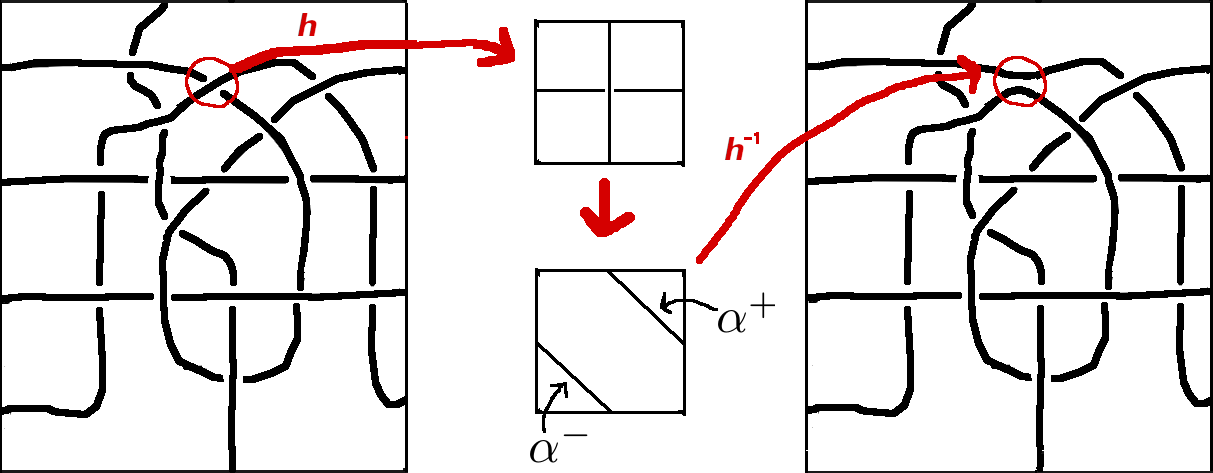}
  \caption{Illustration of the local smoothing performed in the proof of
  Claim~\ref{claim:Smoothing} of Lemma~\ref{lemma:SequenceOfSmoothingsAndBridge}.}
  \label{fig:LocalSmoothing}
\end{figure}

    We have proved Claim~\ref{claim:Smoothing} in the case where $C$
    contains no crossings of $S$.  Suppose we have proved this for
    the case when $C$ contains $q$ crossings of $S$, and suppose that
    $C$ contains $q+1$ crossings. Fix a vertical splitting curve $J_0$
    witnessing $\MIN_{\vert}(S)=p+1$.  Let $x\in C\cap\Cross(T)$ and
    let $U$ be a neighborhood of $x$ witnessing the transversality
    (see Definition~\ref{def:BasicImmersion}). Select $U$ to be also
    so small that $U\cap |J_0|=\es$. This is possible, because by the
    assumption of general position, $J_0$ does not contain $x$ and
    $J_0$ is compact, so $d(x,J_0)>0$. Further we can assume that
    $U\cap \partial\square=\es$. Let $\e_0$ be so small that
    for all $\e<\e_0$ we can choose $U=B(x,\e)$ and also smaller
    than any $\e_0$ occurred in the previous induction steps.
    By the assumption that $U$ witnesses
    transversality, there is a homeomorphism $h\colon U\to [-1,1]^2$
    such that
    $h[C]\subset ([-1,1]\times \{0\})\cup (\{0\}\times
    [-1,1])$. Consider $h[A\cap U]$ where $A$ is as defined in the
    beginning of the proof of the lemma.  By \eqref{eq:CinBdA} there
    is a point $y\in h[A\cap U]$ such that $y$ is not in any of the
    coordinate axes. Without loss of generality we may assume that $y$
    is in the positive-positive quarter, i.e. $y\in \openint{0,1}^2$.
    Then we have
    $$(\{0\}\times [0,1])\cup ([0,1]\times \{0\})\subset hC$$
    Let $S'$ be obtained from $S$ by replacing $S\cap U$ with
    $h^{-1}[\alpha^+\cup \alpha^-]$ where $\a^+=[(0,1),(1,0)]$ and
    $\a^-=[(-1,0),(0,-1)]$, see Figure~\ref{fig:LocalSmoothing}. Note that one has to
    reparametrize the immersion of $S'$ in an obvious way to make it a
    t-diagram (see Section~\ref{ssec:Reparam}).  We did not remove any endpoints in the process,
    so we have $S'\subset^{\e}_{\max} S$.  Let $A'$ and $C'$ be the sets
    defined for $S'$ in the same way as $A$ and $C$ were defined for
    $S$.  It remains to show that the induction hypothesis (for $q$)
    can be applied to $S'$ and $C'$. Once we do that, for any
    $\e<\e_0$ we will obtain a maximal $\e$-smoothing of $S'$ which
    contains a $(p+1)$-bridge. Since a maximal $\e$-smoothing of $S'$
    is a maximal smoothing of $S$ (Lemma~\ref{lemma:SmoothingTrans}),
    it will conclude the proof. Thus, we need to show that
    $\Cross(S')\cap C'$ has a smaller cardinality than
    $\Cross(S)\cap C$ and that $\MIN_{\vert}(S')=p+1$.  First let us
    prove some preliminary claims.
    \begin{ClaimL}\label{claim:Ap}
      $\left(\,\openint{-1,0}\times \openint{0,1}\,\right)\cap hA=\left(\,\openint{0,1}\times
      \openint{-1,0}\,\right)\cap hA=\es$
    \end{ClaimL}
    (Recall that $\openint{a,b}$ denotes the open interval from $a$ to~$b$.)
    \begin{proof}[Proof of Claim \ref{claim:Ap}]
      We will show that
      $\openint{-1,0}\times \openint{0,1}\cap hA=\es$ and the other
      part follows by symmetry. Suppose for a contradiction that there
      is $z\in \openint{-1,0}\times \openint{0,1}\cap hA$ and let
      $z'=h^{-1}(z)$.  Let $y$ be as above, and $y'=h^{-1}(y)$.  Let
      $\gamma_z$ and $\gamma_y$ be paths connecting $z'$ and $y'$
      respectively to $H^+$ and let $\f$ be the line segment
      connecting the endpoints of $\gamma_z$ and $\gamma_y$ in $H^+$.
      Let $\psi$ be the line segment connecting $z$ and $y$.  If
      $\gamma_z$ and $\gamma_y$ intersect, modify them by perform
      smoothing operations within $A$ on their crossings to make them
      disjoint.  Then $h^{-1}[\psi]\cup \gamma_1\cup \f\cup \gamma_2$
      is a Jordan curve and it is entirely in $\bar A$. Also, the part
      $\gamma_1\cup \f\cup \gamma_2$ does not intersect $S$.  But the
      remaining part, $h^{-1}[\psi]$ intersects $S$ only in one point.
      The strand of $S$ which enters the Jordan region must end in
      $\partial \square$, so this is a contradiction with the Jordan
      Curve Theorem.
    \end{proof}
    \begin{ClaimL}\label{claim:CpSubsetC}
      $C'\setminus U\subset C\setminus U$
    \end{ClaimL}
    \begin{proof}[Proof of Claim~\ref{claim:CpSubsetC}]
      Suppose $z\in C'\setminus U$ and suppose $\gamma$ is a curve
      connecting $H^+$ to $z$.  If $\gamma$ lies outside of $U$, then
      it witnesses that $z\in C$. Otherwise, w.l.o.g. assume general
      position between $\gamma$ and $\partial U$, so $\gamma^{-1}[U]$
      is a finite union of closed intervals. Let $[t_0,t_1]$ be the
      first one of them (i.e. such that $t_0$ is the smallest point
      with $\gamma(t_0)\in U$). Then $\gamma\rest[0,t_0]$ witnesses
      that $\gamma(t_0)\in A$, so by Claim~\ref{claim:Ap}
      $h(\gamma(t_0))$ is either
      in $([0,1]\times \{1\})\cup (\{1\}\times [0,1])$ or in
      in $([-1,0]\times \{-1\})\cup (\{-1\}\times [-1,0])$.
      By the Jordan Curve Theorem $t_1$ must be in the same of those
      sets, so by a local ambient isotopy it is possible to remove $\gamma\rest [t_0,t_1]$
      from $U$ thereby reducing the number of components of $\gamma^{-1}[U]$
      by one. Inductively we can remove all of them,
      so we reduce the argument to the one where $\gamma$ is outside of~$U$.
    \end{proof}
    By Claim~\ref{claim:CpSubsetC} and the fact that
    $\Cross(S')\setminus U=\Cross(S')\setminus \{x\}=\Cross(S)\setminus \{x\}$
    we have
    $$\Cross(S')\cap C'\setminus U\subset\Cross(S)\cap C\setminus \{x\}$$
    so $x\in (\Cross(S)\cap C)\setminus(\Cross(S')\cap C')$.  Now let
    us show that $\MIN_{\vert}(S')=p+1$.
    \begin{ClaimL}\label{claim:ApsubA}
      $A'\subset A$.
    \end{ClaimL}
    \begin{proof}[Proof of Claim \ref{claim:ApsubA}]
      Use the same argument as in the proof of
      Claim~\ref{claim:CpSubsetC} to show that
      $A'\setminus U\subset A\setminus U$.  Suppose $z\in A'\cap
      U$. Then $h(z)$ cannot be between the lines $\a^+$ and
      $\a^-$ (Figure~\ref{fig:LocalSmoothing}) because then we would be able to find a point
      in $hA$ which contradicts Claim \ref{claim:Ap}.
      So $h(z)\in [0,1]^2$ or $h(z)\in
      [-1,0]^2$. Suppose $h(z)\in [0,1]^2$ and let $\gamma$ be a path
      connecting $z$ to $H^+$. By the same argument as above, $h\left[\gamma\cap U\right]$
      cannot pass between $\alpha^+$ and $\alpha^-$ and also
      $\gamma\cap S'=\es$.  But this implies that $\gamma$ in fact
      witnesses that $z\in A$. The case $h(z)\in [-1,0]^2$ follows in
      the same way.
    \end{proof}
    Recall that $J_0$ was defined to witness that
    $\MIN_{\vert}(S)= p+1$ and lies outside $U$, so
    also witnesses $\MIN_{\vert}(S')\le p+1$.
    Thus, it is enough to show that
    $\MIN_{\vert}(S')\ge p+1$. Suppose $J$ is a vertical splitting
    curve such that $\intersect(S',J)=\MIN_{\vert}(S')$. If $J$ does
    not intersect $U$, then it is a splitting curve for $S$, so
    $\intersect(S',J)=\intersect(S,J)\ge p+1$. Otherwise, there is
    $t$ with $J(t)\in U$.  Let $t$ be the maximum such
    $t$. W.l.o.g. by a general position argument, $J(t)\notin |S'|$.
    If $J(t)\in A'$, let $\gamma$ be a curve connecting $H^+$ to
    $J(t)$ which lies in $A'$. By Claim~\ref{claim:ApsubA} it lies in
    $A$, and so outside $S$. Let $J'=\gamma\cat (J\rest [t,1])$.
    $$\intersect(S',J)\ge \intersect(S',J')=\intersect(S,J')\ge\MIN_{\vert}(S)=p+1.$$
    Suppose that $J(t)$ is in the middle area, i.e.  $h(J(t))$ is
    between the lines $\a^+$ and $\a^-$. Then by Claims \ref{claim:Ap}
    and \ref{claim:ApsubA}, $J(t)\notin A'$. Thus, there must be at
    least one $t_0<t$ such that $J(t_0)\in C'\subset |S'|$,
    so $\intersect(S',J\rest [t_0,1])+1\le \intersect(S',J)$
    On the other hand there it is not hard to find $\gamma$  which connects $H^+$ to $J(t_0)$
    and has exactly one intersection with $S'$ and exactly one intersection with $S$.
    Now
    \begin{equation}
      \intersect(S',J)\ge 1+\intersect(S'J\rest [t,1])
      \ge \intersect(S',\gamma)+\intersect(S',J\rest [t,1])
      =   \intersect(S',\gamma\cat (J\rest [t,1])).\label{eq:plusone}    
    \end{equation}
    Since $t$ is the maximum point in $J^{-1}(t)$ and by the choice of $\gamma$, we have
    $$=\intersect(S,\gamma\cat (J\rest [T,1]))\ge\MIN_{\vert}(S)=p+1.$$
    Suppose now that $t$ is such that $h(J(t))$ is below $\a^-$ and
    that area is not in $A'$. But then the component of $J^{-1}[U]$ to which
    $t$ belongs must be either contained entirely to that area, in which case we can
    remove it by a local ambient isotopy
    out from $U$ and reduce the number of components in $J^{-1}[U]$
    (similarly as we did in the proof of Claim~\ref{claim:CpSubsetC}),
    or it intersects $\a^-$ in which case we apply a similar argument as in
    \eqref{eq:plusone}
    replacing ``$1+$'' with ``$2+$''.
  \end{proof}
  This completes the proof of the lemma.
\end{proof}

\begin{Lemma}\label{lemma:Ifk-1Thenk}
  Suppose $T$ is a t-diagram of type $(0,k,0,k)$ and
  $\MIN_{\vert}(T)\ge k-1$. Then $\MIN_{\vert}(T)\ge k$.
\end{Lemma}
\begin{proof}
  Let $J$ be a vertical splitting curve which has the minimum
  intersection number with~$T$. Assume for a contradiction
  that $\intersect(T,J)=k-1$.
  Let $\e$ be so small that
  $$J\cap \Cup_{x\in\Cross(T)}B(x,\e)=\es$$
  and so that the conclusion of
  Lemma~\ref{lemma:SequenceOfSmoothingsAndBridge} holds for $T$.  Then
  there is a maximal $\e$-smoothing $T^*$ of $T$ and a
  $(k-1)$-bridge $\beta\subset T^*$ which has $k-1$ strands. But there
  is one more endpoint on $V^-$ and one on $V^+$ which are not
  connected by strands of $\beta$. So there must be a strand $s$ in
  $T^*$ connecting them. Let $\beta'=\beta\cup s$. Then $\beta'$ is a
  $k$-bridge. Since $J$ goes outside of a neighborhood where $T$ has
  been modified,
  $\intersect(T,J)=\intersect(T^*,J)\ge\intersect(\beta',J)\ge\MIN_{\vert}(\beta')=
  k$ where the last equality is by Lemma~\ref{lemma:MINofBridge}. A
  contradiction.
\end{proof}

\begin{figure}
  \centering
  \fbox{\includegraphics[width=0.3\textwidth]{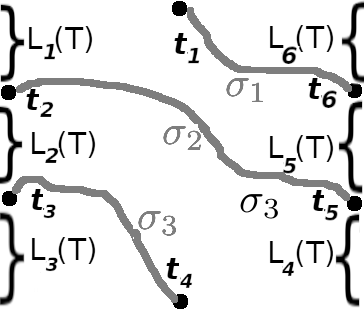}}
  \caption{Illustration of Lemma \ref{lemma:SoneStwo} and the proof of
    Lemma \ref{lem:ConnectingLemma1}. The gray lines are strands of
    $\Sigma_1$. In this example $k=3$ and $n=6$. It is evident from
    the picture that if, say $i=3$ and $j=5$, then any splitting curve
    going from $L_i(T)$ to $L_j(T)$ must intersect at least
    $|f(j)-i|+1=|6-5+1-3|+1=2$ strands of~$\Sigma_1$.}
  \label{fig:IllustrlemmaSigma}
\end{figure}

Recall the function $f(j)=n-j+1$ from Section~\ref{ssec:Matrices}.
\begin{Lemma}\label{lemma:SoneStwo}
  Let $k\in\N_+$ and $n=2k$.  Let $T$ be a t-diagram of type
  $(1,k-1,1,k-1)$ such that $\MIN_{\vert}(T)=k-1$. Let
  $t_1,\dots,t_n$ be the endpoints of $T$ enumerated counterclockwise
  starting from the only point on $H^+$. Let $s_1,\dots,s_n$ be the
  same points enumerated clockwise starting from the same point,
  i.e. $s_1=t_{1}$ and $s_i=t_{n-i+2}$ for all $i\in \dint{2,n}$.
  Then there is $\e_0$ such that for all $\e<\e_0$ there are two
  trivial t-diagrams $\Sigma_1$ and $\Sigma_2$ satisfying the
  following conditions.
  \begin{enumerate}
  \item The $i$-th strand of $\Sigma_1$ connects $t_i$ to $t_{f(i)}$ for all $i\in \dint{1,k}$
  \item The $i$-th strand of $\Sigma_2$ connects $s_i$ to $s_{f(i)}$ for all $i\in \dint{1,k}$
  \item $\Sigma_1,\Sigma_2\subset^\e_{\max} T$.
  \end{enumerate}  
\end{Lemma}
\begin{Remark}
  Figures \ref{fig:Smoothings} and \ref{fig:IllustrlemmaSigma}
  illustrate this lemma.
\end{Remark}
\begin{proof}
  Clearly proving the existence of $\Sigma_1$ implies the existence of
  $\Sigma_2$ by mirroring $T$ along the middle vertical axis.  Note
  that by the definition of type $(1,k-1,1,k-1)$ and the
  counterclockwise enumeration of $t_i$, $t_1$ is the only point on
  $H^+$ and $t_{k+1}$ is the only point on $H^-$.  There is an
  orientation-preserving self-homeomorphism $h$ of $\square$ taking
  $t_1$ to $V^-$, $t_{k+1}$ to $V^+$, fixing all other points $t_i$ as
  well as the top-right corner and lower-left corners of $\square$.
  In particular
  \begin{equation}
    h^{-1}H^{\pm}\subset H^{\pm}.\label{eq:HpmHpm}  
  \end{equation}
  For any vertical splitting curve $J$ we have
  $\intersect(hT,J)=\intersect(T,h^{-1}J).$
  If $J$ is a vertical splitting curve, then so is $h^{-1}J$, by
  \eqref{eq:HpmHpm}.  Thus, $\MIN_{\vert}(hT)\ge \MIN_{\vert}(T)=k-1$
  and $hT$ is of type $(0,k,0,k)$. By Lemmas \ref{lemma:Ifk-1Thenk} 
  and \ref{lemma:UpperBoundOnMIN} we
  have $\MIN_{\vert}(hT)=k$. Let $\beta$ be the $k$-bridge given by
  Lemma \ref{lemma:SequenceOfSmoothingsAndBridge}.  Then
  $\Sigma_1=h^{-1}\beta$ is as desired.
\end{proof}

For a sequence $a_1,\dots,a_l\in\partial\square$
with $l\ge 3$ denote by $\cordop{a_1,\dots,a_l}$ the statement that
these points are arranged counterclockwise.

\begin{Fact}\label{fact:MustIntersect}
  For splitting curves $J$ and $\sigma$ with $\cordop{J(0),\sigma(0),J(1),\sigma(1)}$
  we have $\intersect(\sigma,J)\ge 1$.
\end{Fact}
\begin{proof}
  Direct application of \cite[Theorem III.3.A]{Bi83}.
\end{proof}

We are now ready to prove Lemma~\ref{lem:ConnectingLemma1}.

\begin{proof}[Proof of Lemma~\ref{lem:ConnectingLemma1}]
  Let $J$ be an $ij$-splitting curve for $T$. W.l.o.g. $i\le j$.
  If $i=j$, $w(i,j)=0$, so the statement is trivial. So we are left
  with the case $i<j$.
  Suppose first that
  $i<j\le k$.  Since $T$ is of type $(1,k-1,1,k-1)$, $L_i(T)$ and
  $L_j(T)$ are now subsets of $V^-$.  We have to show that
  $\intersect(T,J)\ge w(i,j)=|j-i|=j-i$. Suppose not. Let $J_1$ be the
  line from the top endpoint of $V^-$, $(0,1)$, to $J(0)$. Then $\intersect(T,J_1)=i-1$.  Let
  $J_2$ be the line from $J(1)$ to the lower endpoint of $V^-$, $(0,0)$. Then
  $\intersect(T,J_2)=k-j$.  Let $J_*=J_1\cat J\cat J_2$. Then
  $\intersect(T,J_*)=\intersect(T,J_1)+\intersect(T,J)+\intersect(T,J_2)=i-1+\intersect(T,J)+
  k-j< i-1+j-i+k-j= k-1.$ By pushing $J_*$ slightly away from $V^-$
  into $\square$, obtain a vertical splitting curve $J_*'$ with
  $\intersect(T,J_*')<k-1$, a contradiction.  Same argument for
  $i>j>k$ noting that then $L_i(T),L_j(T)\subset V^+$. Suppose
  \begin{equation}
    i\le k<j.\label{eq:ikj}  
  \end{equation}
  In this case $i\in V^-$ and $j\in V^+$.  
  We consider separately the cases $f(j)\le i$ and $f(j)>i$.  Suppose
  first that 
  \begin{equation}
    \label{eq:fjgei}
    f(j)\le i
  \end{equation}
  Figure~\ref{fig:IllustrlemmaSigma} illustrates the forthcoming
  technicalities. Let $\Sigma_1\subset^\e_{\max} T$ be as given by
  Lemma~\ref{lemma:SoneStwo}. Let $\sigma_1,\dots,\sigma_k$ be the
  strands of $\Sigma_1$ such that $\sigma_l$ connects $t_l$ to
  $t_{f(l)}$, i.e., we can think of $\sigma_l$ as
  a splitting curve $\sigma_l\colon I\to \square$ with
  \begin{equation}
    \sigma_l(0)=t_l\quad\text{ and }\quad \sigma_l(1)=t_{f(l)}\label{eq:sigmatl}  
  \end{equation}
  where $t_1,\dots,t_n$ are also as given in the formulation of
  Lemma~\ref{lemma:SoneStwo}. Note that
  \begin{equation}
    \cordop{t_1,\dots,t_n}.\label{eq:cordopt}
  \end{equation}
  By assuming that $\e$ is sufficiently small and using the fact that
  $\Sigma_1\subset^\e T$, by Lemma~\ref{lemma:CrossingsWithSmoothings}
  we can assume that $\intersect(T,J)\ge \intersect(\Sigma_1,J)$, so
  it is enough to bound the latter.  Let
  $$A=\{l\in\dint{1,k}\mid \cordop{\sigma_l(0),J(0),\sigma_l(1),J(1)}\}$$
  Then by Fact~\ref{fact:MustIntersect}
  $$\intersect(\Sigma_1,J)\ge \sum_{l\in A}\intersect(\sigma_l,J)\ge \Card(A).$$
  So it is enough to compute a lower bound for $\Card(A)$. By \eqref{eq:sigmatl},
  $$A=\{l\in\dint{1,k}\mid \cordop{t_l,J(0),t_{f(l)},J(1)}\}$$
  By definition of $L_i(T)$ and $L_j(T)$,
  we have that
  $\cordop{t_i,J(0),t_{i+1}}$ and
  $\cordop{t_j,J(1),t_{j+1}}$, so
  \begin{equation}
    \cordop{t_l,J(0),t_{f(l)},J(1)}\label{eq:cordopJ}  
  \end{equation}
  holds if $l\le i$ and $f(l)\le j$. The latter is equivalent to
  $f(j)\le l$ by~\eqref{eq:involution}.  Thus, \eqref{eq:cordopJ}
  holds if $l\in \dint{1,i}\cap \dint{f(j),k}=\dint{f(j),i}$, thus
  \begin{equation*}
    \Card(A)\ge \Card(\dint{f(j),i})\stackrel{\eqref{eq:dintcard},\eqref{eq:fjgei}}{=}|f(j)-i|+1\stackrel{\eqref{def:w2}}{=}w(i,j).
  \end{equation*}
  For the case $f(j)>i$ use $\Sigma_2$ and $s_1,\dots,s_n$ in a symmetric fashion.
\end{proof}

\subsection{Bounding the Number of Crossings Between Overlapping Tangle Diagrams}

Fix $k\in\N_+$ and $n=2k$. Next we will define the matrix $e(T,S)$
for particular t-diagrams $T$ and $S$. Viewed as a multigraph,
the vertices of $e(T,S)$ are the sets $L_1(T),\dots,L_n(T)$ and the
edges are the strands of~$S^*$ where $S^*$ is as
given by Lemma \ref{lemma:SequenceOfSmoothingsAndBridge} for~$S$.

\begin{Def}\label{def:eTS}
  Suppose that $T$ is a t-diagram of type $(1,k-1,1,k-1)$ and $S$ a
  t-diagram of type $(0,l_1,0,l_2)$ for some $l_1,l_2\in\N$. We will
  define the $(n\times n)$-matrix $e(T,S)$ as follows.  Let $\e_0$ be
  as given by Lemma~\ref{lemma:SequenceOfSmoothingsAndBridge} for $S$
  and let $\e<\e_0$ so small that for any $x,y\in \Cross(T\times S)$,
  $B(x,\e)$ and $B(y,\e)$ are disjoint and witness the transversality
  of the respective crossing. Let $S^*$ be an $\e$-smoothing $S$ such
  that $S^*$ is a trivial t-diagram and there is a trivial bridge of
  type $(0,m,0,m)$ as a subset of $S^*$ where $m=\MIN_{\vert}(S)$.
  Such smoothing exists by
  Lemma~\ref{lemma:SequenceOfSmoothingsAndBridge} and the choice
  of~$\e$ above. For any $(i,j)\in \dint{1,n}^2$, if $i\ne j$, define
  $e(T,S)(i,j)$ to be the number of strands of $S^*$ that connect
  $L_i(T)$ to $L_j(T)$.  Define $e(T,S)(i,i)$ to be $2x$ where $x$ is
  the number of strands which start and end at~$L_i(T)$.
\end{Def}

\begin{Remark}
  A careful observation will reveal that $e(T,S)$ is not unique,
  because it depends on the choice of $S^*$. But this does not matter
  to us.  Any choice of $e(T,S)$ will work for our purposes as long as
  it is obtained from any such~$S^*$.
\end{Remark}

\begin{Prop}[Overlapping t-Diagrams]\label{prop:Hybrid}
  Let $l_1,l_2\in\N$. Suppose that $T$ and $S$ are t-diagrams of types
  $(1,k-1,1,k-1)$ and $(0,l_1,0,l_2)$ respectively. Suppose that
  $e(T,S)$ satisfies property $(*)$
  (Definition~\ref{def:PropertyStar}).  Then
  $\intersect(T,S)\ge \MIN_{\vert}(S)+k-1.$
\end{Prop}
\begin{proof}
  Let $\e$ and $S^*$ be as in
  Definition~\ref{def:eTS}. In particular, $S^*$ contains an
  $m$-bridge for $m=\MIN_{\vert}(S)$.
  By Lemma~\ref{lemma:CrossingsWithSmoothings},
  $\intersect(T,S)\ge \intersect(T,S^*)$.
  Since $S^*$ and $S$ have the same
  number of endpoints, they have the same number of strands which is
  equal to~$l=(l_1+l_2)/2$.
  Let $\mSS$ be the set of strands of $S^*$,
  $\mSS=\{s_1,\dots,s_l\}$. 
  Then 
  $$\intersect(T,S^*)\ge \sum_{s\in\mSS}\intersect(T,s).$$
  (The inequality may be strict, if there are loops in $S^*$ which are not
  counted in the sum.)
  
  For all $(i,j)\in \dint{1,n}^2$
  let $\mSS_{ij}\subset\mSS$ be the set of those strands
  which go from $L_i(T)$ to $L_j(T)$. Then $\{\mSS_{ij}\}_{1\le i<j\le n}$
  is a partition of a (possibly proper) subset of $\mSS$ (we are not counting strands in $\mSS_{ii}$). Therefore
  $$\intersect(T,S^*)\ge \sum_{i=1}^n\sum_{j=i+1}^n\sum_{s\in \mSS_{ij}}\intersect(T,s).$$
  For each $s\in \mSS_{ij}$ we have $\intersect(T,s)\ge v_T(i,j)$ by
  Definition \ref{def:MatrixV_T} of~$v_T$. On the other hand,
  $\Card(\mSS_{ij})=e(T,S)(i,j)$ by the definition of $e(T,S)$, so by
  Lemma~\ref{lem:ConnectingLemma1},
  $$\intersect(T,S^*)\ge \sum_{i=1}^n\sum_{j=i+1}^n e(T,S)v_T(i,j)\ge \sum_{i=1}^n\sum_{j=i+1}^n e(T,S)w(i,j)=\xi(e(T,S)).$$
  By the assumption that $e(T,S)$ satisfies property $(*)$ and by
  Lemma~\ref{lemma:KeyLemma} we now have
  \begin{equation}
      \intersect(T,S^*)\ge \mu(e(T,S))+k-1.\label{eq:intmuTS}
  \end{equation}
  So it remains to show that $\mu(e(T,S))\ge m=\MIN_{\vert}(S)$.
  By definition,
  $$\mu(e(T,S))=\sum_{i=1}^k\sum_{j=k+1}^n e(T,S)(i,j)=\sum_{i=1}^k\sum_{j=k+1}^n \Card(\mSS_{ij}).$$
  The latter quantity counts exactly the number of strands of $S^*$ which go
  from $V^-$ to $V^+$. Since $S^*$ contains an $m$-bridge, it is at least~$m$.
\end{proof}

\begin{figure}
  \centering
  \includegraphics[width=0.7\linewidth]{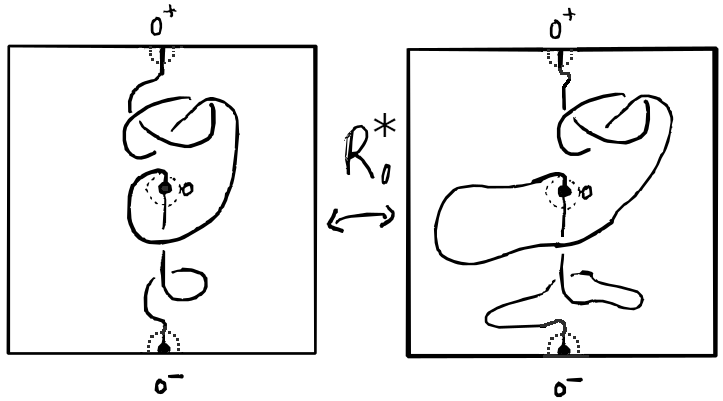}
  \caption{An $R^*_0$-move keeping a neighborhoods of $o$, $o^-$, and $o^+$ fixed.}
  \label{fig:SpecialRmoves}
\end{figure}

\section{Bounding the Number of Hybrid Crossings}
\label{sec:HybridCrossings}

As mentioned in the beginning of Section~\ref{sec:BoundingTangles},
in this section we will give a lower bound for $\intersect(\widehat K^-,\widehat K^+)$
which will be the missing piece for the proof of Theorem~\ref{thm:Main}.

\subsection{Subdivided Tangle Diagrams}
\label{sec:SubdividedTD}

In order to state the Fan Lemma \ref{lemma:Fan} in a sufficiently
general form so that it can later be used for links and not only for knots,
we need to extend the definitions of the partial knot diagrams $K^{\pm}$ to tangle diagrams
more generally.

\begin{Def}[Subdivided Tangle Diagram]\label{def:SubdividedTangleDiagram}
  A tuple $(T,T^+,T^-,N,M,\ell)$ is a \emph{subdivided t-diagram}
  if the following hold:
  \begin{enumerate}[label={\upshape ($PT_{\arabic*}$)}, leftmargin=3pc]
  \item $(T,N,M,\ell)$ is a t-diagram of type $(1,0,1,0)$ \label{def:PTD1}
  \item $N=\{1\}$ is a singleton
  \item There is a partition $M=M_+\cup M_-$ such that
    $T^+=T\rest \big((I^+\times \{1\}) \cup (S^1\times M_+)\big)$ and
    $T^-=T\rest \big((I^-\times \{1\}) \cup (S^1\times M_-)\big)$.
  \end{enumerate}
  A subdivided t-diagram $(T,T^+,T^-,\{1\},M,\ell)$ is
  \emph{centered} if the knot diagram defined by
  $T\rest (I\times \{1\})$ is centered.  The lower and upper minimal
  splitting numbers $\MIN_{\pm}(T^{\mp})$ are defined in the same way
  as for (partial) knot diagrams. Consistent with the terminology for
  knot diagrams, call $T^{\pm}$ the \emph{partial tangle diagrams}, or
  \emph{partial t-diagrams}. The $\sim^*_\RR$, 
  $\sim^*_0$, and $\sim_{\RR}$-equivalences are extended to subdivided
  t-diagrams in an obvious way. A subdivided t-diagram is \emph{permitted} if for
  all $x\in T^+\cap T^-$, $\ell^{-1}(x)\in (I^+\cup (S^1\times M_+))$.
   Extend the definitions of $hT$ and
  $T\capdot G$ (see Definition \ref{def:TcapG-tangle} and
  Section~\ref{sec:ApplyHomeo}) for partial and subdivided tangle
  diagrams in an obvious way, but note that we will treat $T\capdot G$ as
  a t-diagram, not as a subdivided one, because it might not be of type
  $(1,0,1,0)$ as condition \ref{def:PTD1} above requires.  Also, if
  $T^+$ is a partial t-diagram, then $T\capdot G$ is a t-diagram when
  $o\notin G$. All situations where we will use this notion will be
  ``non-problematic'' in this sense.
\end{Def}

\subsection{Upper and Lower Splitting Curves}

\begin{Def}\label{def:lowersplitting}
  A \emph{lower splitting curve} is an embedding
  $J^-\colon \left[0,\frac{1}{2}\right]\to \square^-$ such that
  $J^-(0)=o^-$ and $J^-\left(\frac{1}{2}\right)=o$. An \emph{upper splitting curve}
  is symmetrically an embedding $J^+\colon \left[\frac{1}{2},1\right]\to \square^+$ such that
  $J^+\left(\frac{1}{2}\right)=o$ and $J^+(1)=o^+$.
  Given a t-diagram $T$, define the \emph{lower splitting number of $T^+$} by
  $$\MIN_-(T^+)=\min\{\intersect(T^+,J^-)-1\mid J^-\text{ is a lower splitting curve}\}$$
  and symmetrically, define the
  \emph{upper splitting number of $T^-$} by
  $$\MIN_+(T^-)=\min\{\intersect(T^-,J^+)\mid J^+\text{ is an upper splitting curve}\}.$$
  As in the definition of $\MIN_{\ray}$, note that
  the point $o$ is not counted in $\intersect(T^\pm,J^\mp)$
\end{Def}

\begin{Lemma}\label{lemma:Raypm_ineq}
  For a subdivided t-diagram $T$, for all 
  $\pm,\mp\in\{+,-\}$,
  $\MIN_{\pm}(T^\mp)\ge \MIN_{\ray}(T^\mp)$.
\end{Lemma}
\begin{proof}
  All upper and lower splitting curves are rays.
\end{proof}

\subsection{Straightening Rays}

Let $v^-=\left(0,\tfrac{1}{2}\right)$ and
$v^+=\left(1,\tfrac{1}{2}\right)$ be midpoints of $V^-$ and $V^+$
respectively. Denote by
$H_{\tfrac{1}{2}}=[0,1]\times \left\{\tfrac{1}{2}\right\}=[v^-,v^+]$
the middle horizontal line segment.
Let $b_{\pm}(\delta)$ be the point $o\pm (\delta/2,0)$.
For $r>0$, let $\square(r)$ be the square with side-length $r$
and center $o$:
$$\square(r)=\left[\frac{1-r}{2},\frac{1+r}{2}\right]^2.$$
For instance, $\square=\square(1)$.

\begin{figure}
  \centering
  \begin{subfigure}[t]{0.40\textwidth}
    \centering
    \fbox{\includegraphics[width=\linewidth]{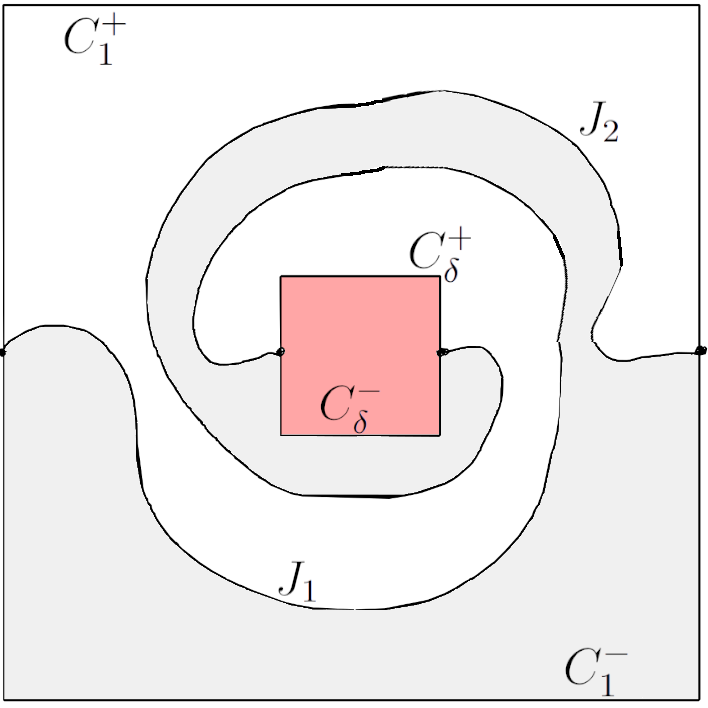}}
    \caption{
    }
  \end{subfigure}
  \hspace{10pt}
  \raisebox{20\height}{$\Rightarrow$}
  \hspace{10pt}
  \begin{subfigure}[t]{0.40\textwidth}
    \centering
    \fbox{\includegraphics[width=\linewidth]{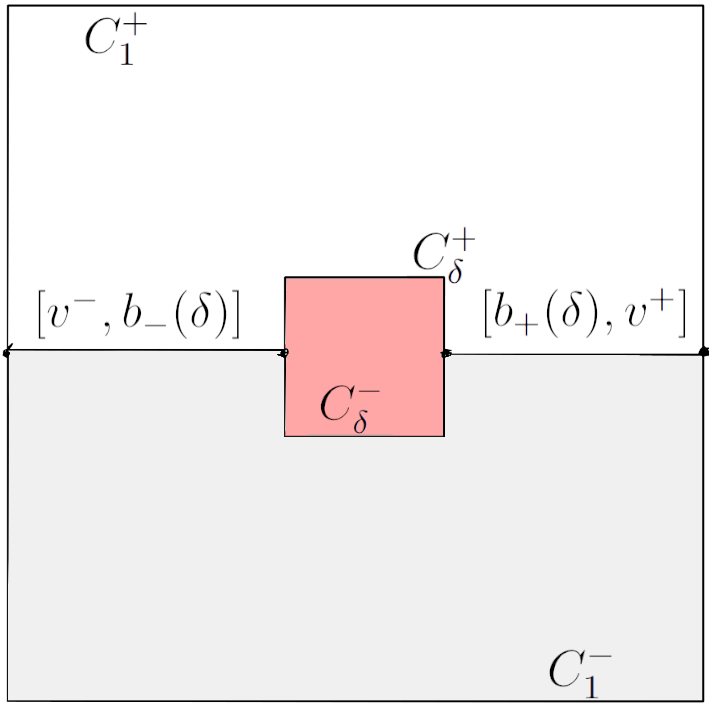}}
    \caption{}    
  \end{subfigure}

  \caption{Illustration of Lemma~\ref{lemma:StraightenTwoCurves}. The white
  area in (a) is bounded by the Jordan curve $C^+_1\cup J_1\cup C^+_\delta\cup J_2$. It is mapped by a canonical homeomorphism to the Jordan curve $C^+_1\cup [v^-,b_-(\delta)]\cup C^+_\delta\cup [b_+(\delta),v^+]$ in (b). This
  homeomorphism extends to the white area by the Jordan Curve Theorem. The light gray area in (a) is bounded
  by $C^-_1\cup J_1\cup C^-_\delta\cup J_2$ and is mapped to the corresponding area in (b). Together these homeomorphisms extend to a self-homeomorphism
  of the square which fixes the light red area (which is $\square(\delta)$ pointwise. 
  One way to visualize this homeomorphism is to imagine an isotopy which rotates the light red square
  $360^\circ$ clockwise.}
  \label{fig:JordanSq}
\end{figure}

\begin{Lemma}\label{lemma:StraightenTwoCurves}
  Suppose $0<\delta<\frac{1}{8}$.
  Suppose $J_1$, $J_2\colon I\to\square(1-\delta)$ are embeddings such that 
  $J_1\cap J_2=\es$,
  $J_1(0)=b_-(\delta)$, $J_1(1)=v^-$,
  $J_2(0)=b_+(\delta)$, $J_2(1)=v^+$, and for $t\in\ocirc{I}$, $J_l(t)\in \ocirc{\square}\setminus \square(\delta)$,
  $l\in \{1,2\}$. Then there is a homeomorphism $h\colon\square\to\square$ such that
  $|h\circ J_1|=[v^-,b_-(\delta)]$, $|h\circ J_2|=[b_+(\delta),v^+]$ and $h$ is identity
  on $\partial\square\cup \square(\delta)$.
\end{Lemma}
\begin{proof}
  For $r\le 1$, let $C^{\pm}_r$ be the upper and lower halves of
  $\partial\square(r)$,
  $C^{\pm}_r=\partial\square(r) \cap \square^{\pm}$.  Then
  $C^{\pm}_{\delta}$ are curves connecting $b_{-}(\delta)$ to
  $b_+(\delta)$ and $C^{\pm}_{1}$ are curves connecting $v^-$ to $v^+$.
  Now $C^+_{1}\cup J_1\cup C^+_\delta\cup J_2$ 
  and $C^-_1\cup J_1\cup C^-_\delta\cup J_2$ are
  Jordan curves. Let
  $f_1\colon |J_1|\to [v^-,b_-(\delta)]$ be given by
  $f_1(J_1(t))=(1-t)v^-+tb_-(\delta)$. Similarly let
  $f_2\colon |J_2|\to [b_+(\delta),v^+]$ be given by
  $f_2(J_2(t))=(1-t)b_+(\delta)+tv^+$. Denote by $\Gamma^\pm_r$ the
  identity map on $C^\pm_r$. Then
  $\Gamma^\pm_1\cup f_1\cup f_2\cup \Gamma^\pm_\delta$ is a homeomorphism
  between two Jordan curves, so it extends to a homeomorphism $h^\pm$
  between the bounded regions. Let $h_0=h^+\cup h^-$.
  This homeomorphism is identity on 
  $C^+_\delta\cup C^-_\delta=\partial\square(\delta)$ and on
  $C_1^+\cup C^-_1=\partial\square$. So  it extends to a
  self-homeomorphism $h$ of $\square$ which keeps $\partial\square\cup \square(\delta)$
  fixed pointwise. See Figure \ref{fig:JordanSq} for an illustration.
\end{proof}

Denote $d_{\pm}(\delta)=o\pm ((1-\delta)/2,0)=b_{\pm}(1-\delta)$.
The following follows by simple rescaling and the assumption $\delta<1/8$:
\begin{Cor}\label{cor:Rescale}
    Lemma~\ref{lemma:StraightenTwoCurves} applies with $\square$
    replaced by $\square(1-\delta)$ and $v^\pm$ replaced by $d_{\pm}(\delta)$. \qed
\end{Cor}

For $0<r<1$, let
$$U(r)=\ocirc{\square}\setminus \square(1-r).$$
Clearly $U(r)$ is open and its boundary is
$\partial \square(1-r)\cup\partial\square$.

\begin{Lemma}[Straightening]\label{lemma:StraightenJ}
  Suppose 
  that $J\colon I\to \square$ is a ray such that $J(1)=v^+$ and for
  some $0<\d<\frac{1}{8}$, 
  $$J\cap \square(2\delta)=[0,b_+(2\delta)]\quad\text{ and }\quad
  J\cap U(2\delta)=[d_+(2\delta),v^+].$$
  Then there is a boundary-fixing self-homeomorphism $h$ of $\square$
  such that $|h\circ J|=[o,v^+]$ and $h$ is the identity on
  $U(\delta)\cup \square(\delta)$.
\end{Lemma}
\begin{proof}
  We will work in $\square(1-\delta)$ and show that there is a self-homeomorphism $h$
  of $\square(1-\delta)$ which fixes $\square(\delta)\cup \partial\square(1-\delta)$
  and straightens $J$. 
  The lemma then follows by extending $h$ by identity 
  to $\square\setminus \square(1-\delta)$.
  W.l.o.g. assume that $J(1-\delta)=d_+(\delta)$.  Using a tubular
  neighborhood of $J$ find a parallel curve $J'$ such that
  $J'(0)\in \partial \square(2\delta)$,
  $J'(1)\in \partial\square(1-2\delta)$, and $J'\cap J=\es$. 
  Then $J'(1)\in\partial U(2\delta)$.

  The set $A=\square(2\delta)\setminus (\square(\delta)\cup J)$ is 
  connected and $b_-(\delta)$ is on the 
  boundary of $A$. So there is $\gamma_1\subset A\cup \{b_-(\delta)\}$ 
  which starts at $b_-(\delta)$ and ends at $J'(0)$. 
  The set $B=U(2\delta)\setminus (U(\delta)\cup J)$ is also connected
  and $d_-(\delta)$ is on the boundary of $B$, so there
  is a curve $\gamma_2\subset B\cup \{d_-(\delta\}$ which starts at
  $J'(1)$ and ends at $d_-(\delta)$. 
  Let $J_1=J\rest [0,1-\delta]$ and let
  $$J_2=[o,b_-(\delta)]\cup \gamma_1\cup J'\cup \gamma_2$$
  Then $J_1$ and $J_2$ satisfy the assumptions of Lemma~\ref{lemma:StraightenTwoCurves}
  for $\square(1-\delta)$ instead of $\square$. So by Corollary~\ref{cor:Rescale}
  there is a required homeomorphism.
\end{proof}

\begin{Lemma}\label{lemma:WLOGJ}
  Let $T$ be a centered subdivided t-diagram. Let $J_0$ be the ray
  $J_0(t)=(1-t)o+tv^+$.  Then there is an $\widehat T\sim^*_0 T$ such that
  $$\intersect(J,\widehat T^-)=\MIN_{\ray}(T^-)=\MIN_{\ray}(\widehat T^-).$$
\end{Lemma}
\begin{proof}
  Let $J$ be a ray such that $\intersect(J,T^-)=\MIN_{\ray}(T^-)$.
  Let $\delta>0$ be such that $B(o,\delta)\cap T^-$ is one line
  segment and also $U(\delta)\cap T^-$ is one line segment. Just in
  case assume also that $\delta<1/8$ so that $\square_{\delta}$ and
  $U(\delta)$ are disjoint.  Let $t$ be the maximum of the set
  $J^{-1}[\square_{\delta}]$.  Let $\gamma$ be a curve connecting
  $b_+(\delta)$ to $J(t)$ inside
  $(\square_\delta\setminus T^-)\cup \{b_+(\delta)\}$.  Let
  $J'=[o,b_+(\delta)]\cup \gamma\cup J\rest [\delta,1]$.  Then $J'$ is
  still a ray witnessing $\MIN_{\ray}(T^-)$. Let $t'$ be the minimum
  of $(J')^{-1}[U(\delta)]$. Let $\gamma'$ be a curve in
  $(U(\delta)\setminus T^-)\cup \{v^+\}$ connecting $J'(t')$ to $v^+$.
  Let $J''=(J'\rest [0,t'])\cup \gamma'$.  Then $J''$ still witnesses
  $\MIN_{\ray}(T^-)$. After reparametrization one can make sure that
  $J''$ satisfies the assumptions of the Straightening Lemma
  \ref{lemma:StraightenJ}, so let $h$ be the homeomorphism given by
  that Lemma. Then $h$ witnesses an $R^*_0$-move which results in
  $\widehat T=hT$ and $hJ=[0,v^+]$ witnesses
  $\MIN_{\ray}(T^-)=\MIN_{\ray}(\widehat T^-)$ by
  Lemma~\ref{lemma:invariance_under_R0}.
  %
  %
\end{proof}

\begin{figure}
  \centering
  \begin{subfigure}[t]{0.3\textwidth}
    \centering
    \includegraphics[width=\linewidth]{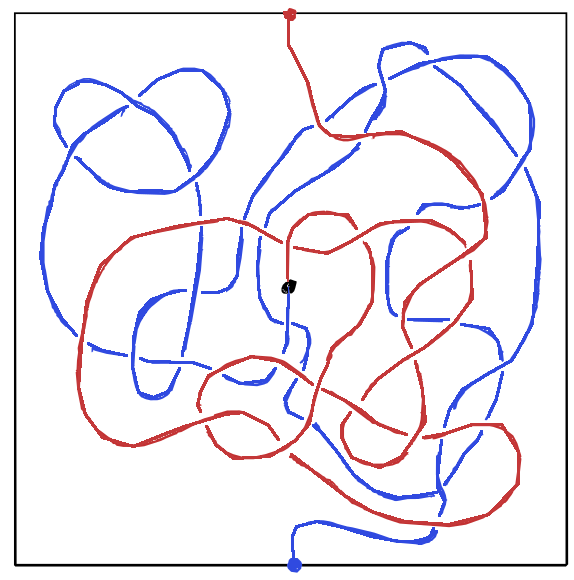}
    \caption{}
  \end{subfigure}
  \hfill
  \begin{subfigure}[t]{0.3\textwidth}
    \centering
    \includegraphics[width=\linewidth]{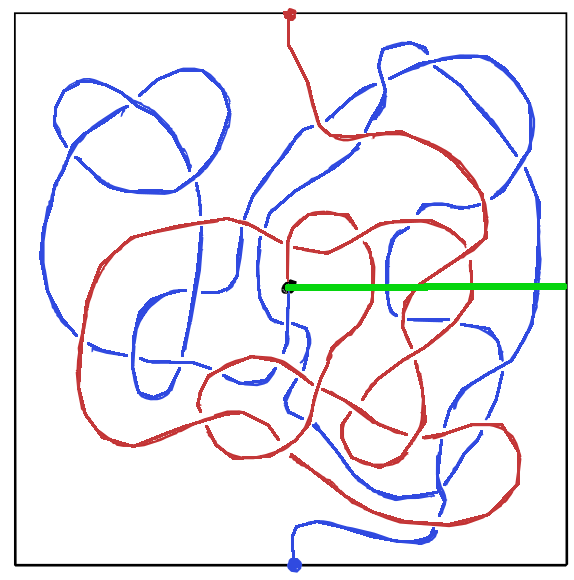}
    \caption{}
  \end{subfigure}
  \hfill
  \begin{subfigure}[t]{0.3\textwidth}
    \centering
    \includegraphics[width=\linewidth]{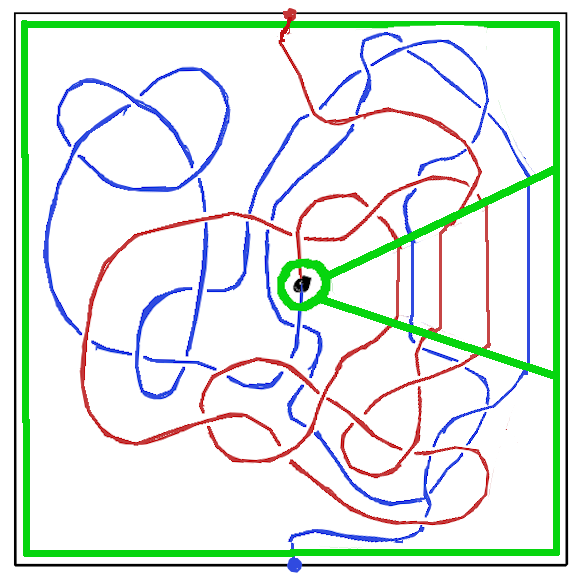}
    \caption{}
  \end{subfigure}

  \vspace{1em} 

  \begin{subfigure}[t]{0.3\textwidth}
    \centering
    \includegraphics[width=\linewidth]{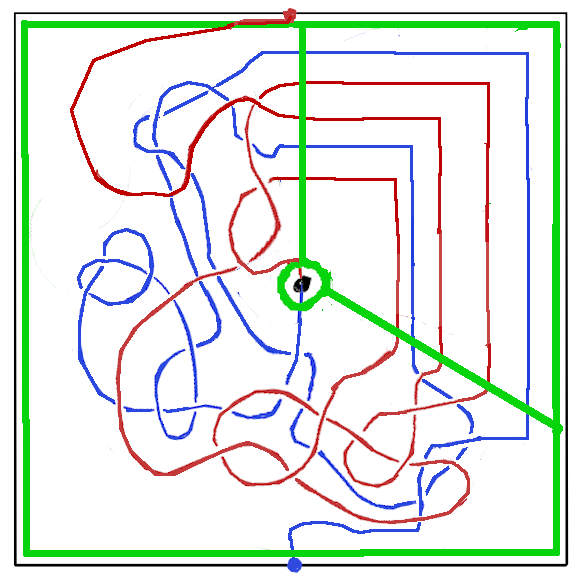}
    \caption{}
  \end{subfigure}
  \hfill
  \begin{subfigure}[t]{0.3\textwidth}
    \centering
    \includegraphics[width=\linewidth]{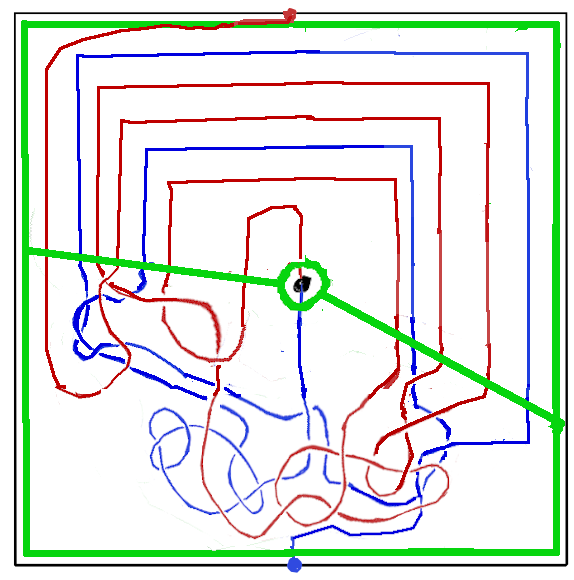}
    \caption{}
  \end{subfigure}
  \hfill
  \begin{subfigure}[t]{0.3\textwidth}
    \centering
    \includegraphics[width=\linewidth]{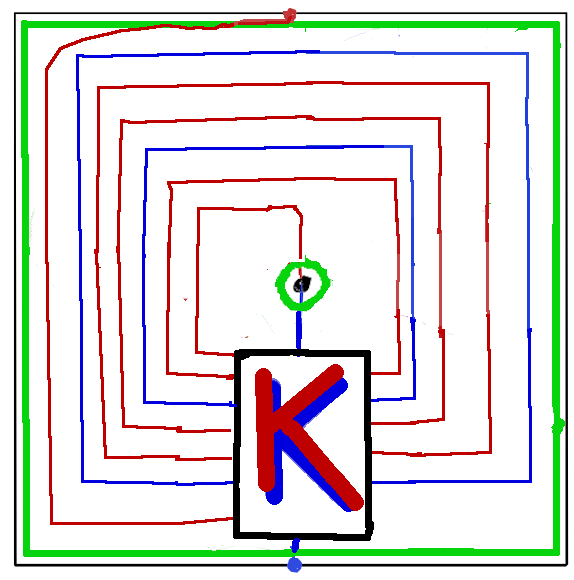}
    \caption{}
  \end{subfigure}

  \vspace{1em} 

  \begin{subfigure}[t]{0.49\textwidth}
    \centering
    \includegraphics[width=\linewidth]{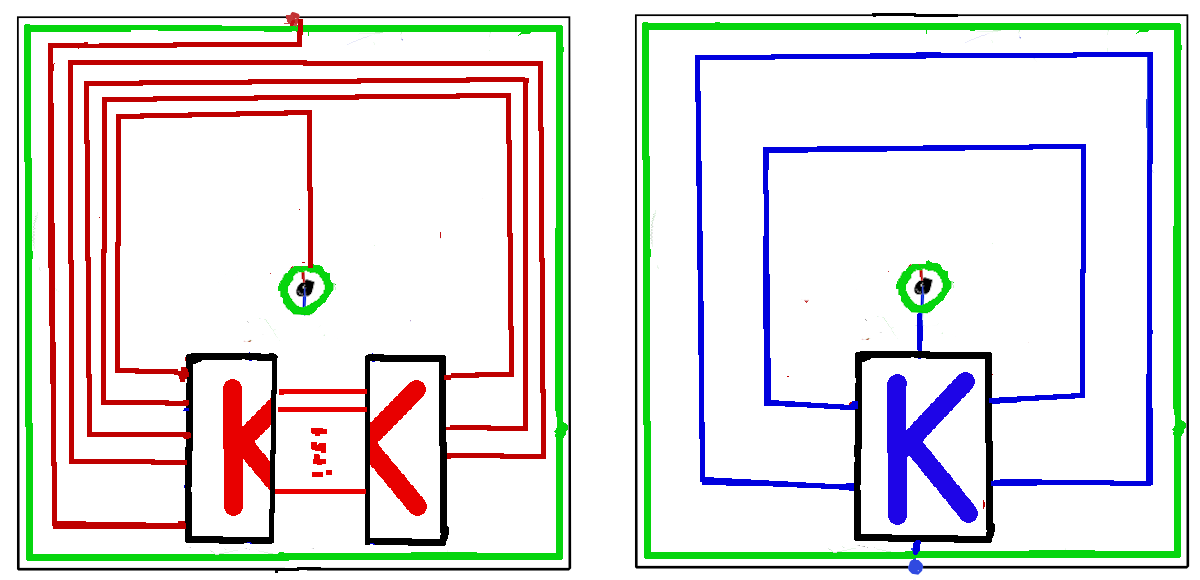}
    \caption{}
  \end{subfigure}
  \hfill
  \begin{subfigure}[t]{0.49\textwidth}
    \centering
    \includegraphics[width=\linewidth]{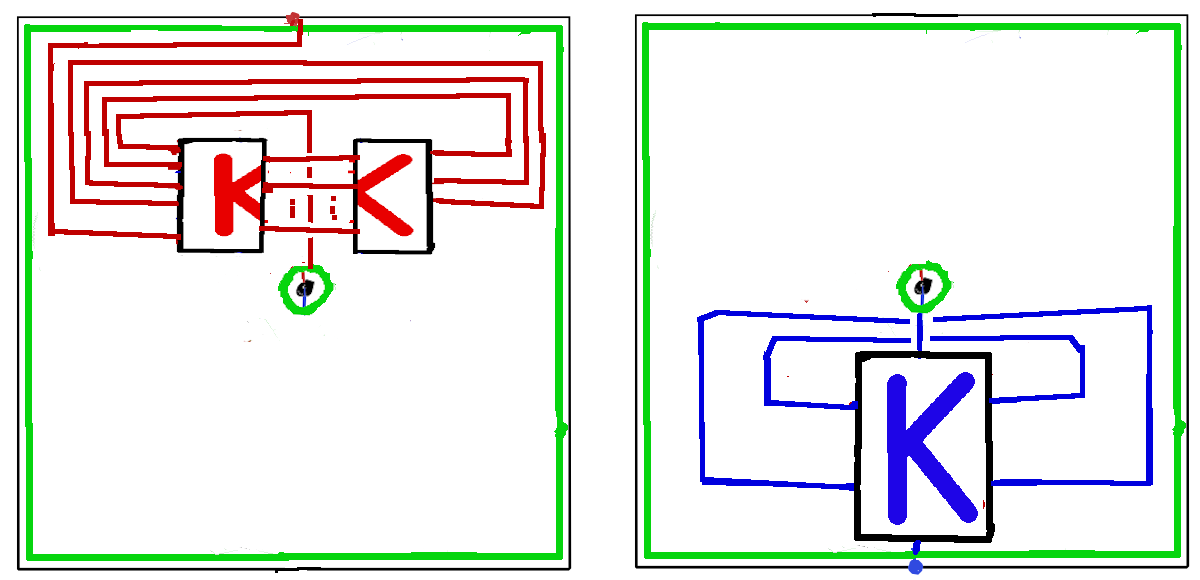}
    \caption{}
  \end{subfigure}

  \vspace{1em} 

  \begin{subfigure}[t]{0.3\textwidth}
    \centering
    \includegraphics[width=\linewidth]{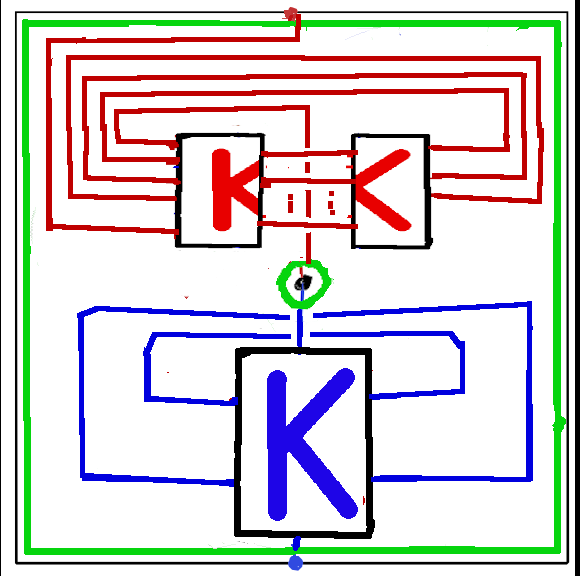}
    \caption{}
  \end{subfigure}
  \caption{}
  \label{fig:Anim}
\end{figure}

\subsection{The Fan Lemma}

Recall the notions of axis-aligned rectangle and homeomorphism, a
$G$-tangle, $T\capdot G$ from Definitions~\ref{def:G-tangle} and
\ref{def:TcapG-tangle} and their adaptations to subdivided t-diagrams
from Definition~\ref{def:SubdividedTangleDiagram}. 
Fix $G$ to be the square
$$G=\left[\tfrac{3}{8},\tfrac{5}{8}\right]\times \left[\tfrac{1}{8},\tfrac{3}{8}\right]$$
Note that it is axis-aligned (Definition~\ref{def:AxisAligned}).
Fix also $g$ to be the unique axis-aligned homeomorphism $g\colon G\to \square$
(cf. Definition~\ref{fact:AxisAligned}).
Let us introduce the notation
$\chi_h(T)=\intersect(T^-,T^+)$ to denote the number of ``hybrid''
crossings in a tangle diagram.

\begin{Lemma}[The Fan Lemma]\label{lemma:Fan}
  Suppose that $T=(T,T^+,T^-,\{1\},M,\ell)$ is a centered subdivided
  t-diagram with one strand and $M$ loops. Then there is a centered
  subdivided t-diagram $\Fan(T)\sim^*_0 T$,
  $$\Fan(T)=(\Fan(T),\Fan(T)^+,\Fan(T)^-,\{1\},M,\Fan(\ell)),$$
  with the following properties. The diagrams $\Fan(T)$, $\Fan(T)^+$,
  and $\Fan(T)^-$ are in general position with respect to~$G$. Denote
  $S^-=g(\Fan(T)^-\capdot G)$ and $S^+=g(\Fan(T)^+\capdot G)$. Then:
  \begin{enumerate}[label={\upshape ($F{\arabic*}$)}, leftmargin=3pc]
  \item $\Cross(\Fan(T))\subset G$ and $\Cross_h(\Fan(T))\subset G$\label{fan1}
  \item $\MIN_+(\Fan(T)^-)=\MIN_{\ray}(T^-)=\MIN_{\ray}(\Fan(T)^-)$ \label{fan2}
  \item $\MIN_{\vert}(S^-)=\MIN_{+}(\Fan(T)^-)$ \label{fan3}
  \item $S^-$ is of type  $(1,k-1,1,k-1)$ where $k=\MIN_{\vert}(S^-)+1$ \label{fan4}
  \item $\chi_h(\Fan(T))=\intersect(S^+, S^-)$ \label{fan5}
  \item $S^+$ is of type $(0,l+2,0,l)$ for some $l$. \label{fan6}
  \item $\MIN_{\vert}(S^+)=\MIN_-(\Fan(T)^+)$ \label{fan7}
  \item $e(S^-,S^+)$ satisfies property $(*)$ (Definitions~\ref{def:eTS} and \ref{def:PropertyStar}).
    \label{fan8}
  \end{enumerate}
\end{Lemma}
\begin{proof}
  By Lemma \ref{lemma:WLOGJ} assume without loss of generality that
  $T$ is the $\Fan(T)$ of Lemma \ref{lemma:WLOGJ}, i.e.
  $\MIN_{\ray}(T^-)$ is realized by $J_0=[o,v^+]$.  We will describe the
  process depicted in Figures~\ref{fig:Anim}(b)-(f). For any
  $\delta>0$ denote
  $$\tau(\delta)= [b_+(\d),v^+-(\d,0)]\times \left[\tfrac{1}{2}-\d,\tfrac{1}{2}+\d\right].$$
  By general position and centeredness of $T$, we can assume
  without loss of generality that for a sufficiently small $\d$ 
  the set $T\cap \tau(\delta)$ is a collection of
  vertical straight lines and the number of those lines equals
  $\MIN_{\ray}(T)$.  Assume that $\delta$ is also so small that
  each of the (four) sets
  $\square_{\delta}\cap T^{\pm}$ and
  $U(\delta)\cap T^{\pm}$ is one line segment where $U(\delta)$ is as
  defined before the Straightening Lemma~\ref{lemma:StraightenJ}, see Figures~\ref{fig:Xi_1}(c)
  and \ref{fig:LastHomeo}(a).
  We will now construct a homeomorphism $\psi\colon\square\to\square$
  which will realize the fan-move depicted on
  Figures~\ref{fig:Anim}(b)-(f).  Let
  Let $\eta'\colon \partial\square\to S^1$ be the homeomorphism
  given by $\eta'(x)=(x-o)/|x-o|$ and let
  $\eta\colon\square\to D^2=\{(x,y)\in\R^2\mid x^2+y^2\le 1\}$ be 
  be the radial extension of $\eta'$ such that $\eta(o)=\bar 0$.
  Note that for all $r\le 1$, $\partial\square_r$ maps to $\partial B(\bar 0,r)$.

  \begin{figure}
    \centering
    \begin{subfigure}{0.45\textwidth}    
      \centering
      \fbox{\includegraphics[width=\linewidth]{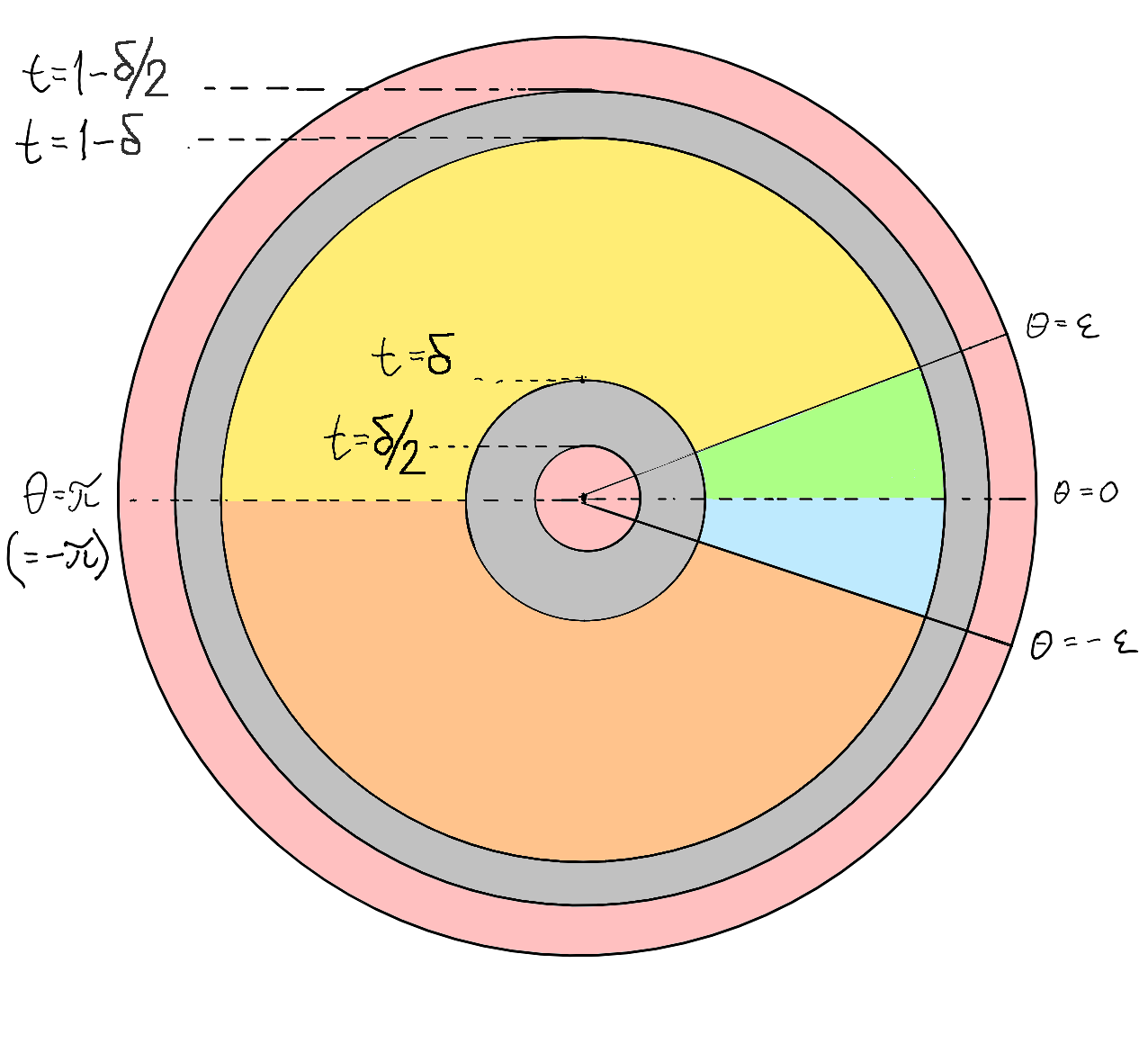}}
      \caption{}
    \end{subfigure}
    \raisebox{8\height}{$\bigstackarrowww{\psi}$}
    \begin{subfigure}{0.45\textwidth}
      \centering
      \fbox{\includegraphics[width=\linewidth]{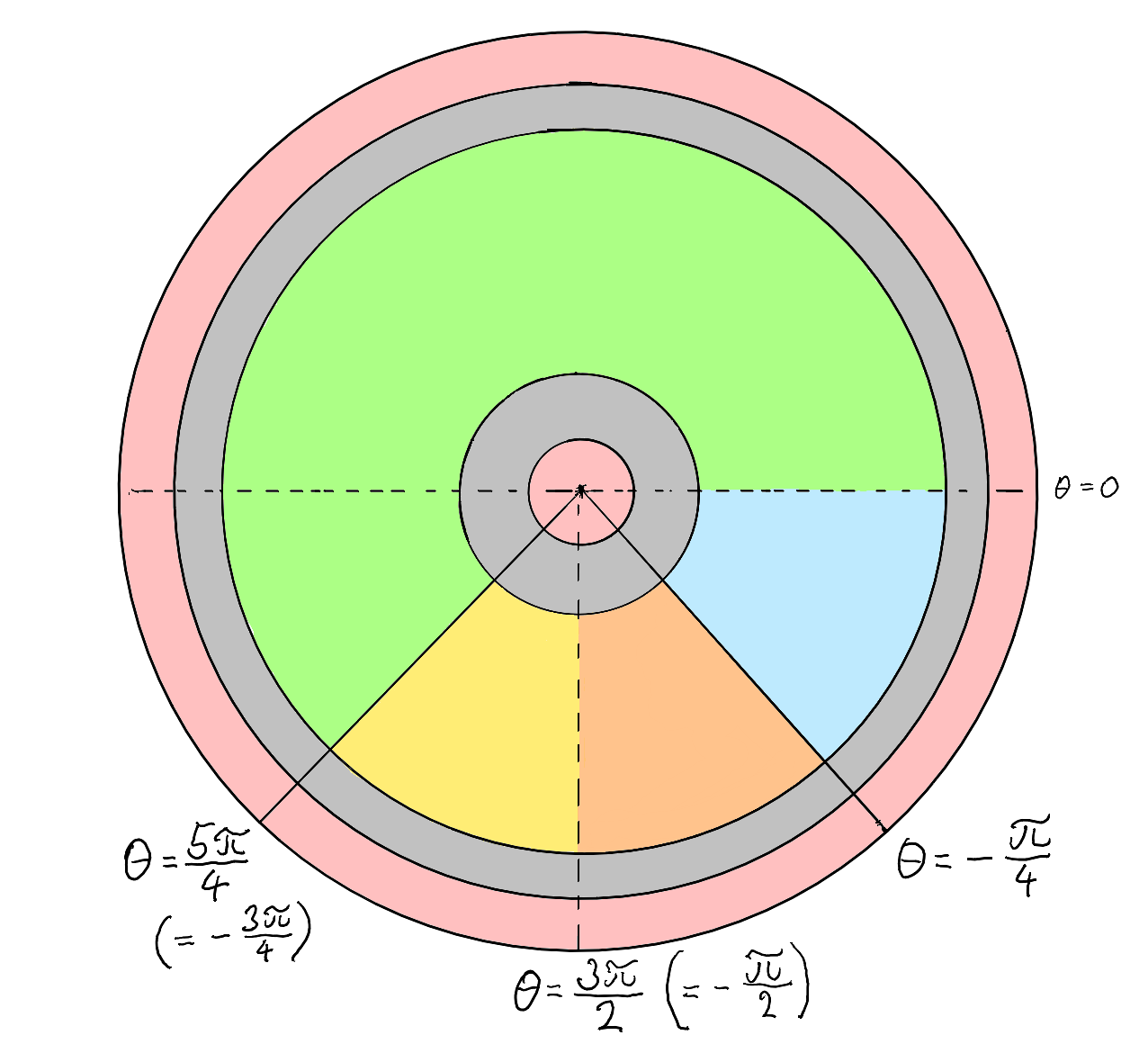}}
      \caption{}
    \end{subfigure}
    
    \vspace{1em}
    
    \begin{subfigure}{0.46\textwidth}
      \includegraphics[width=\linewidth]{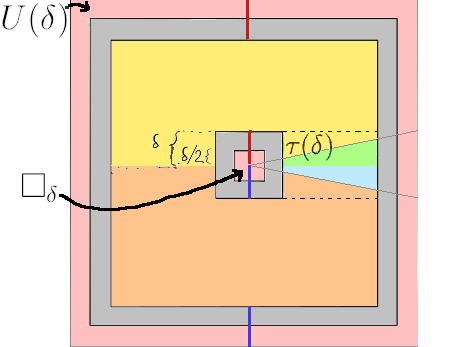}\hfill
      \caption{}
    \end{subfigure}
    \hspace{0em}
    \raisebox{6\height}{$\bigstackarrowww{\eta^{-1} \psi \eta}$}
    \hspace{0em}
    \begin{subfigure}{0.42\textwidth}
      \centering
      \includegraphics[width=\linewidth]{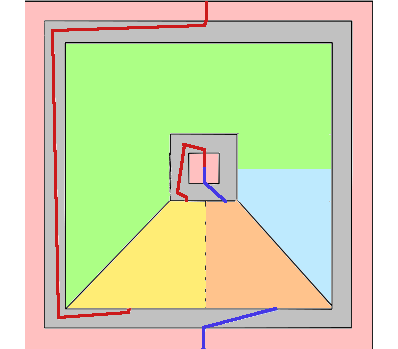}
      \caption{}
    \end{subfigure}

    \caption{Illustrating the maps $\psi$ and
      $\eta^{-1}\circ\psi\circ\eta$ from the proof of
      Lemma~\ref{lemma:Fan}.  Figure (a) is the domain and (b) the
      range of $\psi$. The colored regions are mapped respectively to
      the same colored regions. The light-red area is where $\psi$ is
      the identity. The yellow and orange areas are squeezed and the
      blue and green areas are stretched. The gray areas is where
      $\psi$ is defined by interpolating between the light-red and
      other colored areas. In (c) and (d) the domain and range of
      $\eta^{-1}\circ\psi\circ\eta$ is shown. The green and blue areas
      form an area inside of $\tau(\delta)$ according to
      \eqref{eq:InsideTau}. Figures (c) also shows the line segments
      $K^+\cap U(\delta)$ and $K^+\cap \square_\delta$ (red) as well
      as $K^-\cap U(\delta)$ and $K^-\cap \square_\delta$ (blue). In
      (d) the images of these line segments are shown. Notice how the
      gray area ``drags'' them into the lower region due to the
      interpolation.  }
    \label{fig:Xi_1}
  \end{figure}

  For convenience, denote $\la t,\theta\ra=te^{i\theta}\in D^2$ (for $t\le 1$).
  There is $\e>0$ such that
  \begin{equation}
    \eta[\tau(\delta)]\supset \{te^{i\theta}\mid t\in [\delta,1-\delta],\theta\in [-\e,\e]\}.\label{eq:InsideTau}  
  \end{equation}
  Let $\psi\colon D^2\to D^2$ be given by
  $$
  \psi(\la t,\theta\ra)=
  \begin{cases}
    \left\la t,\pi\dfrac{5\theta}{4\e}\right\ra&\text{ if }t\in [\delta,1-\delta]\text{ and }\theta\in [0,\e]\\
    \left\la t,\pi\dfrac{\theta}{4\e} \right\ra&\text{ if }t\in [\delta,1-\delta]\text{ and }\theta\in [-\e,0]\\
    \left\la t,\pi\dfrac{\theta-\e}{2(\pi-\e)}+1\right\ra&\text{ if }t\in [\delta,1-\delta]\text{ and }\theta\in [\e,\pi]\\
    \left\la t,\pi\dfrac{\theta+\e}{2(\e-\pi)}-\frac{1}{4}\right\ra&\text{ if }t\in [\delta,1-\delta]\text{ and }\theta\in [-\pi,-\e]\\
    \left\la t,\theta\right\ra&\text{ if }t<\delta/2\text{ or }t>1-\delta/2
  \end{cases}
  $$
  Note that $\psi(\la t,\theta\ra)$ is not yet defined for $t\in [\delta/2,\delta]\cup [1-\delta,1-\delta/2]$. For these values of $t$, $\psi(t)$ will be defined by
  a ``circular'' interpolation as follows. For $t\in [\delta/2,\delta]$ let
  $$
  \psi(t,\theta)=\Big\la t,\ (2-2t/\delta)\Theta\big(\psi(\la \delta/2,\theta\ra)\big)+(2t/\delta-1)\Theta\big(\psi(\la \delta,\theta\ra)\big) \Big\ra
  $$
  and for $t\in [1-\delta,1-\delta/2]$ let
  $$
  \psi(t,\theta)=\Big\la t,\ (2(1-t)/\delta-1)\Theta\big(\psi(\la 1-1/8,\theta\ra)\big)+(2(t-1)/\delta+2)\Theta\big(\psi(\la 1-\delta/2,\theta\ra)\big) \Big\ra
  $$
  where $\Theta(\la t,\theta \ra)=\theta$. See Figure~\ref{fig:Xi_1} for a visual explanation.

  Consider the map $\eta^{-1}\circ\psi\circ\eta$. It performs the
  corresponding shift in $\square$ and stretches what used to be short
  vertical lines in $\tau$ in a fan-like manner to go around $o$ as
  shown in Figures~\ref{fig:Anim}(b)-(f), although we are still
  missing the final step. We should package the lower area into the
  rectangle $G=[1/4,3/4]\times [0,1/2-\delta]$. We will not explicitly construct this last
  homeomorphism $\zeta$ and will be content with
  Figure~\ref{fig:LastHomeo}. Let also $g\colon G\to \square$
  be the (unique) axis-aligned linear homeomorphism.

  Both homeomorphisms $\eta^{-1}\psi\eta$ and $\zeta$ keep
  $\partial\square$ and a neighborhood of $o$ fixed, so
  $h=\zeta\eta^{-1}\psi\eta$ realizes an $R^*_0$-move.  Let
  $\Fan(T)=\zeta\eta^{-1}\psi\eta T$. The final result without extra
  colors is depicted on Figure~\ref{fig:squareMtoG}(a).  If $x$ is a
  crossing of $\Fan(T)$, then there is a crossing $y$ of $T$ such that
  $h(y)=x$. Since there are no crossings in
  $\tau\cup \square_\delta\cup U(\delta)$, $y$ must be such that if
  $\eta(y)=\la t,\theta\ra$, then $t\in [\delta,1-\delta]$ and
  $\theta\notin [-\e,\e]$ (yellow or orange areas in Figures
  \ref{fig:Xi_1}(a),(c)).  But then it is easy to check that
  $\psi(\eta(y))$ ends up in the area
  $\{\la t,\theta\ra\mid t\in [\delta,1-\delta],\theta\in
  [-3\pi/4,\pi/4]\}$ (yellow and orange in
  Fig.~\ref{fig:Xi_1}(b),(d)).  This is in turn which is mapped by
  $\eta^{-1}$ to an area which is then squeezed into $G$ by~$\zeta$.
  Thus, $x\in G$ which proves the first part of \ref{fan1}.  But
  $\Cross_h(T)\subset\Cross(T)$, so it proves also the second part.
  The fact that $\MIN_{\ray}(T)=\MIN_{\ray}(\Fan(T))$ follows from
  Lemma~\ref{lemma:invariance_under_R0}.  By
  Lemma~\ref{lemma:Raypm_ineq},
  $\MIN_+(\Fan(T)^-)\ge \MIN_{\ray}(\Fan(T))$ so we need to prove
  $\MIN_+(\Fan(T)^-)\le \MIN_{\ray}(T)$.  Let $J$ be the vertical line
  $[o,o^+]$. It is clear from the construction that $h^{-1}J$
  intersects $T^-$ exactly in the vertical lines of $T^-\cap \tau$.
  In fact $h^{-1}J\subset \square_\delta\cup J_0\cup U(\delta)$ and
  $U(\delta)\cap T^-\cap h^{-1}J=\es$ and
  $\square_\delta\cap T^-\cap h^{-1}J=\{o\}$.  Here $J_0$ is from the
  very beginning of this proof.  Thus,
  $\intersect (J,\Fan(T)^-)=\intersect(h^{-1}J,T^-)=\intersect(J_0,T^-)=\MIN_{\ray}(T^-)$.
  This proves \ref{fan2}. Before proving \ref{fan3}, let us note that
  the areas in the vertical edges of $G$ which intersect $T^-$ are
  images under $h$ of the horizontal edges of $\tau$.  On the other
  hand the horizontal edges of $G$ intersect $T^-$ each in one point.
  The lower edge at $o^-$ and the top edge because it is the image of
  the lower edge of $\square_\delta$ under $h$. It therefore follows
  that
  \begin{equation}
    S^-\text{ is of type }(1,\MIN_{\ray}(T^-),1,\MIN_{\ray}(T^-))\label{eq:IsOfType}  
  \end{equation}
  Let us now prove \ref{fan3}. By Lemma~\ref{lemma:UpperBoundOnMIN}
  it follows from \eqref{eq:IsOfType} that
  $\MIN_{\vert}(S^-)\le \MIN_{\ray}(T^-)$ which by \ref{fan2}
  implies
  $\MIN_{\vert}(S^-)\le \MIN_{+}(T^-)$. Let us prove the other direction.
  If $J_v$ is a vertical splitting curve
  of $S^-$ in $\square$, then $g^{-1}J_v$ is a vertical splitting
  curve going vertically through $G$. It can be extended at the top
  and connected to $o$ without introducing new intersections with
  $T^-$ other than at $o$ which does not count (see
  Definition~\ref{def:lowersplitting}). Fundamentally, the reason is
  that $T^-\cap \square_\delta$ is a single line segment, so
  $((\square_{\delta}\setminus T^-)\cap\square^-)\cup \{o\}$ is
  path-connected, see Figure~\ref{fig:LastHomeo}(d). In this way we
  obtain a ray so it intersects $T^-$ in at least $\MIN_{\ray}(T^-)$
  many points, so we have shown
  $\MIN_{\vert}(S^-)\ge\MIN_{\ray}(T^-)$. By \ref{fan2} this proves
  $\MIN_{\vert}(S^-)\ge\MIN_{+}(T^-)$. This completes the proof of \ref{fan3}.
  Next, \ref{fan4} follows from \ref{fan3} together with \eqref{eq:IsOfType}.
  Let us prove \ref{fan5}. Since $S^+$ and $S^-$ are complete tangle diagrams (there are no loose ends), 
  $\ocirc{S}^+\cap \ocirc{S^-}=S^+\cap S^-$, so $\intersect(S^+,S^-)=\Card(S^+\cap S^-)$
  (see Definition~\ref{def:MutualGP}\eqref{def:Intersect} and note that
  the ``interior'', denoted by $\ocirc{\cdot}$, means here that we simply remove the endpoints of the strands, not interior of the curves as sets which would render them empty).
  For $\Fan(T)^\pm$ denote the interiod by $\Fan^\circ(T)^{\pm}$. We have,  
  $$\Fan(T)^+\cap \Fan(T)^-=\Big(\Fan^\circ(T)^+\cap \Fan^\circ(T)^-\Big)\cup \{o\},$$
  but $o\notin G$,
  so 
  $$\Fan(T)^+\cap \Fan(T)^-\cap G=\Fan^\circ(T)^+\cap \Fan^\circ(T)^-\cap G$$ 
  and therefore,
  by \ref{fan1}, $\Cross_h(\Fan(T))=\Fan(T)^+\cap \Fan(T)^-\cap G$
  We now have:
  \begin{align*}
    \intersect(S^+,S^-)&=\Card(g(\Fan(T)^+\cap G)\cap g(\Fan(T)^-\cap G))&&\text{by definition of }S^{\pm}\text{ and the above}\\
                       &=\Card((\Fan(T)^+\cap G)\cap (\Fan(T)^-\cap G))&&\text{because }g\text{ is a homeomorphism}\\
                       &=\Card(\Fan(T)^+\cap \Fan(T)^-\cap G)\\
                       &=\Card(\Cross_h(\Fan(T)))&&\text{by the above}\\
                       &=\chi_h(\Fan(T))
  \end{align*}
  For \ref{fan6}, let $l$ be the number of lines in $T^+\cap \tau$. From the construction
  it follows that this is the same as the number of lines of $\Fan(T)^+$ entering
  $G$ through the right vertical edge. On the left vertical edge two additional
  lines of $\Fan(T)^+$ come into $G$ from the areas $h\square_\delta$ and $hU(\delta)$
  (see Figure~\ref{fig:LastHomeo}(c),(d) and the caption of that figure).
  No $\Fan(T)^+$ lines intersect $G$ at the horizontal edges. This proves~\ref{fan6}.

  Let $J_v$ be a vertical splitting curve for $S^+$ with
  $i(J_v,S^+)=\MIN_{\vert}(S^+)$. Then $g^{-1}J$ can be modified in
  the same manner as in the proof of \ref{fan3} above to be connected
  to $o$ without introducing new intersections with $\Fan(T)^+$. Further, it can be modified within $U(\delta)\cap \square^-$
  to be connected to $o^-$ still without new intersections with
  $\Fan(T)^+$ (recall that $(U(\delta)\setminus T^+)\cup \{o^-\}$ is
  path-connected by the choice of~$\delta$, and so is this set
  intersected with~$\square^-$). This shows that
  \begin{equation}
    \MIN_{\vert}(S^+)\ge \MIN_{-}(\Fan(T)^+).\label{eq:ineqMINvertgeMin-}  
  \end{equation}
  Suppose that $J^-$ is a
  lower splitting curve for $T^+$ such that $\intersect(J^-,\Fan(T)^+)=\MIN_-(\Fan(T)^+)$.
  If $g(J^-\cap G)$ is a vertical
  splitting curve, then
  \begin{align*}
    \MIN_{\vert}(S^+)&\le \intersect(g(J^-\cap G),S^+)\\
    &=\intersect(J^-\cap G,T^+\cap G)&&\text{because }g\text{ is homeo and def. of }S^+\\
    &\le \intersect(J^-,T^+)&&\text{because }J^-\cap G\cap T^+\cap G\subset J^-\cap T^+.\\
    &=\MIN_-(\Fan(T)^+)&&\text{by the choice of }J^-.                               
  \end{align*}
  Otherwise, note that there is a homeomorphic transformation $\f$ of $\square^-$ into
  its subset such that any lower splitting curve is mapped onto a curve
  whose intersection with $G$ is a vertical splitting curve of~$G$,
  and such that $\f[\Fan(T)^+]=\Fan(T)^+\cap G$. The construction of $\f$ is illustrated in
  Figure~\ref{fig:squareMtoG}.  Let $J=g^{-1}\f J^-$. Then by the
  construction of $g$ (Figure~\ref{fig:squareMtoG}) $J$ is a vertical
  splitting curve for $S^+$ and we have
  $$(g^{-1}\f J^-)\cap S^+=\f J^-\cap \Fan(T)^+\cap G=J^-\cap \f^{-1}\Fan(T)^+=J^-\cap \Fan(T)^+,$$
  so
  $\intersect(J,S^+)= \intersect(J^-,\Fan(T)^+)= \MIN_{-}(\Fan(T)^+)$,
  so it again follows that $\MIN_{\vert}(S^+)\le \MIN_-(\Fan(T)^+)$. Together with \eqref{eq:ineqMINvertgeMin-} this proves~\ref{fan7}.
  
  It remains to prove \ref{fan8}, so let us show that $e(S^-,S^+)$
  satisfies property $(*)$ (Definition~\ref{def:PropertyStar}). To be
  consistent with Definition~\ref{def:eTS}, denote $T=S^-$ and
  $S=S^+$.  Let $S^*$ be a smoothing of $S$ as in
  Definition~\ref{def:eTS}. Then endpoints of $S^*$ are the same as
  the endpoints of $S$. Recall the notation $L_i(T)$
  (Definition~\ref{def:Region}). These are open subarcs of the left
  and right edge of $\square$ separated by the endpoints of
  $T=S^-$. Via $g^{-1}$ they become subarcs of vertical edges of
  $G$. Since there are $k-1$ endpoints of $T^-$ no each side of $G$
  by (\ref{fan4}), there are $k$ such subarcs on each side.  Let
  $c_{i}$ be the number of endpoints of $S=S^+$ in $L_i(T)$.  Recall
  that the sets $L_1(T),\dots,L_n(T)$ ($n=2k$) are enumerated
  counterclockwise starting from the top left. Since the arcs that go
  around in the fan correspond to the original lines in $\tau\cap T$,
  each point $S\cap L_i(T)$ is matched by a point in
  $S\cap L_{f(i)}(T)$ where $f(i)=n-i+1$ as defined in the beginning
  of Section~\ref{ssec:Matrices}. Figures \ref{fig:LastHomeo}(d),
  \ref{fig:squareMtoG}(a) as well as \ref{fig:Anim}(f) illustrate this
  matching via red arcs. Therefore, if $i\in\dint{2,k-1}$, we have
  $c_i=c_{f(i)}$. But $c_i=s_i(e(T,S))$, so this
  proves~\ref{item:star2}.  When $k>1$, there is one extra endpoint
  of $S$ in $L_1(T)$ corresponding to the arc of $T^+$ which
  originates in $o$ and one extra endpoint in $L_k(T)$ corresponding
  to the arc which originates in $o^+$ (red arcs in the figures).
  This proves \ref{item:star3}. Clearly, if an arc goes from $L_i(T)$
  to $L_j(T)$, the same arc goes from $L_j(T)$ to $L_i(T)$, so the
  matrix is symmetric verifying $e\in \MM^s_n(\N)$ and condition
  \ref{item:star4} follows immediately from Definition~\ref{def:eTS}
  of~$e(T,S)$.
\end{proof}

\begin{figure}
  \centering
  \begin{subfigure}[t]{0.24\linewidth}    
    \centering
    \fbox{\includegraphics[width=\linewidth]{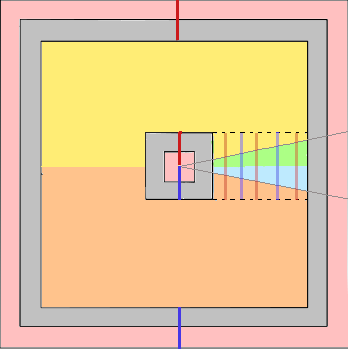}}
    \caption{}
  \end{subfigure}
  \hfill
  \begin{subfigure}[t]{0.24\linewidth}
    \centering
    \fbox{\includegraphics[width=\linewidth]{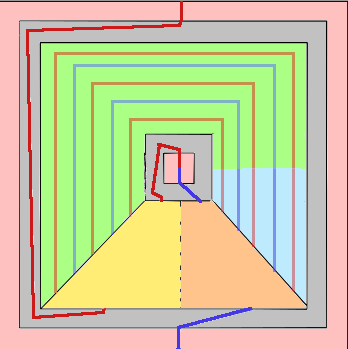}}
    \caption{}
  \end{subfigure}
  \hfill
  \begin{subfigure}[t]{0.24\linewidth}
    \centering
    \fbox{\includegraphics[width=\linewidth]{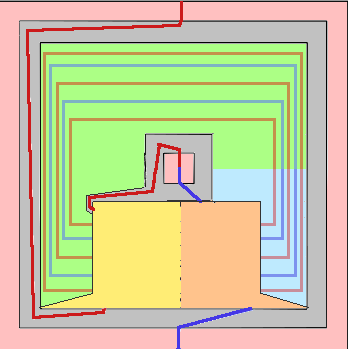}}
    \caption{}
  \end{subfigure}  
  \hfill
  \begin{subfigure}[t]{0.24\linewidth}
    \centering
    \fbox{\includegraphics[width=\linewidth]{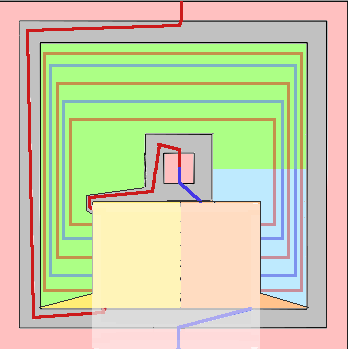}}
    \caption{}
  \end{subfigure}
  \caption{\textbf{(a)}: Same as Figure~\ref{fig:Xi_1}(c) but with the vertical
    lines of $K\cap \tau$ shown. \textbf{(b)}: Same as Figure~\ref{fig:Xi_1}(c)
    but with the stretched lines shown.  Notice how the gray area
    ``drags'' the two red lines into the lower region due to the
    nature of the interpolation.  \textbf{(c)}: The result of the application
    of the homeomorphism $\zeta$ which ``packages'' the lower part of
    the knot into the rectangle $G$. \textbf{(d)}: The rectangle $G$ is
    overlayed in white. All the properties \ref{fan1}--\ref{fan8} of
    the Fan Lemma~\ref{lemma:Fan} can be ``read'' from Figure (d). For
    example the number of blue lines (lines of $K^-$) which go around
    the ``fan'' equals to $\MIN_{\ray}(K^-)=\MIN_{\ray}(\Fan(K)^-)$
    because it is the same as the number of blue lines shown in (a) to
    the right of the small gray square (denoted $\square_\delta$ in
    the proof).  The ``+2'' in \ref{fan6} comes from the two red lines
    which arrive into $G$ from the gray areas.}
  \label{fig:LastHomeo}
\end{figure}

\begin{figure}
  \centering
  \centering
  \begin{subfigure}[t]{0.31\linewidth}    
    \centering
    \fbox{\includegraphics[width=\linewidth]{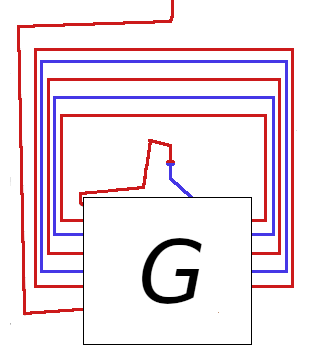}}
    \caption{}
  \end{subfigure}
  \hfill
  \begin{subfigure}[t]{0.31\linewidth}
    \centering
    \fbox{\includegraphics[width=\linewidth]{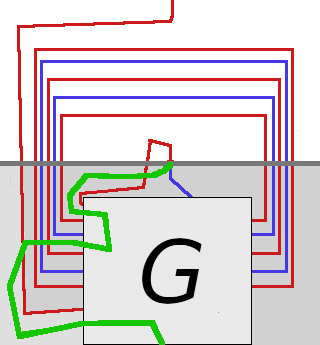}}
    \caption{}
  \end{subfigure}
  \hfill
  \begin{subfigure}[t]{0.31\linewidth}
    \centering
    \fbox{\includegraphics[width=\linewidth]{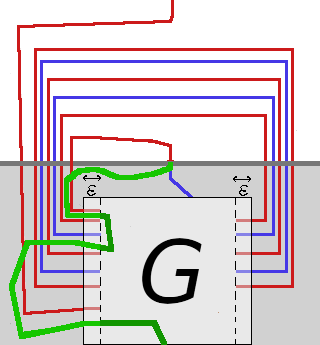}}
    \caption{}
  \end{subfigure}  

  \vspace{1em}
  
  \begin{subfigure}[t]{0.31\linewidth}
    \centering
    \fbox{\includegraphics[width=\linewidth]{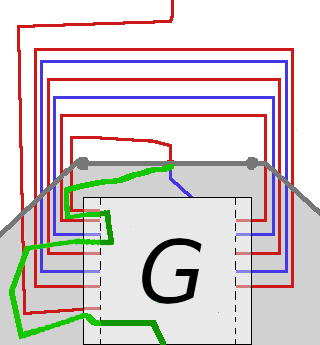}}
    \caption{}
  \end{subfigure}  
  \hspace{10pt}
  \begin{subfigure}[t]{0.31\linewidth}
    \centering
    \fbox{\includegraphics[width=\linewidth]{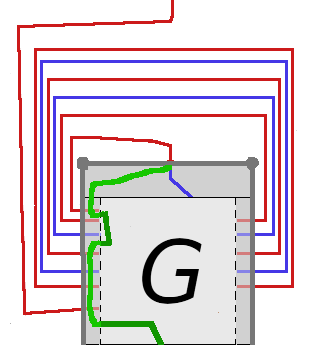}}
    \caption{}
  \end{subfigure}

  \caption{\textbf{(a)}: The result of the procedure depicted in
    Figures~\ref{fig:Xi_1} and \ref{fig:LastHomeo}. Figures (b)-(e)
    show the construction of the homeomorphism $\f$ from the
    proof of Lemma~\ref{lemma:Fan}\ref{fan7}. \textbf{(b)}: The light-gray area is
    $\square^-$ which needs to be transformed such that any lower
    splitting curve becomes a vertical splitting curve in $G$.  We
    have added an example lower splitting curve $J$ in green color.
    It goes partially in $G$ and partially outside of it, but
    necessarily in $\square^-$ by definition. Note also that by
    definition it must start in $o$ and end in $o^-$. \textbf{(c)}: First move
    the red curve (part of $K^+$) coming from $o$ a bit to clear the
    top area. This of course forces $J$ to move too. We have also
    chosen $\e>$ such that in the $\e$-neighborhood of the vertical
    edges of $G$, $\Fan(K)$ is still trivial and consists of horizontal
    continuations of the lines that come into $G$ from the sides.  At
    this point we mark with darker green those parts of $J$ which will
    never need to move for $\f$ will be identity in ``deep'' inside
    $G$.  \textbf{(d)}: Start the transformation by modifying $\square^-$ from
    the sides while keeping the curves intact (i.e. $\f[K]\subset K$),
    but squeezing $J$.  \textbf{(e)}: In the end $J$ is squeezed entirely into
    $G$ and $\square^-$ is squeezed into a rectangle whose vertical
    sides contain the vertical sides of $G$. Just like our green $J$,
    every lower splitting curve must split $G$ vertically. The
    constructed homeomorphism $\f$ keeps the rectangle
    $[1/4+\e,3/4-\e]\times [0,1/2]$ fixed with possible exceptions
    above $G$ (in the beginning of the construction). In this case the
    number of intersection between $J$ and $K^+$ (red) outside $G$
    before the transformation is 4. All these intersections are moved
    inside $G$ by the end of the transformation.}.
  \label{fig:squareMtoG}
\end{figure}

The following is the main result of this section.

\begin{Prop}[A Bound for Hybrid Crossings for Tangle Diagrams]\label{prop:chi_ge_MINMIN}
  Suppose $T$ is a centered subdivided t-diagram. Then
  $$\intersect(T^-,T^+)\ge \MIN_{\ray}(T^-)+\MIN_{\ray}(T^+).$$
\end{Prop}
\begin{proof}
  Let $\Fan(T)$ be the subdivided tangle diagram given by
  Lemma~\ref{lemma:Fan}. Since $\Fan(T)\sim^*_0 T$ and $\MIN_{\ray}$
  is invariant under $\sim^*_0$ (Lemma~\ref{lemma:invariance_under_R0} applies
  to composite tangle as well as knot diagrams),
  we can assume without loss of generality that $\Fan(T)=T$.
  Let $S^-$ and $S^+$ be as defined in the statement of Lemma~\ref{lemma:Fan}.
  By \ref{fan4}, \ref{fan6}, and \ref{fan8}, $S^-$ and $S^+$
  satisfy the assumptions of Proposition~\ref{prop:Hybrid}, so it follows that
  \begin{align*}
    \chi_h(\Fan(T))&=\intersect(S^+,S^-)&& \text{by \ref{fan5}}\\
                  &=\Card(S^-\cap S^+)&& \text{by the definition of }\intersect(\cdot,\cdot)\\
                  &\ge \MIN_{\vert}(S^+)+k-1.&&\text{by Proposition \ref{prop:Hybrid}}\\
                  &=\MIN_{\vert}(S^+)+\MIN_{\vert}(S^-)&&\text{by the definition of $k$ in \ref{fan5}}\\
                  &=\MIN_{+}(\Fan(T)^-)+\MIN_{-}(\Fan(T)^+)&&\text{by \ref{fan3} and \ref{fan7}}\\
                  &\ge \MIN_{\ray}(T^-)+\MIN_{\ray}(T^+)&&\text{by Lemma \ref{lemma:Raypm_ineq}}.\tag*\qedhere
  \end{align*}
\end{proof}

\begin{Cor}[A Bound for Hybrid Crossings of a Knot Diagram]\label{cor:chi_ge_MINMIN}
  Suppose $K$ is a centered knot diagram. Then
  $\intersect(K^-,K^+)\ge \MIN_{\ray}( K^-)+\MIN_{\ray}(K^+)$.  
\end{Cor}
\begin{proof}
  Let $T=(K,K^+,K^-,\{1\},\es,\ell)$. Then $T$ is a centered subdivided t-diagram.
  The result follows from Proposition~\ref{prop:chi_ge_MINMIN}.
\end{proof}

\section{Proof of Theorem~\protect{\ref{thm:Main}}}
\label{sec:Proof1}
\MainThm* 
\begin{proof}
  \begin{align*}
    \chi(\widehat K)&=\chi(\widehat K^+)+\chi(\widehat K^-)+\intersect(\widehat K^-,\widehat K^+)&&\text{by Remark \ref{rem:mirror}}\\
           &\ge \chi(\widehat K^+)+\chi(\widehat K^-)+\MIN_{\ray}(\widehat K^-)+\MIN_{\ray}(\widehat K^+)&&\text{by Corollary \ref{cor:chi_ge_MINMIN}}\\
           &=\chi(\widehat K^+)+\MIN_{\ray}(\widehat K^+) +\chi(\widehat K^-)+\MIN_{\ray}(\widehat K^-)&&\text{rearranging terms}\\
           &\ge c(K_1)+c(K_0)&&\text{by Corollary~\ref{cor:Bound}}. \tag*\qedhere
  \end{align*}
\end{proof}

\section{Proof of Theorem \protect{\ref{thm:MainLinks}}}
\label{sec:Links}
 
\MainThmLinks*
\begin{proof}
  By using subdivided tangle diagrams instead of knots
  (Definition~\ref{def:SubdividedTangleDiagram}) one straightforwardly
  adapts all the proofs of Section \ref{sec:NonHybrid} to links,
  including Corollary \ref{cor:Bound}.  Then the proof is the same as
  that of Theorem~\ref{thm:Main} presented in
  Section~\ref{sec:Proof1}.  (Use Proposition~\ref{prop:chi_ge_MINMIN}
  instead of Corollary~\ref{cor:chi_ge_MINMIN}.)
\end{proof}

\section{Discussion and Open Questions}
\label{sec:Discussion}

\begin{Question}\label{q:1}
  Is there a less \emph{ad hoc} topological interpretation (than that
  presented in Corollary~\ref{cor:Topo2}) of the fact that it is
  impossible to reduce the number of crossings of a composite knot
  using moves in~$\RR^*_{opp1}$?
\end{Question}

\begin{Question}\label{q:PropBound}
  Find invariants of knot diagrams $?^-(K)$ and $?^+(K)$ such that
  Corollary \ref{cor:Bound} can be proved with $?^{\pm}(K)$ instead of
  $\MIN_{\ray}(K^{\pm})$ and with $\sim_{\RR}$ instead of
  $\sim^X_{\RR}$ but also such that
  Proposition~\ref{prop:chi_ge_MINMIN} still holds when $?^{\pm}(K)$ is
  used instead of $\MIN_{\ray}(K^{\pm})$.  This would resolve the
  addivitiy of the crossing number conjecture by an analogue of
  the proof of Theorem~\ref{thm:Main}. 
\end{Question}

\bibliographystyle{alpha}
\bibliography{ref}

\end{document}